\documentclass[11pt,journal]{IEEEtran}
\usepackage{multirow}
\usepackage[dvips]{graphicx}
\usepackage{amsmath}
\usepackage[psamsfonts]{amssymb}
\usepackage{amsxtra}
\usepackage{threeparttable}

\usepackage{stmaryrd}
\usepackage{amsfonts,amsmath,amssymb}

\usepackage{color}
\usepackage[table]{xcolor}

\usepackage{rotating}

\usepackage{cite}

\usepackage{algorithm}
\usepackage{algorithmic}
\usepackage{booktabs}
\usepackage{url}

\usepackage{colortbl}

\definecolor{shadingcolor}{rgb} {.87,    0.92,    .98}
\newcommand{\shadingbox}[1]{   
    \fboxsep 0pt
    \colorbox{shadingcolor}{
        {\hskip -2pt #1}\hskip -2.5pt
    }
}

\newcommand{\minitab}[2][l]{\begin{tabular}{@{}#1}#2\end{tabular}}

\input{def1.set}

\usepackage{mathpazo}

\usepackage[mathscr]{eucal}

\usepackage[font=small,labelfont=bf,tableposition=top]{caption}
\usepackage[font=footnotesize]{subcaption}

\begin{document}
%

\title{Tensor Networks  for Big Data Analytics and
Large-Scale Optimization Problems}

\author{Andrzej CICHOCKI\\
RIKEN Brain Science Institute, Japan \\
and  Systems Research Institute of the  Polish Academy of Science, Poland\\
{\small Part of this work was presented on the  Second International Conference on Engineering and Computational Mathematics (ECM2013), Hong Kong December 16-18, 2013  (invited talk)}}


\maketitle

\begin{abstract}
Tensor decompositions  and  tensor networks  are emerging and promising tools for  data analysis and data mining.
In this paper we review basic and emerging  models and associated  algorithms for large-scale
tensor networks, especially  Tensor Train (TT)  decompositions  using novel mathematical and graphical representations.
We discus the concept of tensorization (i.e., creating very high-order tensors
from lower-order original
data)  and super compression of data achieved via quantized tensor train  (QTT) networks.
The main objective of this paper is
to show how tensor networks can be used to solve a wide class of big data optimization problems (that are far from tractable by classical numerical methods) by applying tensorization and performing all operations using relatively small size matrices and tensors and applying iteratively optimized and approximative tensor contractions.\\


\end{abstract}

\section{\bf Introduction and Motivations}

\begin{figure}[p]
\includegraphics[width=8.99cm,height=8.3cm]{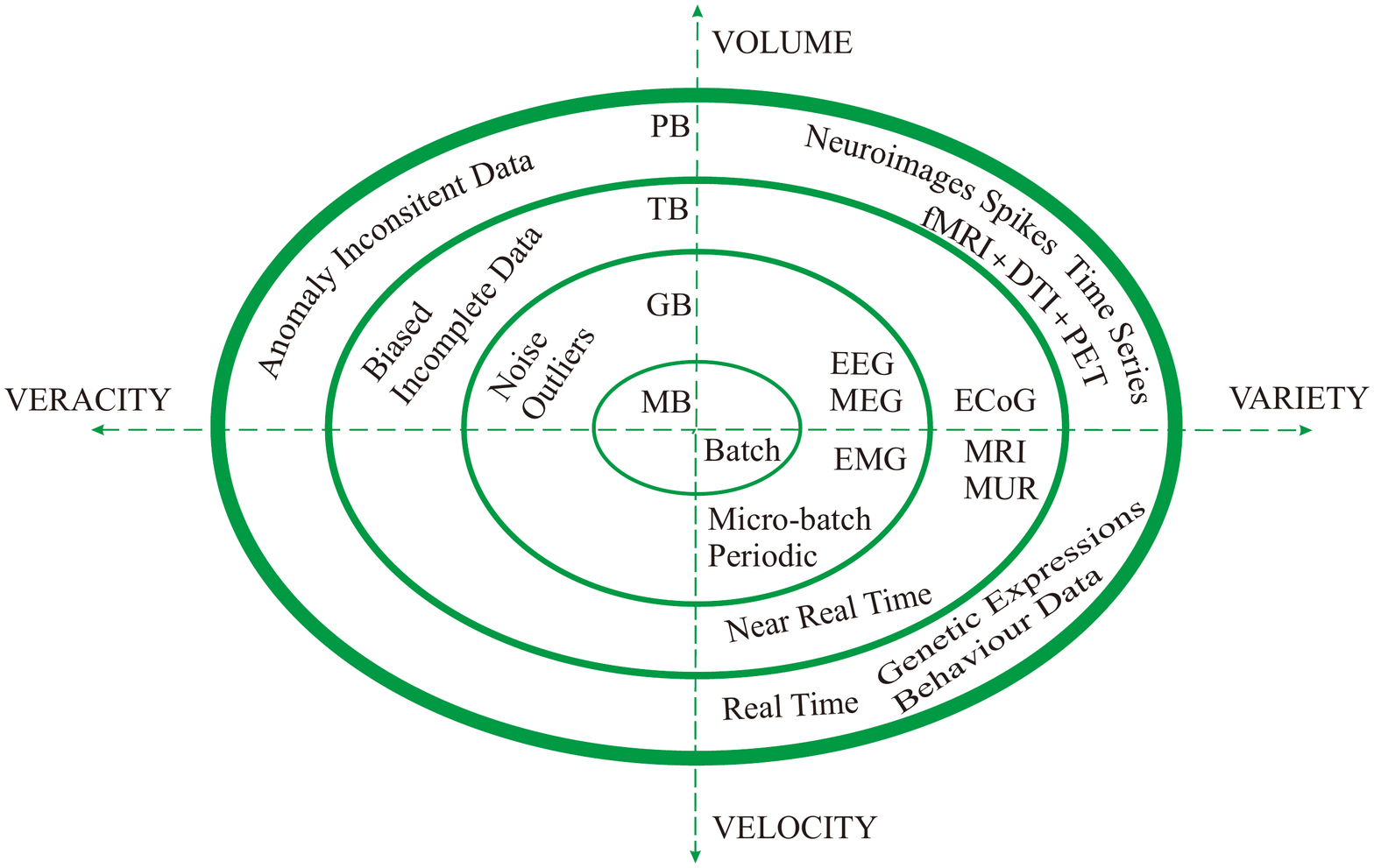}
\includegraphics[width=9.2cm,height=9.8cm]{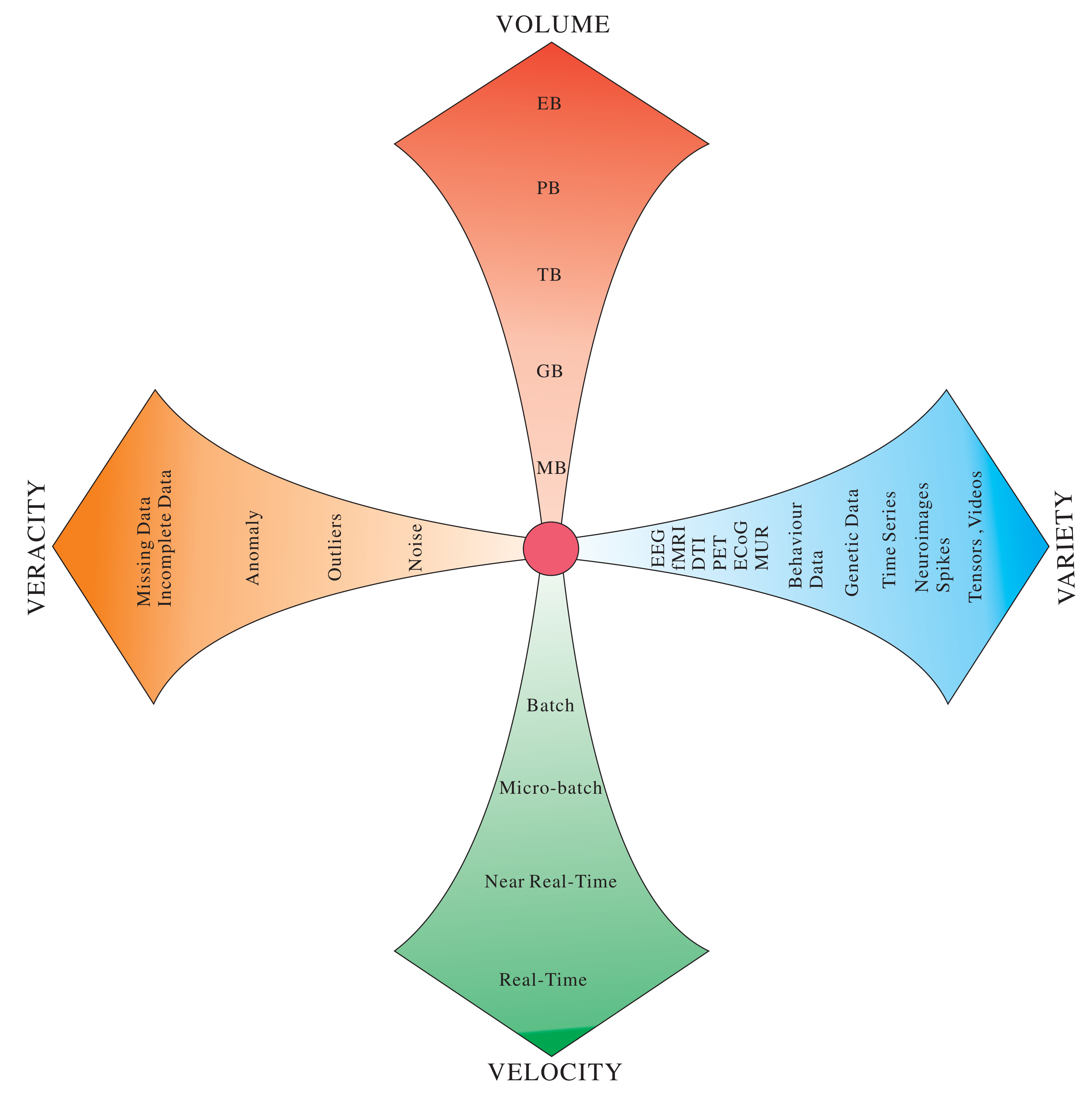}
\caption{Four ``V''s  of big (brain) data: Volume - scale of data, Variety - different forms (types) of  data, Veracity - uncertainty of data, and Velocity - speed at which stream of data is generated and processed. Illustration of challenges for human brain data, which involves analysis of multi-modal, multi-subjects neuroimages, spectrograms, time series, spikes trains, genetic and behavior data. One of the challenges in computational and system neuroscience is to perform fusion or assimilation for  various kinds of data and to understand the relationship and links between them.  Brain data can be recorded by electroencephalography (EEG), electrocorticography (ECoG), magnetoencephalography (MEG), fMRI, DTI, PET, Multi Unit Recording (MUR), to name a few.}
\label{Fig:Big-data}
\end{figure}

Big Data can have a such huge volume and high complexity that  existing standard methods and algorithms become inadequate for the processing and optimization of such data. Big data is characterized not only by big Volume but also  by other specific ``V'' features/challenges: Veracity, Variety, Velocity, Value.  Fig. \ref{Fig:Big-data} illuminates the big data characteristics for brain research related problems. High Volume implies the need for algorithms that are scalable; high Velocity is  related to the processing of stream of data in near real-time; high Veracity  calls for robust and predictive algorithms for noisy, incomplete and/or inconsistent data,  high Variety  require integration across  different types of data, e.g., binary, continuous data, images, time series, etc., and finally  Value refers to   extracting high quality  and consistent data which could lend themselves to  meaningful and interpretable results.

Multidimensional data is becoming ubiquitous across the sciences and engineering because they are increasingly being gathered by information-sensing devices and remote sensing. Big data such as multimedia data (speech, video), and medical/biological data, the analysis  of which critically requires a paradigm shift in order to efficiently process massive datasets within tolerable  time.  Tensors -- multi-dimensional generalizations of matrices, provide often a natural sparse and distributed representation for such data.

 Tensors have been adopted in  diverse branches of data analysis, such as in signal and image processing, 
  Psychometric, Chemometrics, Biometric, Quantum Physics/Information, Quantum Chemistry and Brain Science \cite{Cich-Lath,Cichocki-era,NMF-book,Kolda08,Cichocki-SICE,Hackbush2012,Smilde,Kroonenberg}.
Tensors are particularly attractive for data  which exhibit not only huge volumes but also very high variety, for example, they are suited for problems in
 bio- and neuro-informatics or computational  neuroscience  where data  are collected in various forms of big, sparse tabular,  graphs or networks  with multiple aspects and high dimensionality.

 Tensor decompositions (TDs) provide some extensions of  blind source separation (BSS) and 2-way (matrix) Component Analysis (2-way CA) to  multi-way component analysis (MWCA) methods \cite{Cich-Lath}. Furthermore, TNs/TDs are suitable for dimensionality reduction, they can handle missing values, and noisy data \cite{KressnerSV2013}.
They are also potentially useful
for analysis of linked  (coupled) block of big tensors with millions and even billions of non-zero entries, using the map-reduce paradigm, as well as out-of-core
approaches \cite{Cichocki-era,Wang-out-core05,Suter13,Phan-CP,NLee-Cich14}.
Moreover, multi-block  tensors which arise in numerous important applications (that require the
analysis of  diverse and partially related data) can be decomposed to common (or correlated) and uncorrelated or statistically independent  components.
The effective analysis of coupled tensors requires the development of new models and associated
algorithms and  software that can identify the core relations that may exist
among the different tensors, and scale to extremely large datasets.

\begin{figure} [t]
\centering
\includegraphics[width=3.6cm]{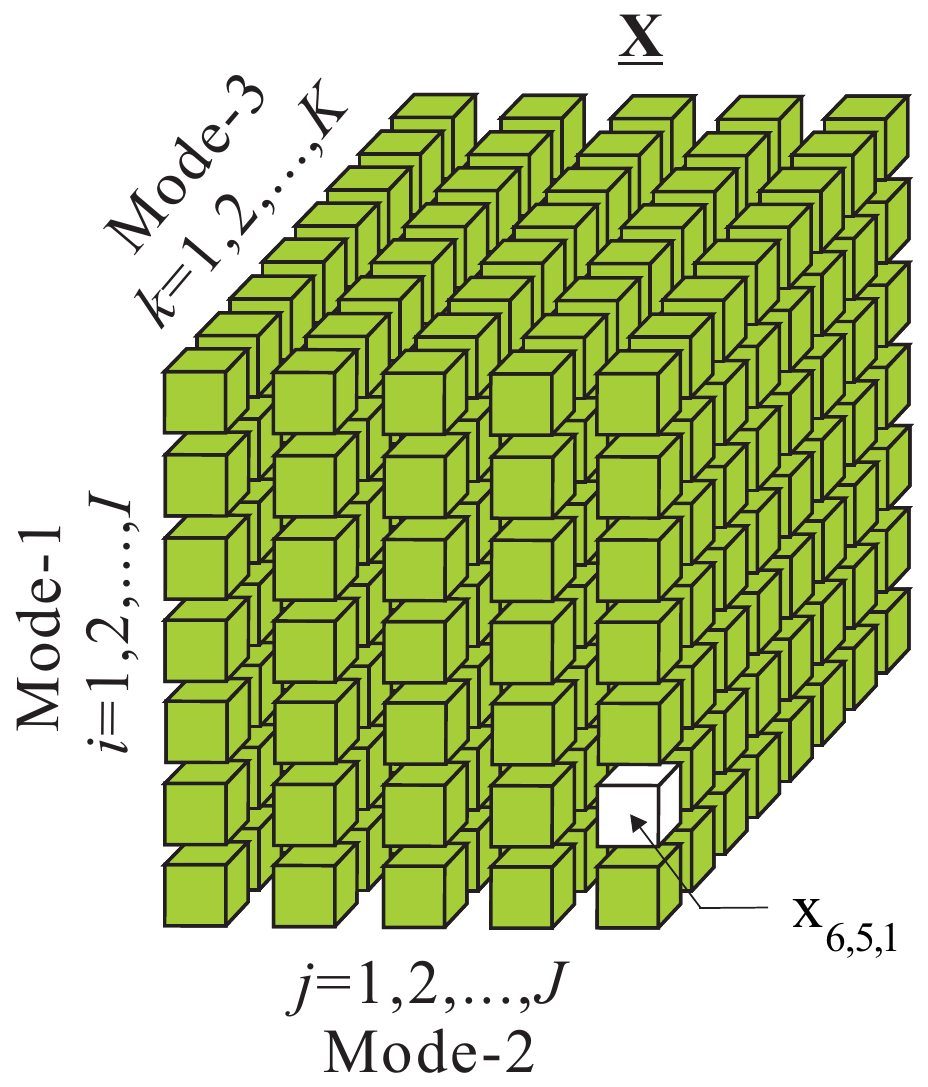}\\
\includegraphics[width=4.6cm]{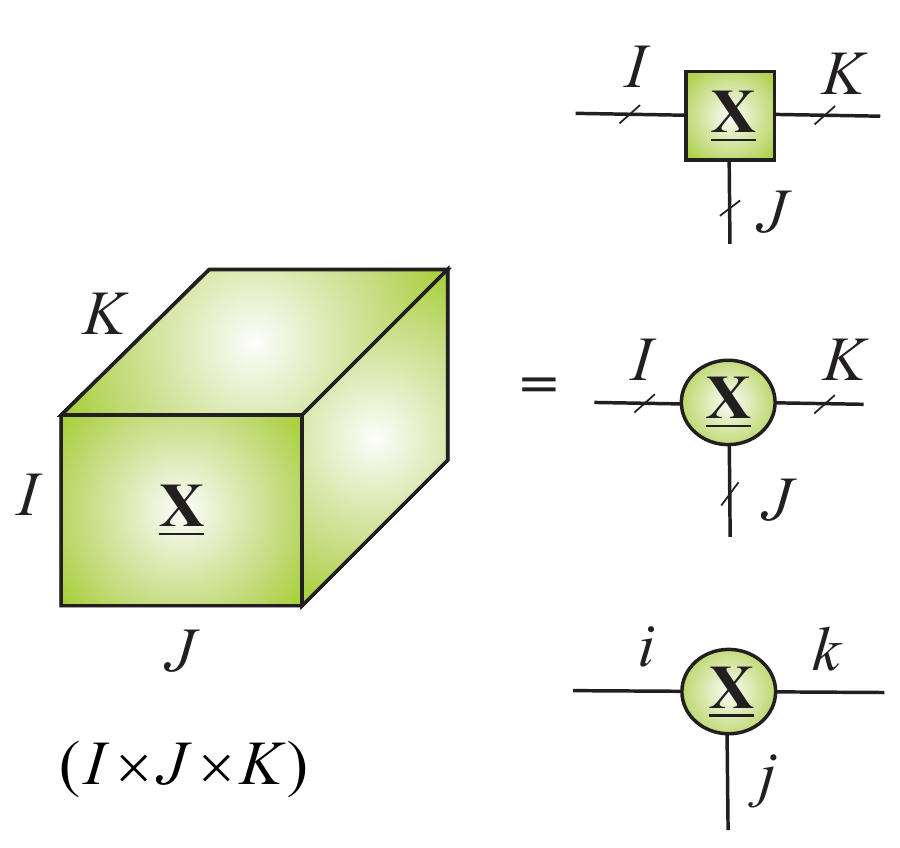}
\caption{A 3rd-order tensor  $\underline \bX \in \Real^{I \times J \times K}$ with
entries $x_{i j k}=\underline \bX(i,j,k)$ and exemplary symbols used in
tensor network diagrams. Each node in the diagram represents a tensor and each edge represents a mode or dimension.  We  indicate maximum size  in each mode by $I,J,K$ or  running indices: $i=1,2,\ldots,I; \;
j=1,2,\ldots,J$ and $k=1,2,\ldots,K$.}
\label{Fig:fibers}
\end{figure}

\begin{figure} [ht]
\centering
\includegraphics[width=8.6cm]{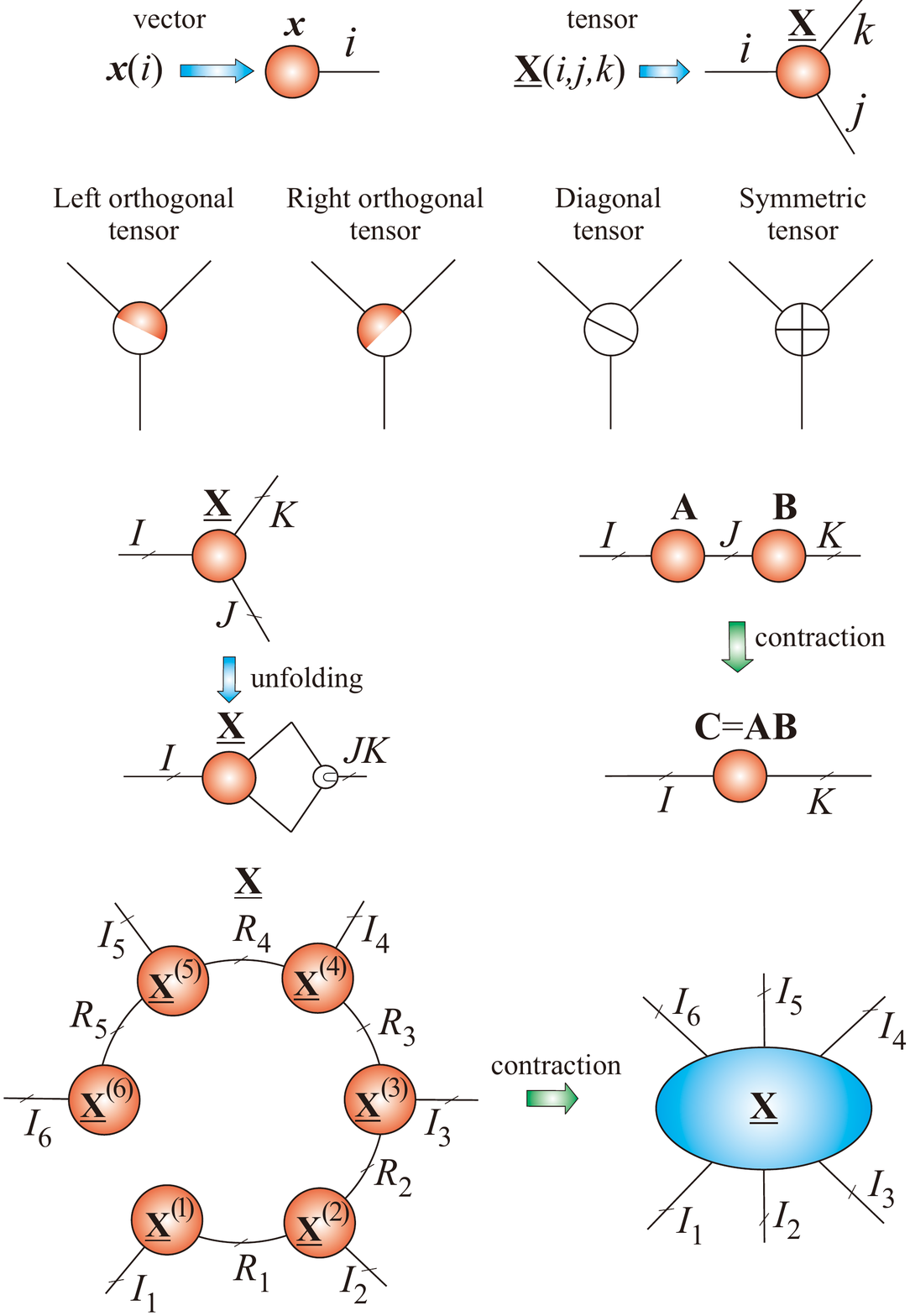}
\caption{Basic symbols and operations for tensor network diagrams. Modes (dimensions) are indicated by running indices ($i,j,k,\ldots$ where $i=1,2,\ldots,I$; \; $j=1,2,\ldots,J$; \; $k=1,2,\ldots,K$,\;$r=1,2,\ldots,R$) in each mode or size of the modes ($I,J,K,R\ldots$). For higher order tensors we will use the symbol $i_n=1,2,\ldots,I_n$ for $n=1,2,\ldots,N$, where $N$ is the order of a tensor.  The minimum set of internal indices $\{R_1,R_2,R_3,\ldots \}$ is called the multilinear rank  od a specific tensor network \cite{Orus11}.}
\label{Fig:symbols}
\end{figure}

Complex interactions and operations  between tensors can be visualized by tensor network diagrams in which  tensors are represented graphically by  nodes or any shapes (e.g., circles,  spheres, triangular, squares, ellipses) and each outgoing edge (line)  emerging from a node represents a mode (a way, a dimension, indices) (see Fig. \ref{Fig:fibers}).  In contrast to  classical graphs, in tensor network diagrams an edge does not need  connect two nodes, but may be connected to only one node. Each such free (dangling) edge corresponds to a (physical) mode that is not contracted and, hence, the order of the entire  tensor network is given by the
number of free  (dangling) edges (see Fig.  \ref{Fig:symbols}).
Tensor  network diagrams are very helpful not only in visualizing tensor
decompositions but also to express complex  mathematical (multilinear) operations of contractions of tensors.
Tensor networks are  connected to   quantum physics, quantum chemistry and
quantum information, which studies the ways to possibly build a quantum computer and to program
it  \cite{Orus11,Orus2013}.

To summarize, the benefits of multiway  (tensor) analysis methods  for big data include:
\begin{itemize}

\item ``Super'' - compression of huge multidimensional data via tensorization and  decompositions of a high-order tensor into factor matrices and/or core tensors of low-rank and low-order;

     \item By performing all mathematical operations in feasible tensor formats \cite{Khoromskij-SC};

\item  Very flexible distributed representations of structurally rich data;


 \item Possibility to operate with noisy and  missing  data  by using powerful low-rank tensor/matrix approximations   and by exploiting robustness and stability of tensor network decomposition algorithms;

\item  A framework to incorporate  various diversities or constraints in different modes or different factors (core tensors) and thus naturally extend the standard (2-way) CA and BSS methods to large-scale multidimensional data;

\item Tensor networks  not only  provide graphically illustrative large distributed networks but also perform complex tensor operations (i.e., tensor contractions and reshaping)  in an intuitive way and  without  using explicitly  mathematical expressions.

\end{itemize}

%

Review and tutorial papers \cite{Cich-Lath,Kolda08,Comon-ALS09,Lu-2011,Morup11,Nikos04}  and books \cite{Smilde,Kroonenberg,NMF-book,Hackbush2012} dealing with TDs  and TNs already exist,
however, they typically focus  on standard models and/or do not  provide explicit links
 to  big data processing topics and/or do not explore  connections to
 wide class of optimization problems.
%
%
This  paper extends beyond the standard tensor decomposition models such as the Tucker and CPD models,
and aims to demonstrate flexibilities  of TNs  in the optimization problems of multi-dimensional, multi-modal data, together with their role as a mathematical backbone for the discovery of hidden structures in  large-scale data \cite{NMF-book,Kolda08}.

 \begin{figure}[ht!]
\centering
5th-order tensor\\
\vspace{0.3cm}
\includegraphics[width=4.2cm]{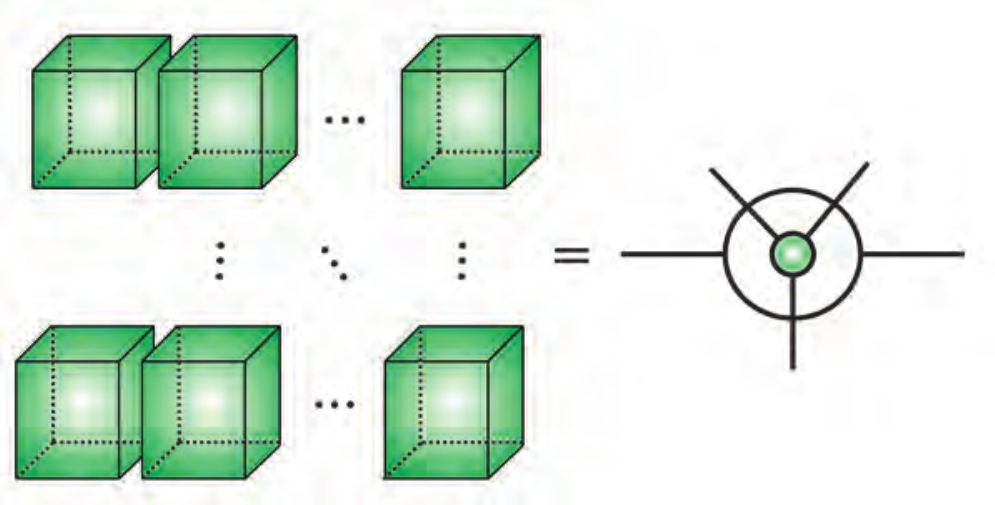}\\
\vspace{0.3cm}
6th-order tensor\\
\vspace{0.3cm}
\includegraphics[width=4.2cm]{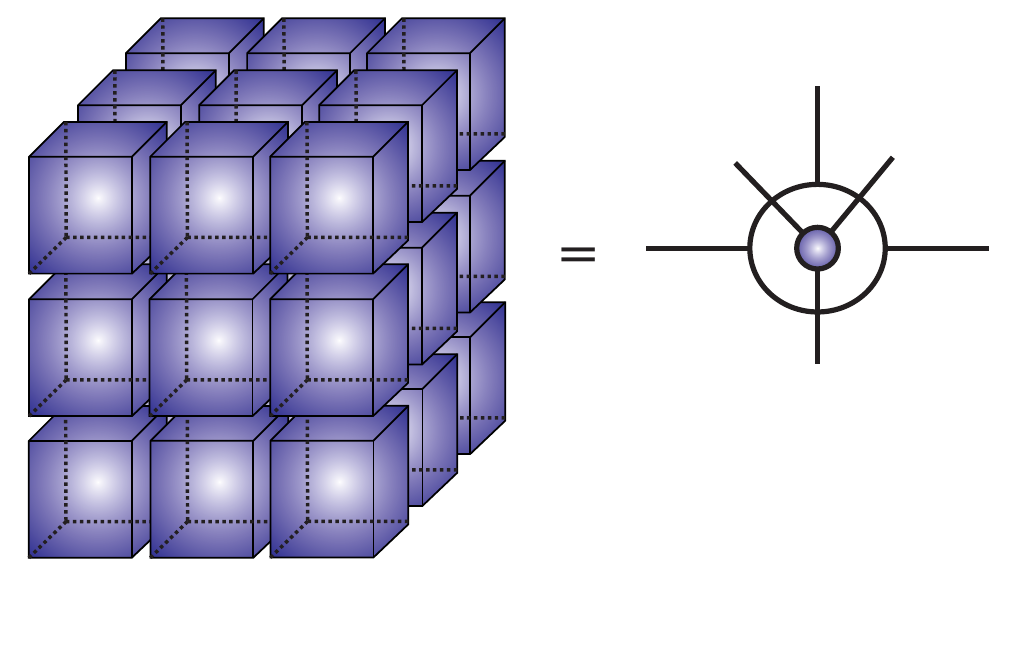}
\caption{Symbols and graphical representations of  higher-order block tensors. Outer circle indicates global structure of a block tensor (e.g., matrix, a 3rd-order tensor), while inner circle indicate the structure  of each element or block of the tensor.}
\label{Fig:symbols2}
\end{figure}

Our objective is to both review tensor models for big data, and to systematically introduce emerging models and associated  algorithms for large-scale  TNs/TDs, together with illustrating the many  potential applications.
 Apart from the optimization framework considered many other challenging problems for big data  related  to  anomaly detection, visualization, clustering, feature extraction and classification   can    also  be solved using  tensor network decompositions and low-rank tensor approximations.

\section{\bf Basic Tensor Operations}

A higher-order tensor can be interpreted as a multiway array of numbers, as illustrated in Figs. \ref{Fig:fibers} and \ref{Fig:symbols}.
 Tensors are denoted  by  bold underlined capital letters, e.g., $\underline \bX  \in \Real^{I_{1} \times I_{2} \times \cdots \times I_{N}}$ (we assume we shall assume that all entries of a tensor are real-valued).
The order  of a tensor is the number of its ``modes'', ``ways'' or ``dimensions'', which include e.g., space, time, frequency,  trials, classes, and dictionaries.
Matrices (2nd-order tensors) are denoted by boldface capital letters, e.g., $\bX$, and vectors (1st-order tensors)  by boldface lowercase letters; for instance  the columns of the matrix $\bA=[\ba_1,\ba_2, \ldots,\ba_R]  \in \Real^{I \times R}$ are denoted by $\ba_r$ and elements of a matrix (scalars) are denoted by lowercase letters, e.g., $a_{ir}$.
Basic tensor and TN notations are given in Table \ref{table_notation1} and illustrated in Figs.
 \ref{Fig:fibers} -- \ref{Fig:symbols2}.
 It should be noted that hierarchical block matrices can be represented by tensors and vice versa. For example, 3rd- and 4th-order tensors  can be represented by block matrices  and all algebraic operations can be equally performed on block matrices \cite{Cichocki-era}.

The most common tensor multiplications are denoted by: $\otimes$ for the Kronecker, $\odot$ for the Khatri-Rao,  $\*$ for the Hadamard (componentwise),  $\circ$ for the outer  and  $ \times_n$ for the mode-$n$  products (see also Table \ref{table_notation1}).
General basic operations, e.g., $vec(\cdot)$, $\diag\{\cdot\}$, are defined as in  MATLAB.
 We refer to \cite{Cichocki-era,NMF-book,Kolda08} for more detail regarding the basic
  notations and tensor operations.

 Subtensors are formed when a subset of indices is fixed. Of particular interest are  {\it fibers} (vectors), defined  by fixing every index but one,  and
{\it slices} which are two-dimensional sections (matrices) of a tensor,
obtained by fixing all the indices but two. 
A matrix has two modes: rows and columns, while an $N$th-order tensor has $N$ modes.

\begin{figure}[h!]
\includegraphics[width=8.6cm]{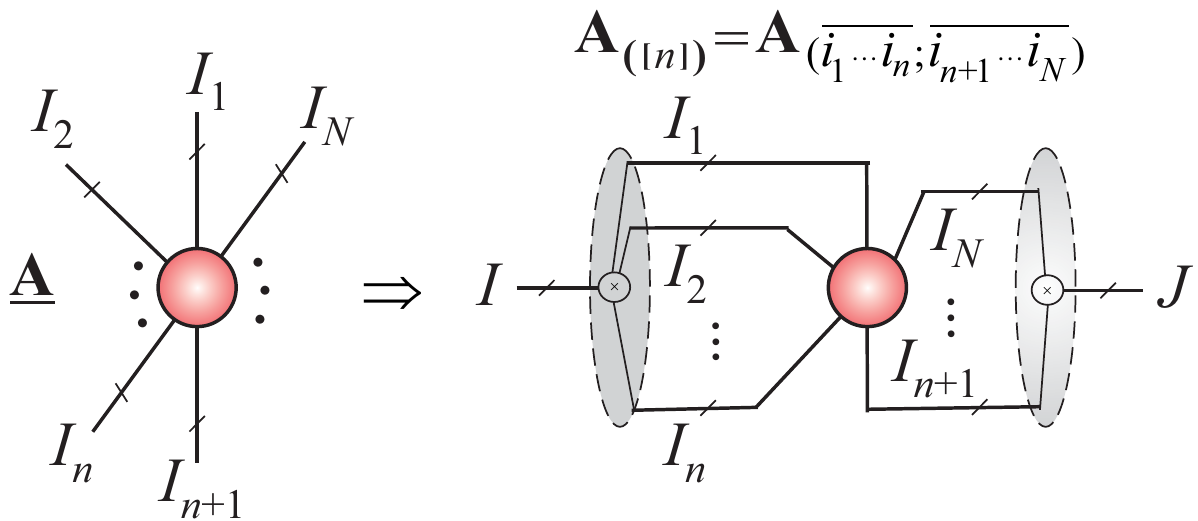}
\caption{Unfolding  the $N$th-order tensor $\underline \bA \in \Real^{I_1 \times I_2 \times \cdots \times I_N}$ \r into a matrix $\bA_{([n])} = \bA_{(\overline{i_1 \cdots i_n}\; ; \; \overline{i_{n+1} \cdots i_N})} \in \Real^{I_1 I_2 \cdots I_n \times I_{n+1} \cdots I_N}$. All entries of an unfolded tensor  are arranged in a specific order. In a more general case, let $\brr=\{m_1,m_2,\ldots,m_R\} \subset \{1,2,\ldots, N\}$ be the row indices and
$\bc=\{n_1,n_2,\ldots, n_C \} \subset \{1,2,\ldots, N\} -\brr$ be the column indices, then the mode-$(\brr,\bc)$ unfolding of $\underline \bA$ is denoted as
$\bA_{(\brr,\bc)} \in \Real^{I_{m_1} I_{m_2} \cdots I_{m_R}  \times I_{n_1} I_{n_2} \cdots I_{n_C}}$.}
\label{Fig:unfolding}
\end{figure}

The process of unfolding (see Fig. \ref{Fig:unfolding}) flattens a tensor into a  matrix \cite{Kolda08}.
In the simplest scenario, mode-$n$ unfolding (matricization,  flattening) of the tensor  $\underline \bA \in \Real^{I_{1} \times I_{2} \times \cdots \times I_{N}}$ yields
a matrix $\bA_{(n)} \in \Real^{I_{n} \times (I_{1} \cdots I_{n-1} I_{n+1} \cdots I_N)}$, with entries $a_{i_n, i_2, \ldots, i_{n-1}, i_{n+1}, \ldots, i_n)}$ such that grouped indices $(i_1,\ldots,i_{n-1},i_{n+1},\ldots,i_N)$ are arranged in a specific order, (in this paper rows and columns are ordered colexicographically). In tensor networks we use, typically a generalized mode-$([n])$ unfolding as illustrated in Fig. \ref{Fig:unfolding}.

\minrowclearance 2ex
\begin{table}[p!]
\caption{Basic tensor notation and matrix/tensor products.} \centering
 { \shadingbox{
    \begin{tabular*}{1.02\linewidth}[t]{@{\extracolsep{\fill}}ll} \hline
{$\underline \bX \in \Real^{I_1 \times I_2 \times \cdots \times I_N}$} & \minitab[p{.4\linewidth}]{$N$th-order tensor of size $I_1 \times I_2 \times \cdots \times I_N$} \\
{$\underline \bG^{(n)},  \; \underline \bX^{(n)}, \; \underline \bS$}   & {core tensors} \\  [-2ex]
{$\underline {\mbi {\Lambda}} \in \Real^{R \times R \times \cdots \times R}$}   & \minitab[p{.42\linewidth}]{\hspace{-0.2em}diagonal core tensor with nonzero $\lambda_r$ entries on diagonal} \\  [-2ex]
{$\bA = [\ba_1,\ba_2,\ldots,\ba_R] \in \Real^{I \times R} $}& \minitab[p{.4\linewidth}]{matrix with  column vectors $\ba_r \in \Real^I$ and entries $a_{ir}$} \\
{$\bA, \, \bB, \,\bC, \; \bB^{(n)}, \; \bU_n $}& {component matrices} \\
$\bX_{(n)} \in \Real^{I_n \times I_1 \cdots I_{n-1} I_{n+1} \cdots I_N}$ & \minitab[p{.65\linewidth}]{mode-$n$ unfolding of $\underline \bX$}\\ [-1ex]
%
$\vtr{\underline \bA}$& vectorization of $\underline \bA$  \\
\\[-2em]
$\underline \bC = \underline \bA \times_n \bB$ & \minitab[p{.5\linewidth}]{\\[-3.2ex] mode-$n$ product of $\underline \bA \in \Real^{I_1 \times I_2 \times \cdots \times I_N}$ and $\bB \in \Real^{J_n \times I_n}$ yields $\underline \bC \in  \Real^{I_1 \cdots \times I_{n-1} \times J_n \times I_{n+1} \cdots \times I_N}$ with entries $c_{i_1 \cdots i_{n-1} \, j \, i_{n+1} \cdots i_N} = \sum_{i_n=1}^{I_n} a_{i_1 \cdots i_n \cdots i_N} b_{j \, i_n}$, and $\bC_{(n)} = \bB \, \bA_{(n)}$} \\
\multicolumn{2}{@{\hspace{-.1ex}}l}{$\underline \bC = \llbracket \underline \bA; \bB^{(1)},  \ldots, \bB^{(N)}\rrbracket =
 {\underline \bA} \times_1 \bB^{(1)} \times_2 \bB^{(2)}  \cdots \times_N \bB^{(N)}$ }\\
%
$\underline \bC = \underline \bA \circ \underline \bB$ & \minitab[p{.51\linewidth}]{\\[-3.2ex] tensor or outer product of $\underline \bA \in \Real^{I_1 \times I_2  \times \cdots \times I_N}$ and $\underline \bB \in \Real^{J_1 \times J_2 \times \cdots \times J_M}$ yields $(N+M)$th-order tensor $\underline \bC$ with entries $c_{i_1 \cdots i_N \,j_1 \cdots j_M} = a_{i_1 \cdots i_N} b_{j_1 \cdots j_M}$} \\
%
%
{$\underline \bX = \ba \circ \bb \circ \bc \in \Real^{I \times J \times K}$}&\minitab[p{.48\linewidth}]{tensor or outer product of vectors  forms a rank-1 tensor with entries $x_{ijk} = a_i b_j c_k$} \\
%
{$\underline \bC = \underline \bA \otimes \underline \bB$}&\minitab[p{.53\linewidth}]{Kronecker product of $\underline \bA \in \Real^{I_1 \times I_2  \times \cdots \times I_N}$ and $\underline \bB \in \Real^{J_1 \times J_2 \times \cdots \times J_N}$ yields $ \underline \bC \in \Real^{I_1 J_1  \times \cdots \times I_N J_N}$ with entries $c_{\overline{i_1,j_1};\ldots;\overline{i_N, j_N}} = a_{i_1, \ldots,i_N} \: b_{j_1, \ldots,j_N}$, where $\overline{i_n, j_n} = j_n+(i_n-1)J_n$} \\[-2ex]
{$\bC = \bA\odot \bB$}&\minitab[p{.45 \linewidth}]{Khatri-Rao product of $\bA \in \Real^{I\times J}$ and $\bB \in \Real^{K\times J}$ yield $\bC \in \Real^{I K \times J}$, with columns $\bc_j = \ba_j \otimes \bb_j$}\\[-1ex]
\\ 
\hline
    \end{tabular*}
    }}
\label{table_notation1}
\end{table}
\minrowclearance 0ex

By a multi-index $i = \overline{i_1,i_2, \ldots, i_N}$ we denote an index which takes all possible
combination of values of $i_1, i_2, \ldots, i_n$,  for $i_n = 1,2,\ldots, I_n$, in a specific order.

{\bf Remark.} The entries of tensors in matricized and/or vectorized form can be ordered  in different forms.
In fact, the multi--index can be defined  using two  different conventions \cite{Dolgov2013alternating2}:

1) Little–-endian convention
\be
\overline{i_1,i_2,\ldots, i_N } &=& i_1 + (i_{2} - 1) I_1+ (i_{3} - 1) I_1 I_2 \notag \\
&\cdots&  + (i_N - 1) I_1 \cdots  I_{N-1},
\ee

2) Big–-endian
\be
&&\overline{i_1,i_2,\ldots, i_N } = i_N + (i_{N-1} - 1)I_N + \notag \\
&&+ (i_{N-2} - 1)I_N I_{N-1}+ \cdots  + (i_1 - 1)I_2 \cdots  I_N. \notag \\
\ee

The little--endian notation
is consistent with the Fortran style of indexing,
while the big--endian notation is similar to numbers written in the positional system and corresponds to reverse lexicographic order \cite{Dolgov2013alternating,Dolgov2013alternating2}.
The definition of unfolding  and  the Kronecker (tensor) product $\otimes$ should be also consistent with the chosen convention{\footnote{The  standard and more popular
definition in multilinear algebra assumes the big--endian convention, which  corresponds  to colexicographic order, while for the development of the efficient program
code, usually,   the little--endian  convention seems to be more convenient
(See more detail the paper of Dolgov and Savostyanov \cite{Dolgov2013alternating2}).}}.
In this paper  we will use the big-endian notation, however to follow this work  it is sufficient to remember that $\bc = \ba \otimes  \bb$ means that $c_{\overline{ij}} = a_i b_j$.

The Kronecker product of two tensors $\underline \bA \in \Real^{I_1 \times I_2 \times \cdots \times I_N}$ and $ \underline  \bB \in \Real^{J_1 \times J_2 \times \cdots \times J_N}$ yields $\underline  \bC = \underline \bA \otimes \underline \bB \in \Real^{I_1 J_1 \times \cdots \times I_N J_N}$ with entries $c_{\overline{i_1, j_1},\ldots, \overline{i_N, j_N}} = a_{i_1, \ldots,i_N} \: b_{j_1, \ldots,j_N}$, where $\overline{i_n, j_n} = j_n+(i_n-1)J_n$ \cite{Phan2012-Kron}.

The outer or tensor product $\underline \bC =\underline \bA \circ \underline \bB$ of the tensors $\underline \bA \in \Real^{I_1 \times \cdots \times I_N}$ and $\underline \bB \in \Real^{J_1 \times \cdots \times J_M}$ is the tensor $\underline \bC \in \Real^{I_1 \times \cdots \times I_N \times J_1 \times \cdots \times J_M }$ with entries $c_{i_1,\ldots,i_N,j_1,\ldots,j_M} = a_{i_1,\ldots,i_N} \; b_{j_1,\ldots,j_M}$.
Specifically, the outer product of two nonzero vectors
$\ba \in \Real^I, \; \bb \in \Real^J$ produces a rank-1 matrix
$\bX=\ba \circ \bb = \ba \bb^T \in \Real^{I \times J}$
and the outer product of  three nonzero vectors: $\ba \in \Real^I, \; \bb \in \Real^J$ and  $\bc \in \Real^K$ produces a 3rd-order rank-1  tensor:
$\underline \bX= \ba \circ \bb \circ \bc \in \Real^{I \times J  \times K}$,
whose entries are $x_{ijk} =a_i \; b_j \; c_k$.
A tensor $\underline \bX \in \Real^{I_1 \times I_2 \times \cdots \times I_N}$ is said to be rank-1 if  it can be expressed exactly as $\underline \bX = \bb^{(1)} \circ \bb^{(2)} \circ \cdots  \circ \bb^{(N)}$ with entries $x_{i_1,i_2,\ldots,i_N} = b^{(1)}_{i_1} b^{(2)}_{i_2} \cdots b^{(N)}_{i_N}$, where $\bb^{(n)} \in \Real^{I_n}$ are nonzero vectors.


The  mode-$n$  product of the tensor
$\underline \bA \in \Real^{I_{1} \times \cdots \times I_{N}}$ and vector  $\bb \in \Real^{I_n}$ is  defined as a tensor
$ \underline \bC = \underline \bA \bar \times_n \bb \in \Real^{I_1 \times \cdots \times I_{n-1}  \times I_{n+1} \times \cdots \times I_N}$, with entries $ c_{i_1,\ldots,i_{n-1},i_{n+1},\ldots, i_{N}} =\sum_{i_n=1}^{I_n} (a_{i_1,i_2,\ldots,i_N})(b_{i_n})$, while a mode-$n$  product of the tensor
$\underline \bA \in \Real^{I_{1} \times \cdots \times I_{N}}$ and a matrix  $\bB \in \Real^{J \times I_n}$
is the tensor
$ \underline \bC = \underline \bA \times_n \bB \in \Real^{I_1 \times \cdots \times I_{n-1} \times J \times I_{n+1} \times \cdots \times I_N}$
with entries $ c_{i_1,i_2,\ldots,i_{n-1},j,i_{n+1},\ldots, i_{N}} =\sum_{i_n=1}^{I_n} a_{i_1,i_2,\ldots,i_N} \; b_{j,i_n}$. This can be also expressed in a matrix form as $\bC_{(n)} =\bB \bA_{(n)}$.

\begin{figure}[ht!]
(a) \hspace{3.9cm} (b)\\
\includegraphics[width=8.6cm,height=4.6cm]{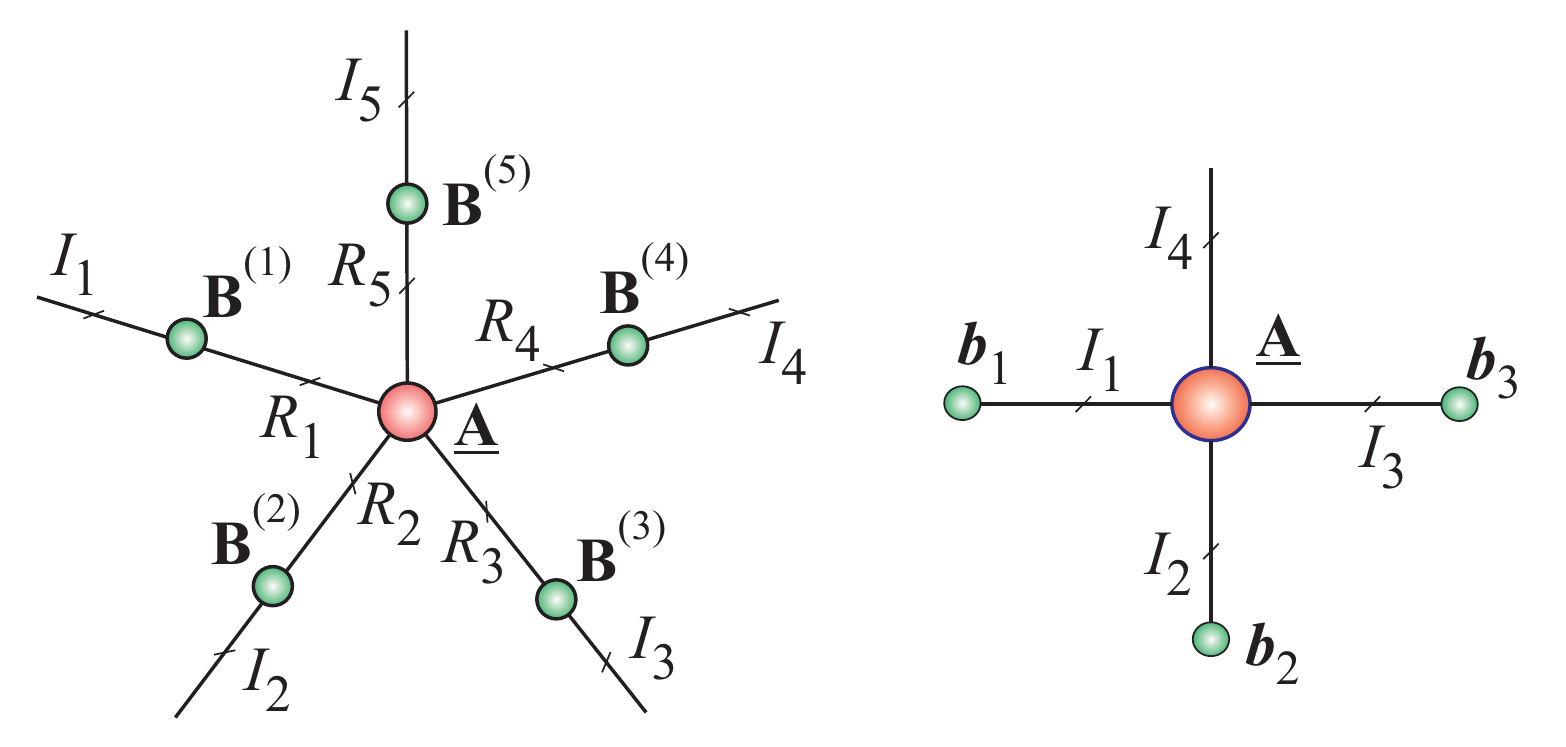}
\caption{(a) Multilinear product of tensor $\underline \bA \in \Real^{R_1 \times R_2 \times \cdots \times R_5}$ and 5 factor (component) matrices $\bB^{(n)} \in \Real^{I_n \times R_n}$ ($n=1,2,\ldots,5$) yields $\underline \bC = \underline \bA \times_1 \bB^{(1)}  \times_2 \bB^{(2)} \times_3\bB^{(3)}  \times_4 \bB^{(4)}   \times_5 \bB^{(5)} \in \Real^{I_1 \times I_2 \times \cdots \times I_5}$ and (b) Multilinear product of tensor $\underline \bA \in \Real^{I_1 \times I_2 \times I_3 \times I_4}$  and  vectors $\bb_n \in \Real^{I_n}$ $(n=1,2,3)$ yields  a vector $\bc= \underline \bA \bar \times_1 \bb_1  \bar \times_2 \bb_2 \bar \times_3  \bb_3 \in \Real^{I_4} $.}
\label{Fig:TO1}
\end{figure}

A full multilinear product
 of a tensor and a set of matrices takes  into account all the modes, and can be compactly written as (see Fig \ref{Fig:TO1} (a)):
\be
\underline \bC &=& \underline \bA \times_1 \bB^{(1)}  \times_2 \bB^{(2)}  \cdots  \times_N \bB^{(N)} \notag \\
&=& \llbracket \underline \bA; \bB^{(1)}, \bB^{(2)}, \ldots, \bB^{(N)} \rrbracket.
\label{mprod}
\ee

\begin{figure}[ht]
\centering
\includegraphics[width=8.6cm]{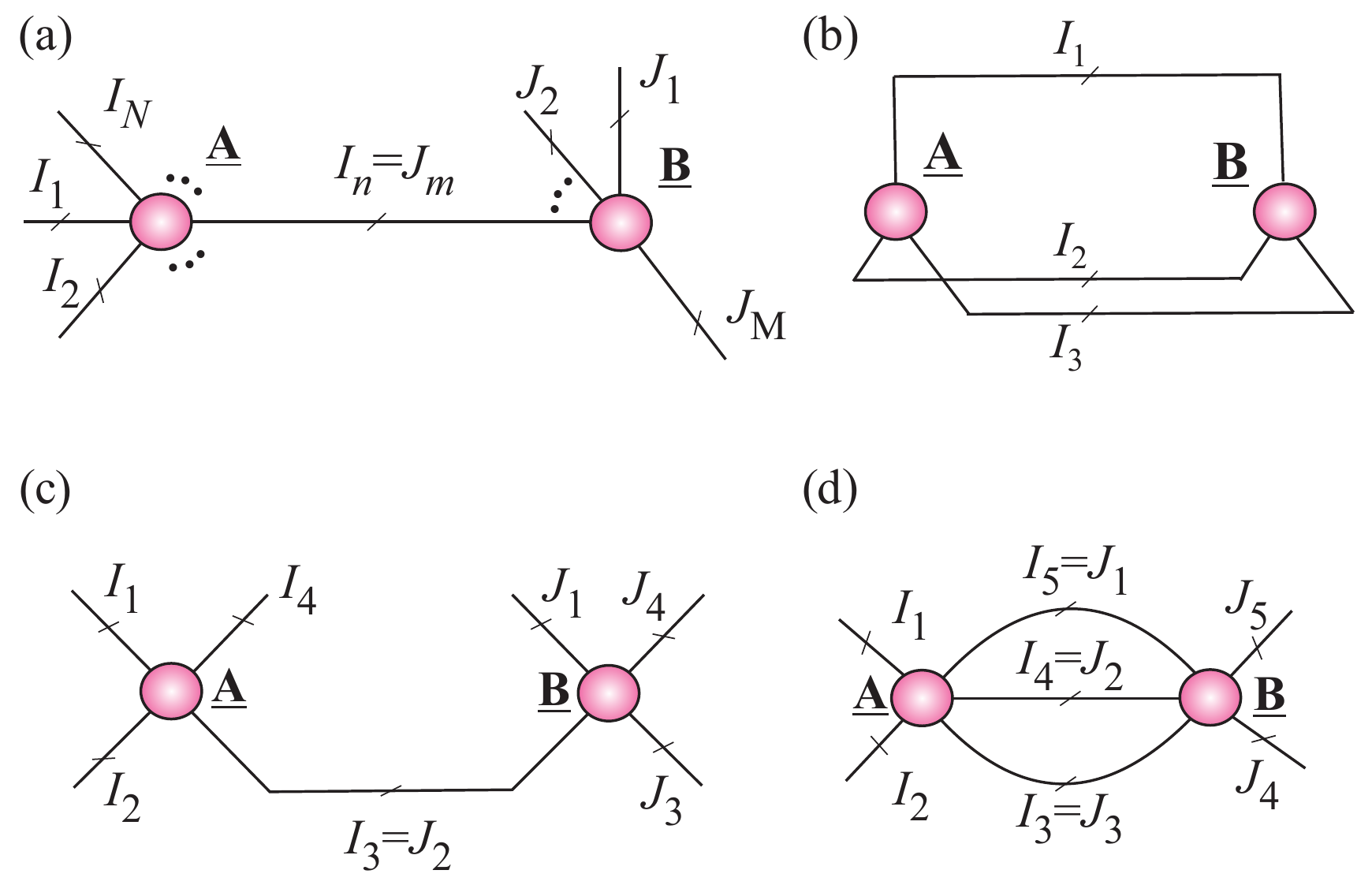}
\caption{Examples of contraction of two tensors: (a) Multilinear product of two tensors is denoted by  $\underline \bC = \underline \bA \; \times_n^m \; \underline \bB$. (b) Inner product of two 3rd-order tensors yields $c=\langle\underline \bA, \underline \bB\rangle=\underline \bA \; \times_{1,2,3}^{1,2,3} \; \underline \bB  =\underline \bA \; \times\; \underline \bB=\sum_{i_1,i_2,i_3} \; a_{i_1,i_2,i_3} \; b_{i_1,i_2,i_3}$.  (c) Tensor contraction of two 4th-order tensors yields the 6h-order tensor $\underline \bC = \underline \bA \; \times_{3}^{2} \; \underline \bB \in \Real^{I_1 \times I_2 \times I_4 \times J_1 \times J_3 \times J_4}$, with entries $c_{i_1,i_2,i_4,j_1,j_3,j_4}= \sum_{i_3} \; a_{i_1,i_2,i_3,i_4}  \; b_{j_1,i_3,j_3,j_4}$. (d) Tensor contraction of two 5th-order tensors yields the 4th-order tensor $\underline \bC = \underline \bA \; \times_{3, 4, 5}^{1, 2, 3} \; \underline \bB \in \Real^{I_1 \times I_2 \times J_4 \times J_5}$, with entries $c_{i_1,i_2,j_4,j_5}= \sum_{i_3,i_4,i_5}\; a_{i_1,i_2,i_3,i_4,i_5}  \; b_{i_5,i_4,i_3,j_4,j_5}$.}
\label{Fig:TO2}
\end{figure}

In a similar way, we can define  the mode-$(^m_n)$ product of two tensors $\underline \bA \in \Real^{I_1 \times I_2 \times \ldots \times I_N}$ and $\underline \bB \in \Real^{J_1 \times J_2 \times \cdots \times J_M}$  with  common modes $I_n=J_m$ that produces  a $(N+M-2)$-order tensor $ \underline \bC  \in \Real^{I_1 \times \cdots I_{n-1} \times I_{n+1} \cdots \times I_N \times J_{1} \times \cdots J_{m-1} \times J_{m+1} \cdots \times J_M}$:
 \be \underline \bC =  \underline \bA \; \times_n^m \; \underline \bB,
 \label{Eq:Tcontrac}
 \ee
  with entries $c_{i_1 \cdots i_{n-1} \, i_{n+1} \cdots i_N, \, j_1, \cdots j_{m-1} \, j_{m+1} \cdots j_M} = \sum_{i=1}^{I_n} a_{i_1 \cdots i_{n-1} \; i \; i_{n+1} \cdots i_N} b_{j_1 \cdots j_{m-1} \; i \; j_{m+1} \cdots j_M}$ (see Fig. \ref{Fig:TO2}) (a).

 When not confusing, the super-index $m$ can be neglected. For example, the mode-1 product of the tensors $\underline \bA \in  \Real^{I_{1} \times I_{2} \times \cdots \times I_{N}}$ and $\underline \bB \in  \Real^{J_{1} \times J_{2} \times \cdots \times J_{M}}$,  with a common first mode $I_1=J_1$ can be written   as
\be
\underline \bC = \underline \bA \; \times_1^1 \; \underline \bB = \underline \bA \times_1 \underline \bB  \in  \Real^{I_{2} \times \cdots \times I_{N} \times  J_2 \times \cdots \times J_{M}},
\ee
with  entries   $c_{\bi_{2:N},\bj_{2:M}}= \sum_{i=1}^{I_1} a_{i, \bi_{2:N}} b_{i,\bj_{2:M}}$, when using  MATLAB notation, $\bi_{p:q}=\{i_p,i_{p+1},\ldots,i_{q-1},i_q\}$.
This operation can be considered as a tensor contraction of two modes.
Tensors can be contracted in several modes or even in all modes
(see Fig. \ref{Fig:TO2}).

Tensor contraction  is a fundamental operation, which  can be considered as a higher dimensional analogue of
inner product, outer product and matrix multiplications, and comprises  computationally
dominant operations in most numerical algorithms.
However, unlike the  matrix by matrix multiplications for which  many
efficient distributed-memory parallel schemes have
been developed, 
for a  tensor contraction we have a rather limited number
of available optimized algorithms \cite{Pfeifer13,Pfeifer2014ncon,ContrationR13}.
In practice, we usually implement approximate tensors contractions
with reduced ranks \cite{ContracPEPS14}.
A significant help in developing effective distributed
tensor contraction algorithms is that the tensors used
in computational models often
exhibit symmetry over all or multiple modes; exploitation
of the symmetry is essential, both in order to save on storage as
well as to avoid unnecessary arithmetic operations \cite{Pfeifer2014ncon,ContrationR13}.

\begin{figure}[t]
(a) \\
\includegraphics[width=8.6cm]{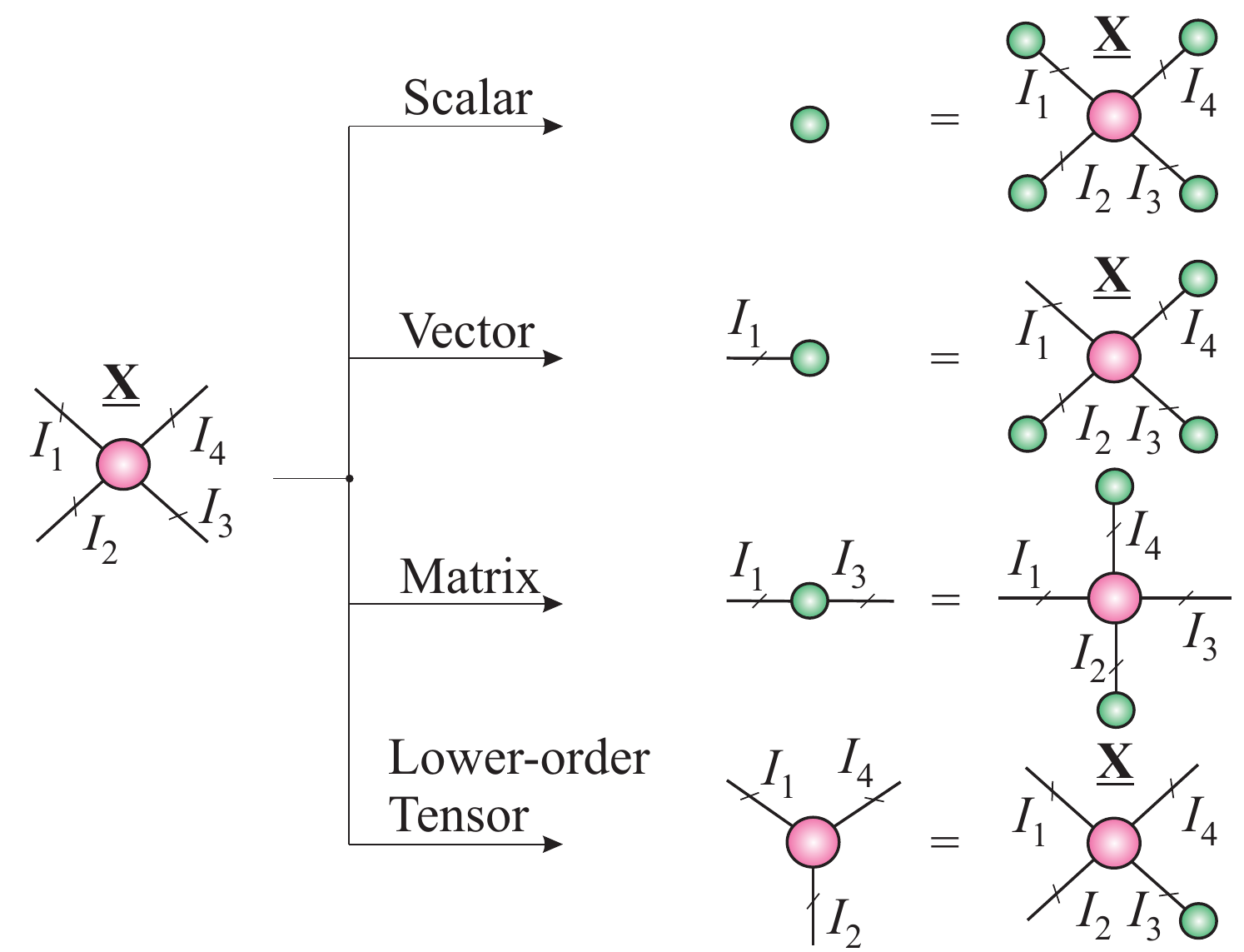}\\
(b)\\
\includegraphics[width=8.6cm,height=6.7cm]{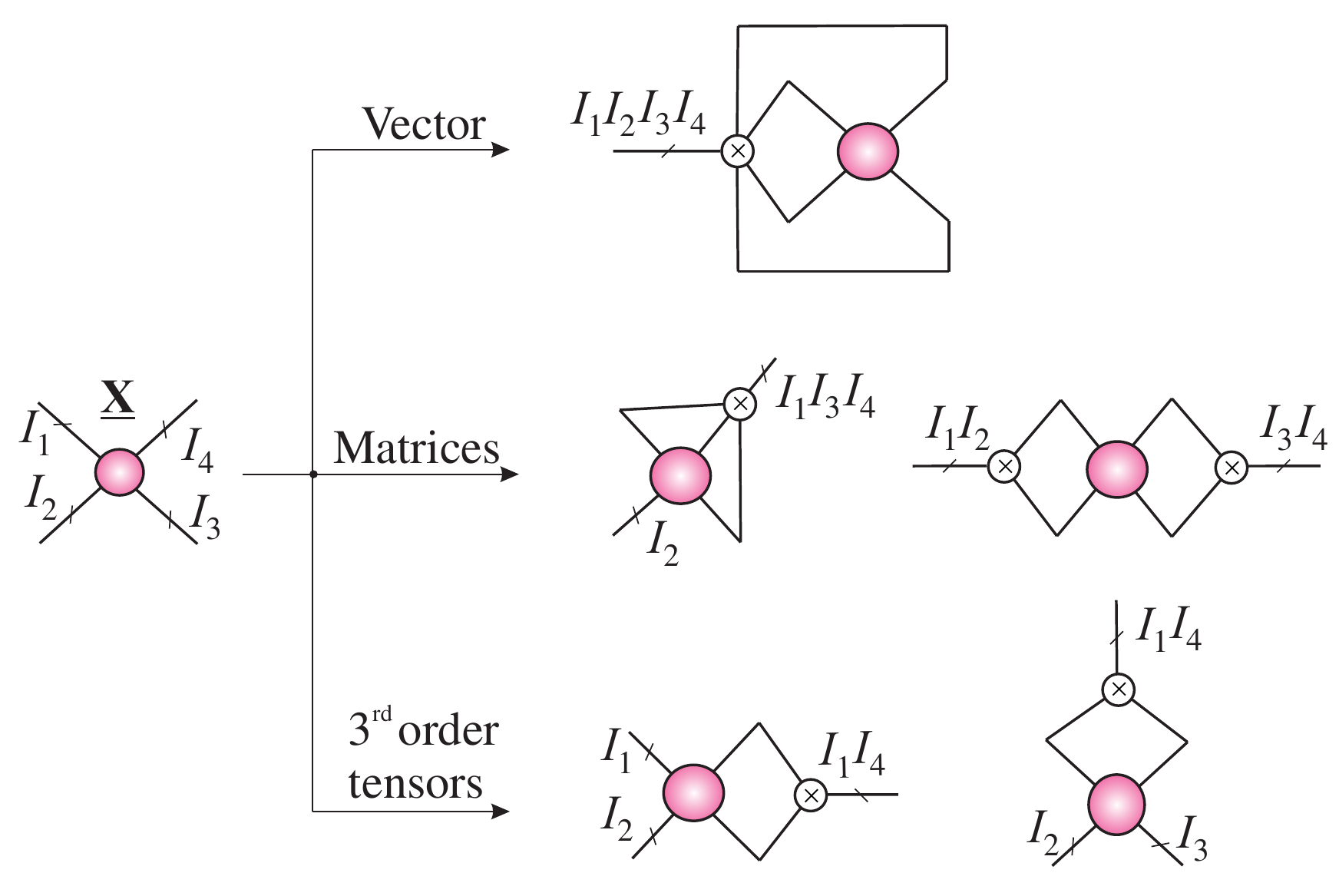}
\caption{(a) Transforming or compressing a 4th-order tensor into scalar, vector, matrix and 3rd-order tensor by multilinear product of the tensor and vectors. (b) Reshaping of a tensor by its vectorization, unfolding and reducing the order by merging the modes.}
\label{Fig:reshaping}
\end{figure}
 \begin{figure}[ht]
\centering
 \includegraphics[width=5.99cm]{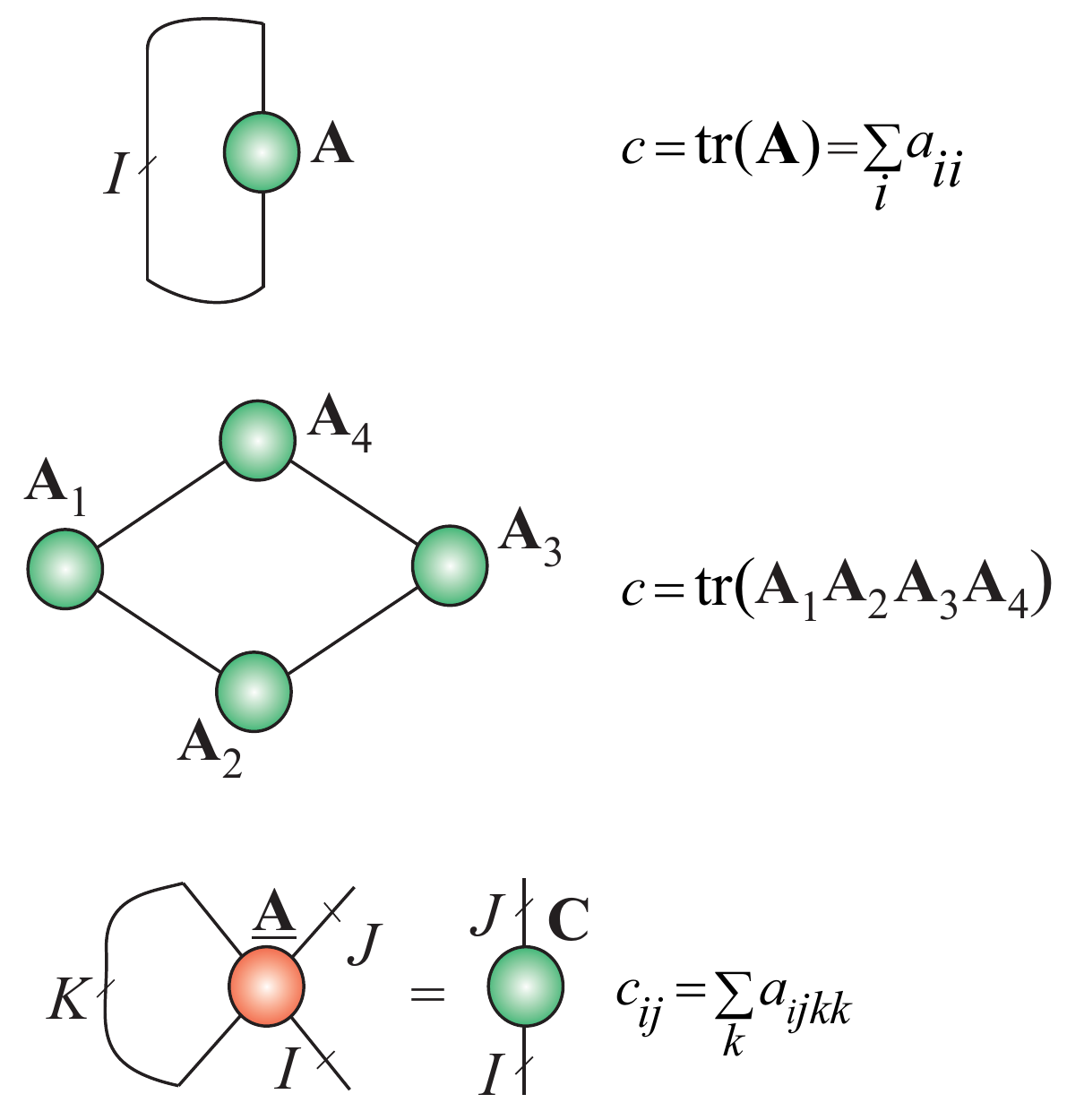}\\
 \vspace{0.6cm}
  \includegraphics[width=8.8cm]{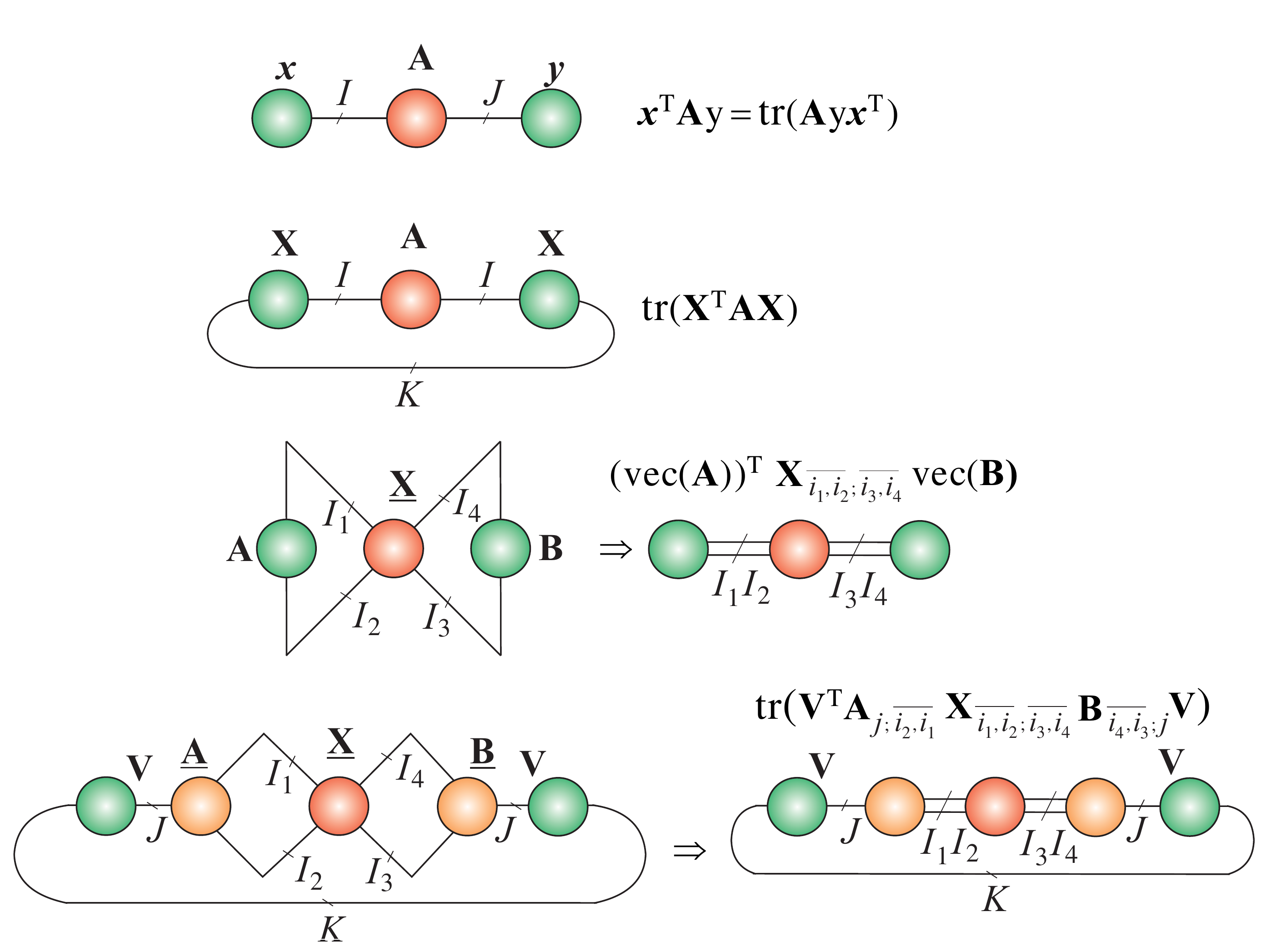}
\caption{Tensor network notation and operations for  traces of matrices and higher-order tensors.}
\label{Fig:trace}
\end{figure}

Tensors often need to  be converted to  traces, scalars, vectors, matrices or  tensors with  reshaped modes and  reduced orders, as illustrated in   Fig. \ref{Fig:reshaping} and Fig. \ref{Fig:trace}.

\section{\bf Low-Rank Tensor Approximations via Tensor Networks}

\subsection{\bf Basic Tensor Network Models}

Tensor networks can be considered as  a new ``language'' for big data tensor decompositions in simulation of large complex systems  (e.g., in condensed matter physics and quantum physics) even with using standard computers \cite{Orus2013,Sachdev09,Espig_2011,Cichocki-era}. In other words,
tensor networks,  can be considered as a diagrammatic language for capturing the internal structure
 of high-order tensor decompositions.

In contrast to the CPD or Tucker decompositions, that have only one single core tensor, TNs decompose a high-order tensor into several lower-order core tensors. The branches (leads, lines, edges)  connecting core tensors between each other correspond to contracted modes (and represent a TN rank), whereas lines that do not go from one tensor to another correspond to physical modes in the TN.
A tensor network  is a set of weakly connected core tensors, where some or all indices  are contracted according to some rules.

\begin{figure}[t!]
\centering
\includegraphics[width=8.6cm]{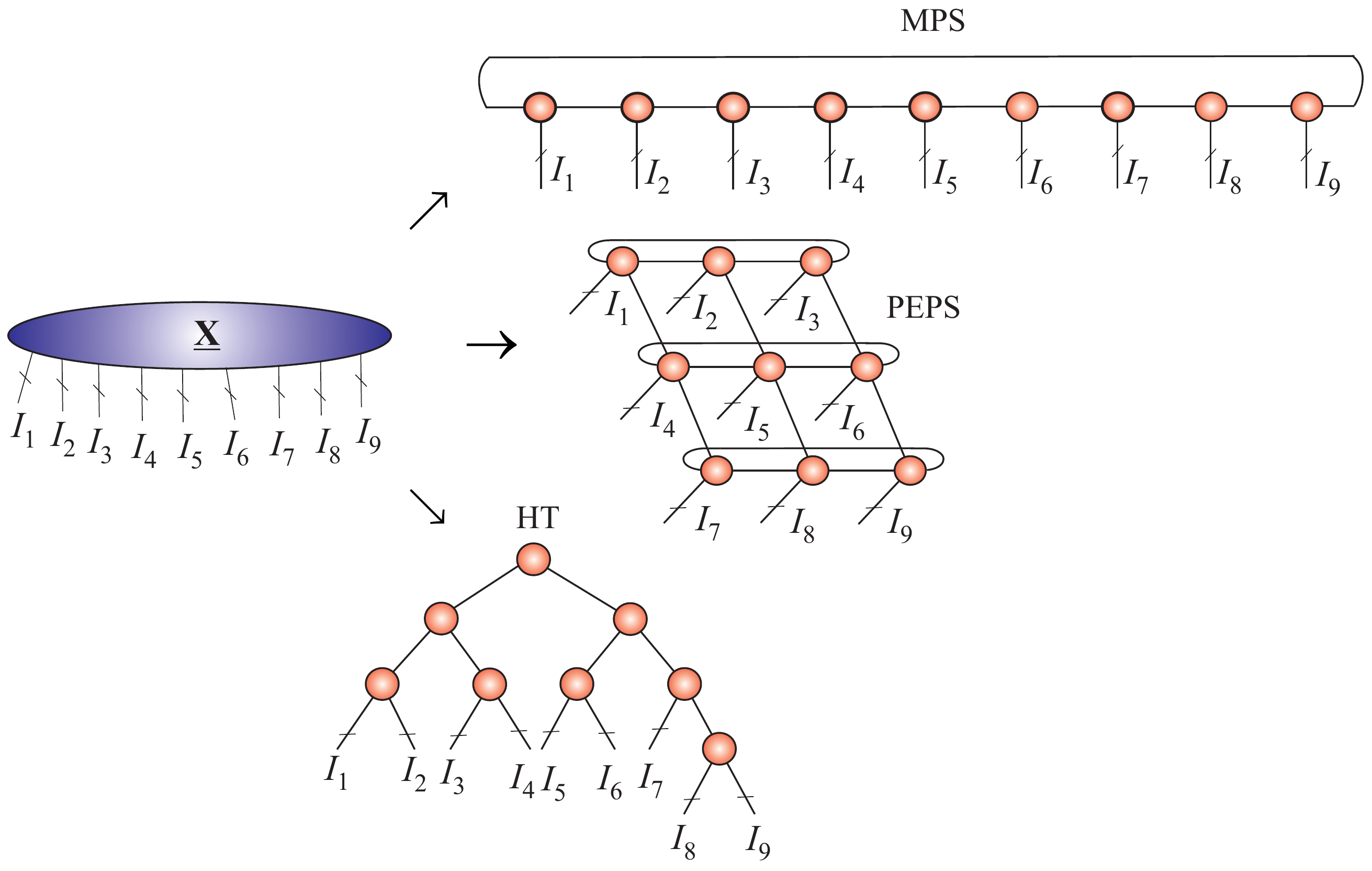}
\caption{Illustration of decomposition of 9th-order tensor $\underline \bX \in \Real^{I_1 \times I_2 \times \cdots \times I_9}$ into different forms of tensor networks (TNs): The Matrix Product State (MPS) with periodic boundary conditions (PBC), called also the Tensor Chain (TC), the Projected Entangled-Pair States (PEPS) with PBC and Hierarchical Tucker (HT) decomposition, which is equivalent to the  Tree Tensor Network State (TTNS).  In general, the objective is to decompose very high-order tensor into sparsely (weakly) connected low-order and low-rank tensors, typically 3rd-order and 4th-order tensors, called cores.}
\label{Fig:TN1}
\end{figure}

\begin{figure}[ht]
\centering
\includegraphics[width=7.6cm]{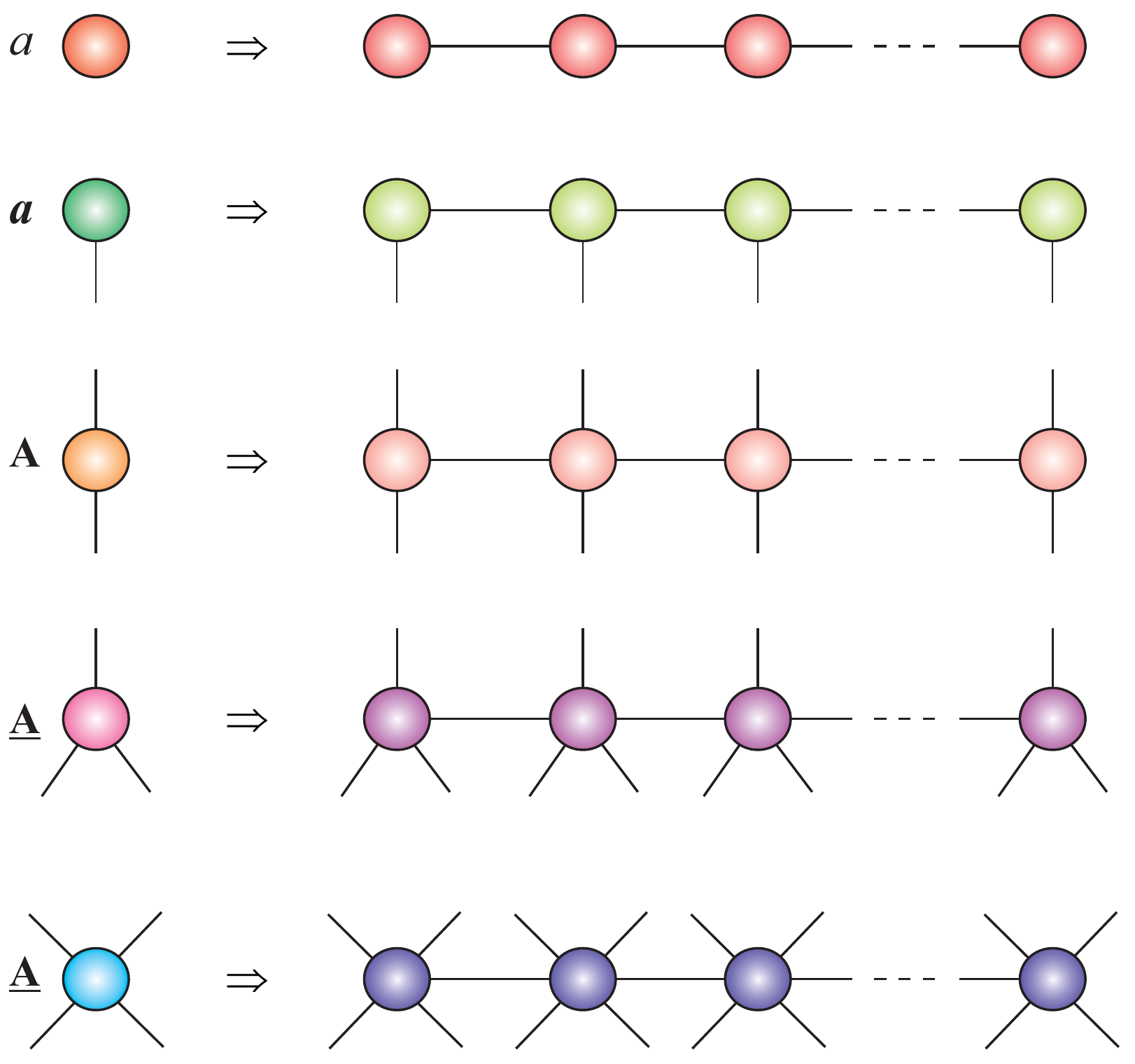}
\caption{Different forms of tensor train decompositions depending on input data: For scalar functions $a$, vectors $\ba$, matrices $\bA$, and 3rd-order and 4th-order tensors $\underline \bA$.}
\label{Fig:VariousTT}
\end{figure}

\begin{figure}
\centering
\includegraphics[width=8.1cm]{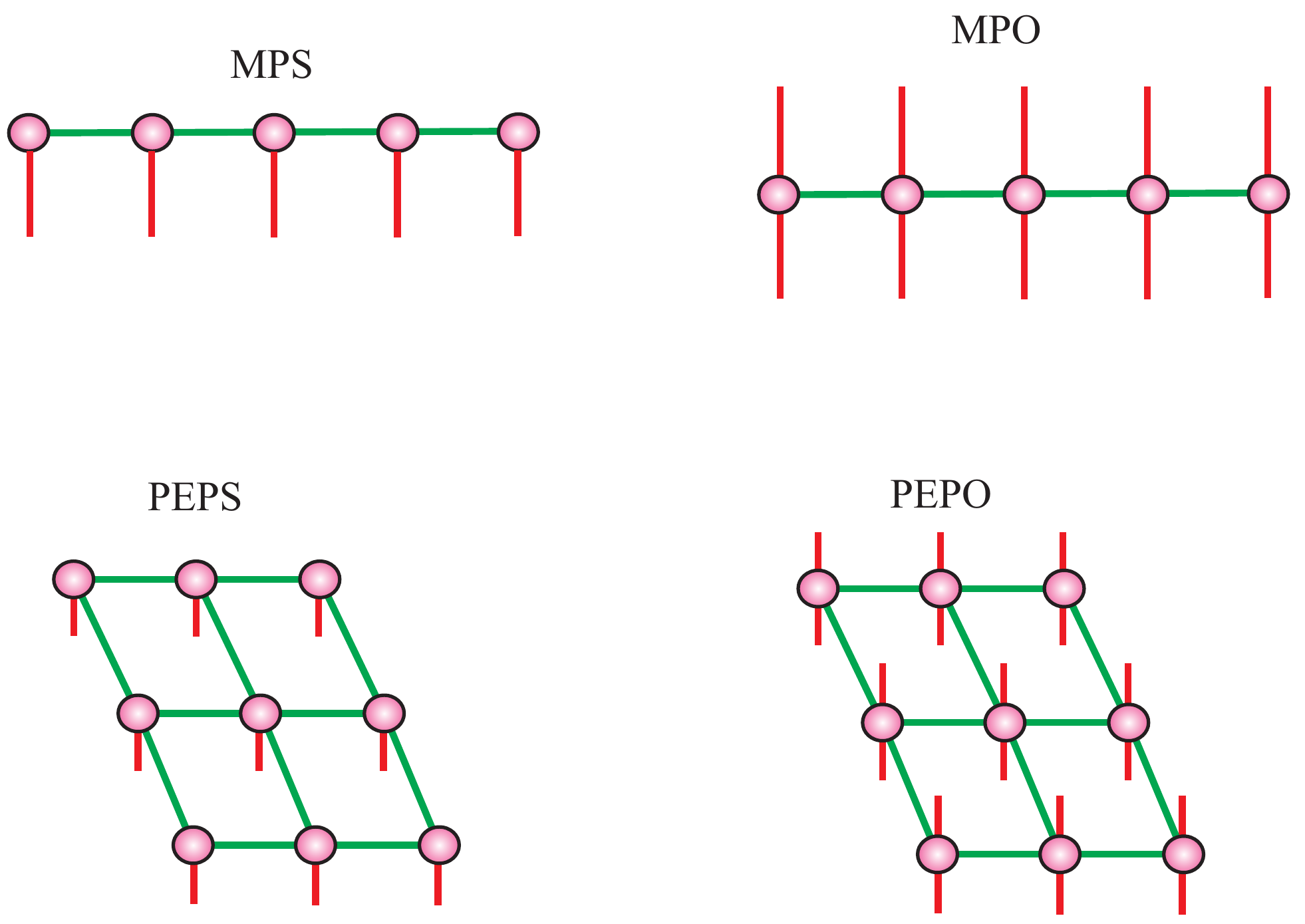}
\caption{Basic tensor networks with open boundary conditions (OBC):  The Matrix Product State (MPS) or (vector) Tensor Train (TT), the Matrix Product Operator (MPO) or Matrix TT, the Projected Entangled-Pair States (PEPS) or Tensor Product State (TPS) and the Projected Entangled-Pair Operators (PEPO).}
\label{Fig:TN11}
\end{figure}

\begin{figure}[ht!]
(a) Hierarchical Tucker (HT) or   Tree Tensor  Network State (TTNS) with 3rd-order and 4th-order cores \\
\begin{center}
\includegraphics[width=8.6cm]{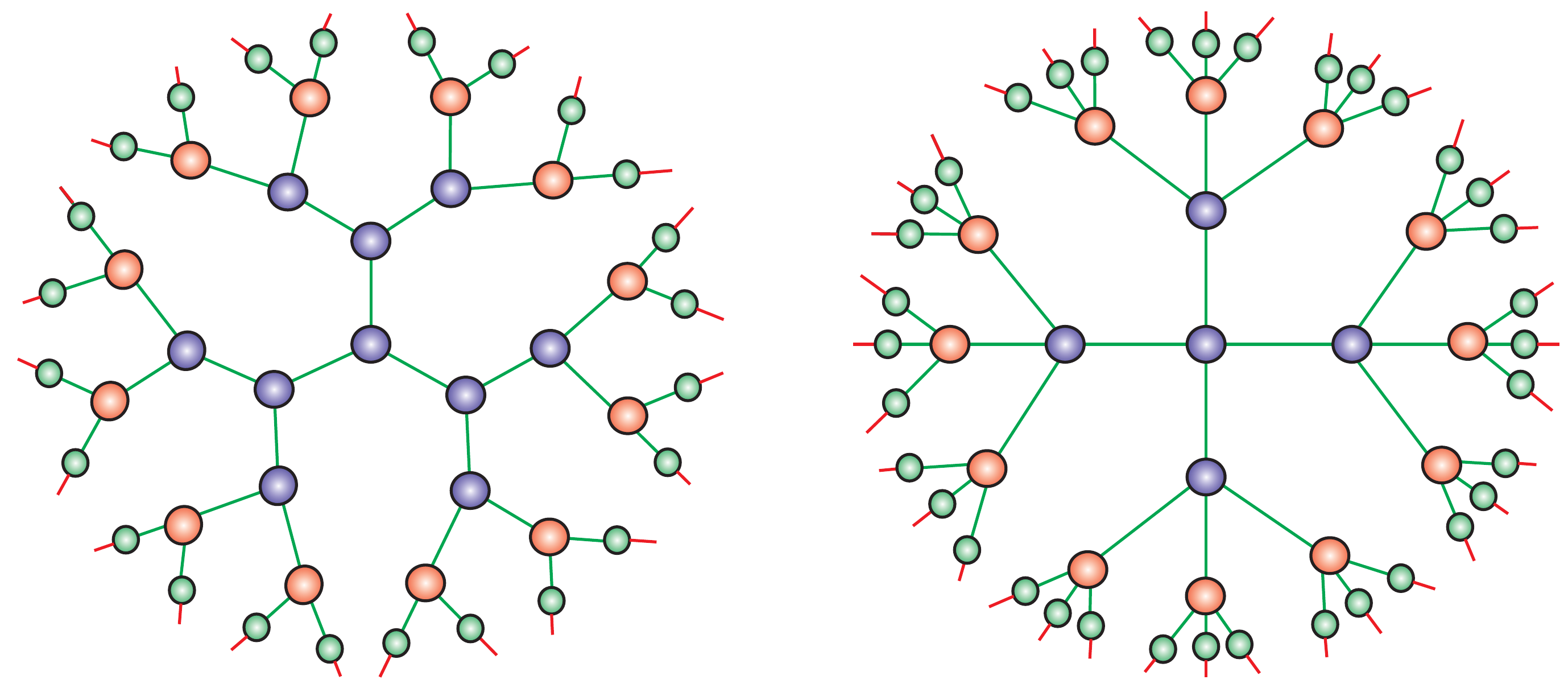}\\
\includegraphics[width=4.0cm]{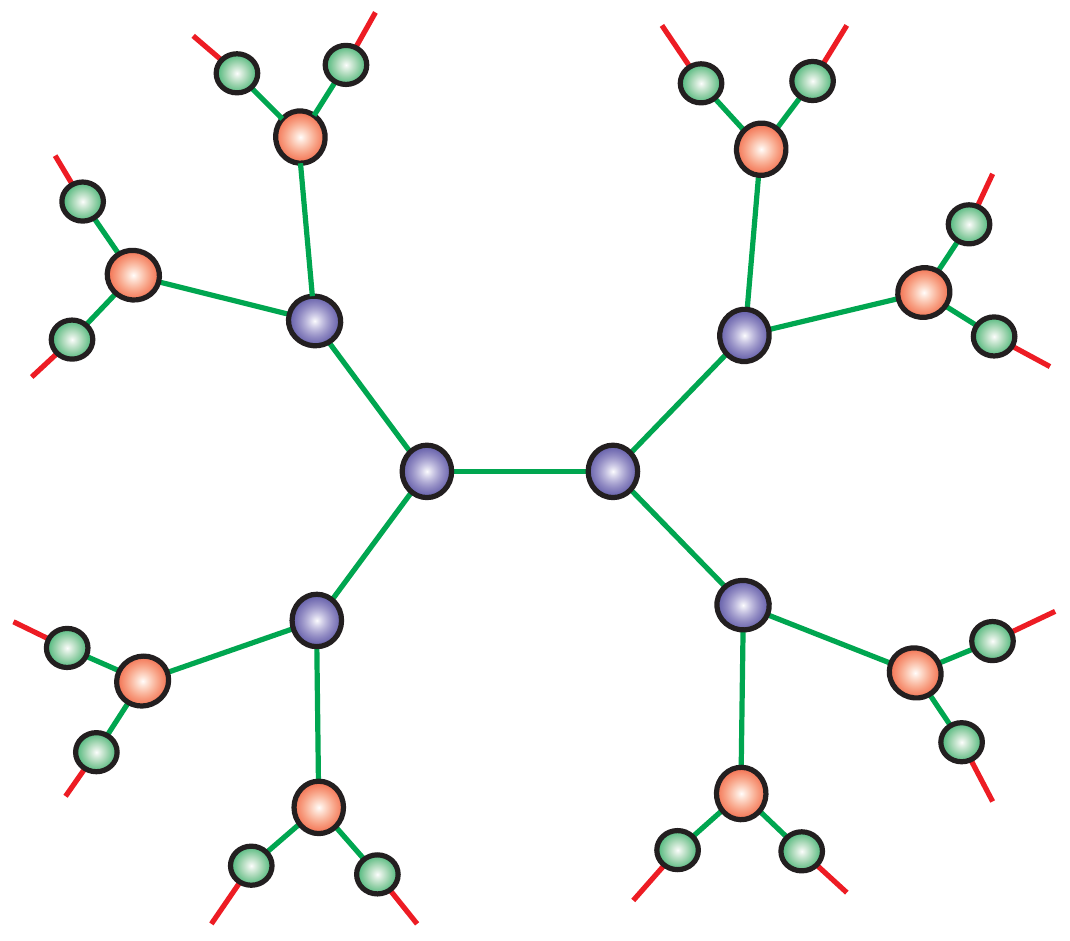}\\
\end{center}
(b) Honey-Comb lattice for a 16th-order data tensor \\
\vspace{0.1cm}
\begin{center}
\includegraphics[width=5.8cm]{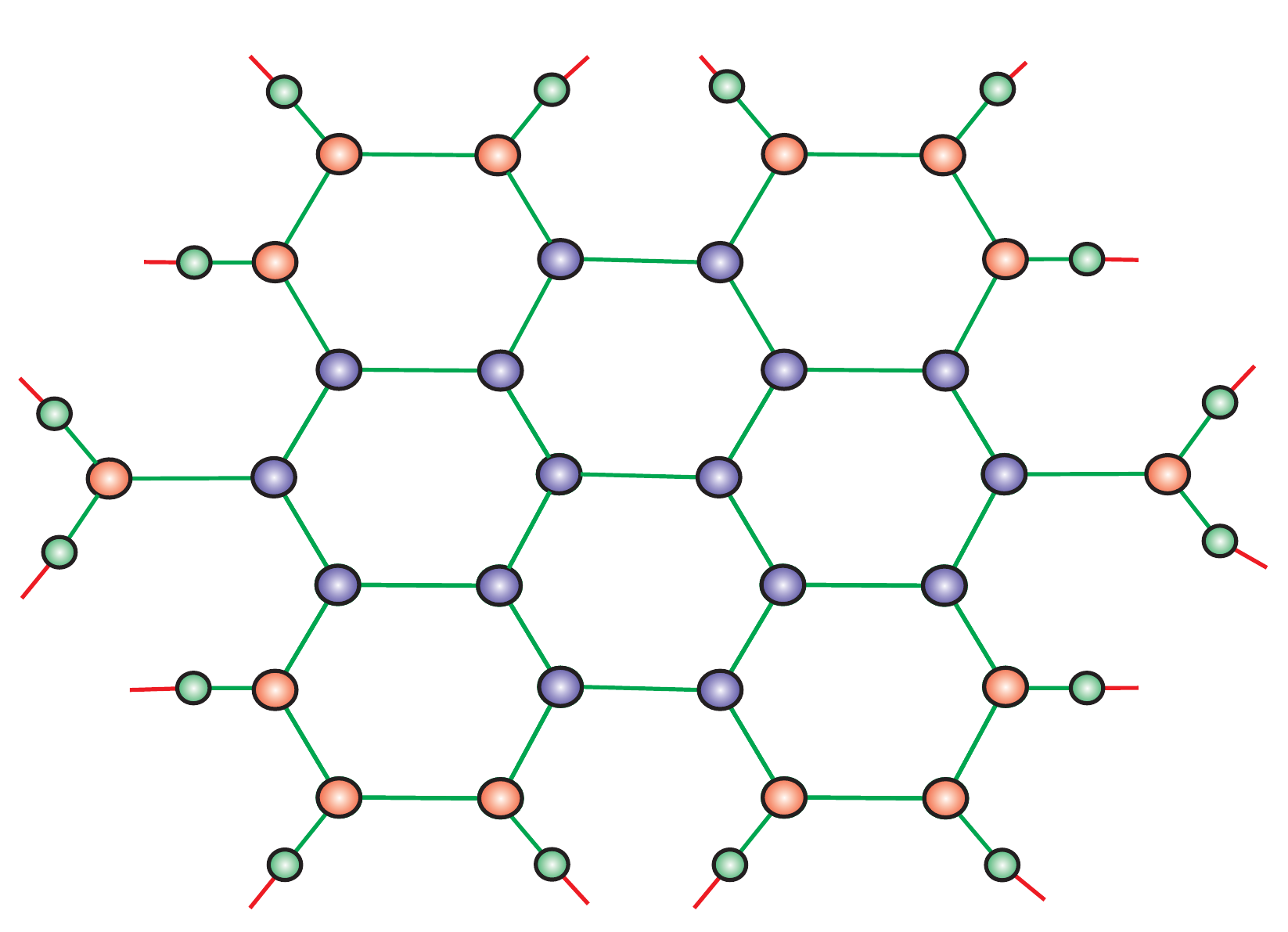}\\
\end{center}
(c) MERA for 8th-order tensor \\
\begin{center}
\includegraphics[width=3.4cm]{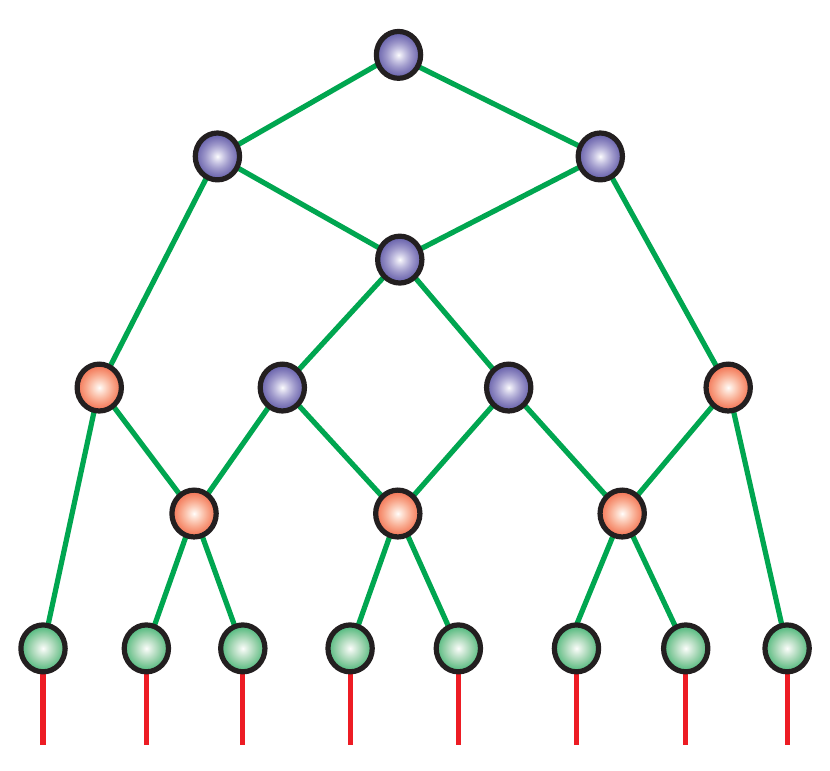}
\end{center}
\caption{Architectures of the fundamental TNs, which can be considered as distributed models of the Tucker-$N$ models. Green nodes denote factor matrices, while blue and red nodes denote cores.}
\label{Fig:COMB}
\end{figure}

 Some examples of basic tensor network diagrams
are given in Figs. \ref{Fig:TN1}, \ref{Fig:VariousTT}, \ref{Fig:TN11}, and \ref{Fig:COMB} \cite{Cichocki-era,Orus11}.
A tensor network may not contain any loops,
i.e., any edges connecting a node with itself.
If a tensor network is a binary tree, i.e., if it does not contain any cycles (loops), each of its edges
splits the modes of the data tensor  into two or more groups, which is  related to the
suitable matricization of the tensor \cite{Grasedyck-rev,Grasedyck-Hrank}.
A tree tensor network, whose all nodes have degree 3 or 4, corresponds
to a Hierarchical Tucker (HT) decomposition  of the  tensor  illustrated in Fig. \ref{Fig:COMB} (a).
The HT decompositions in the
 numerical analysis community have been introduced by Hackbusch and K{\"u}hn  \cite{HackbuschHT09} and Grasedyck \cite{hTucker1} (see also \cite{Grasedyck-rev,Uschmajew-Vander2013,kressner2012htucker,kressner2014htucker,Lubich-Schneider13} and references therein).
The general construction of the  HT decomposition requires a hierarchical splitting of
the modes (with sizes $I_1,I_2,\ldots,I_N$).
The  construction of  Hierarchical Tucker format relies on the
notion of a dimension tree, chosen
{\it a priori}, which specifies the topology of the HT decomposition. Intuitively, the dimension
tree specifies which groups of  modes
are ``separated'' from other groups of modes, where sequential HT decomposition can be performed
via (truncated) SVD applied to unfolded matrices \cite{Grasedyck-rev}.

%
%
%

The Tensor Train (TT) format  proposed  
in the numerical analysis community by Oseledets and Tyrtyshnikow \cite{OseledetsTT09}
 (see also \cite{Vidal03,OseledetsTT11,NLee-Cich14,Kazeev-CME-2014,LubichTITT14,Bigoni-STT-2014}) can be interpreted as
a special case of the HT, where all nodes of the underlying tensor network are aligned
and where, moreover, the leaf matrices are assumed to be identities (and thus need not
be stored). An advantage of the TT format is its simpler practical implementation, as no binary
tree need be involved.
 The Tensor Trains  \cite{OseledetsTT09,OseledetsTT11,OseledetsTZ11}, called also Matrix Product States (MPS) in   quantum information theory \cite{MPS2007,verstraete08-MPS,Orus2013,Schollwock13,Huckle2013}, is the simplest TN model{\footnote{In fact, the TT was rediscovered several times under different names: MPS, valence bond states and density matrix renormalization group (DMRG). The DMRG usually means not only tensor format but also power-full computational algorithms (see \cite{schollwock11-DMRG} and references therein).}}.

For some very high-order data tensors  it has been observed
that the ranks $R_n$ of 3rd-order tensors  increase rapidly with the order of the tensor,
for any choice of tensor network that is a tree (including TT and  HT decompositions) \cite{kressner2012htucker}.

For such cases, PEPS and the Multi-scale Entanglement Renormalization Ansatz (MERA) tensor networks
can be used   which contain cycles, but have  hierarchical structures (see Fig. \ref{Fig:COMB}) (c). For the PEPS and MERA TNs the ranks
can be kept considerably smaller, at the cost of  employing 5th and 4th-order core  tensors and consequently a higher computational complexity w.r.t. their  ranks \cite{Evenbly-Vidal09alg,MERA08}.

Some interesting   connections between  tensor networks and graphical models used extensively in machine learning and statistics as shown  in Table \ref{Table:Morton}
 \cite{Morton12,Morton14,MRF-TT14,Kazeev-CME-2014,Critchthesis}.
Despite clear analogy, more research is needed to find more deep and precise  relationships \cite{Critchthesis}.

\begin{table}[t]
\setlength{\tabcolsep}{2pt}
\renewcommand{\arraystretch}{1.5}
\centering
\caption{Similarities and links between tensor networks (TNs) and
graphical models used in Machine Learning (ML) and Statistics.
The categories are not exactly the same, but they  closely correspond.}
{\shadingbox{
\begin{tabular}{p{3.5cm}|p{5cm}}
\hline
Tensor Networks & Graphical Models in ML/Statistics \\
\hline
TT/MPS & Hidden Markov Models (HMM) \\
HT/TTNS & Gaussian Mixture Model (GMM) \\
TNS/PEPS & Markov Random Field (MRF) and Conditional Random Field (CRF) \\
MERA & Deep Belief Networks (DBN) \\
DMRG and MALS Algs. & Forward-Backward Algs., Block  Nonlinear Gauss-Seidel Methods \\
\hline
\end{tabular}
}}
\label{Table:Morton}
\end{table}

\subsection{\bf Changing the Structure of Tensor Networks}

 One advantage of  a graphical representation of a tensor network is that it
 allows us  to perform even most c complex mathematical operations in intuitive and easy to understand way. Another important advantage
 is the ability to modify or optimize a TN structure, that is, to change  its topology,
 preserving physical modes unchanged.
In   fact, in some applications it is quite  useful to modify  the topology of a tensor network with or  without approximation  by providing  simplified or more convenient graphical representation of the same higher-order data tensor \cite{Handschuh12,Zhao-Xie-RTNS10,hubener2010concatenated}.
For instance, tensor networks may consist of many cycles, those can t be  reduced or  completely eliminated  in order to reduce
computational complexity of contraction of  core tensors and  to provide stability of computation.
Again, observe a strong link with loop elimination in control theory, in addition
  tensor networks having many cycles
 may  not  admit  stable algorithm.
By changing the topology to a tree structure (TT/HT models),
we  can  often  reduce complexity of computation and improve stability of algorithms.

\begin{figure}[p!]
(a)\\
\includegraphics[width=8.9cm]{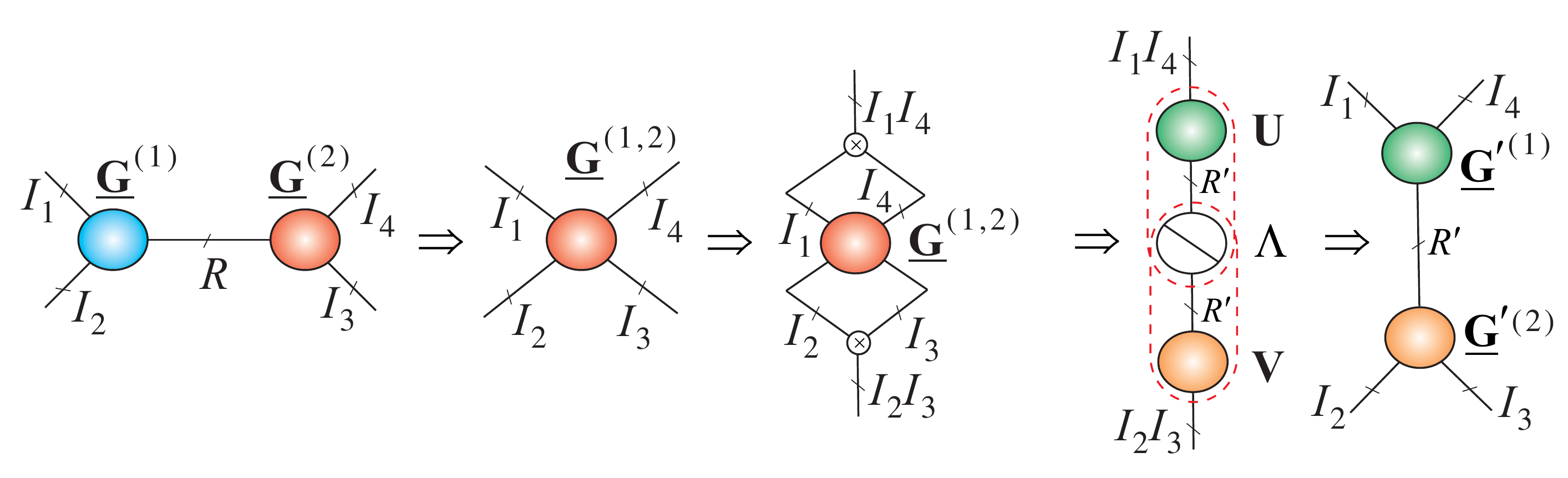}\\
(b)\\
\includegraphics[width=8.6cm]{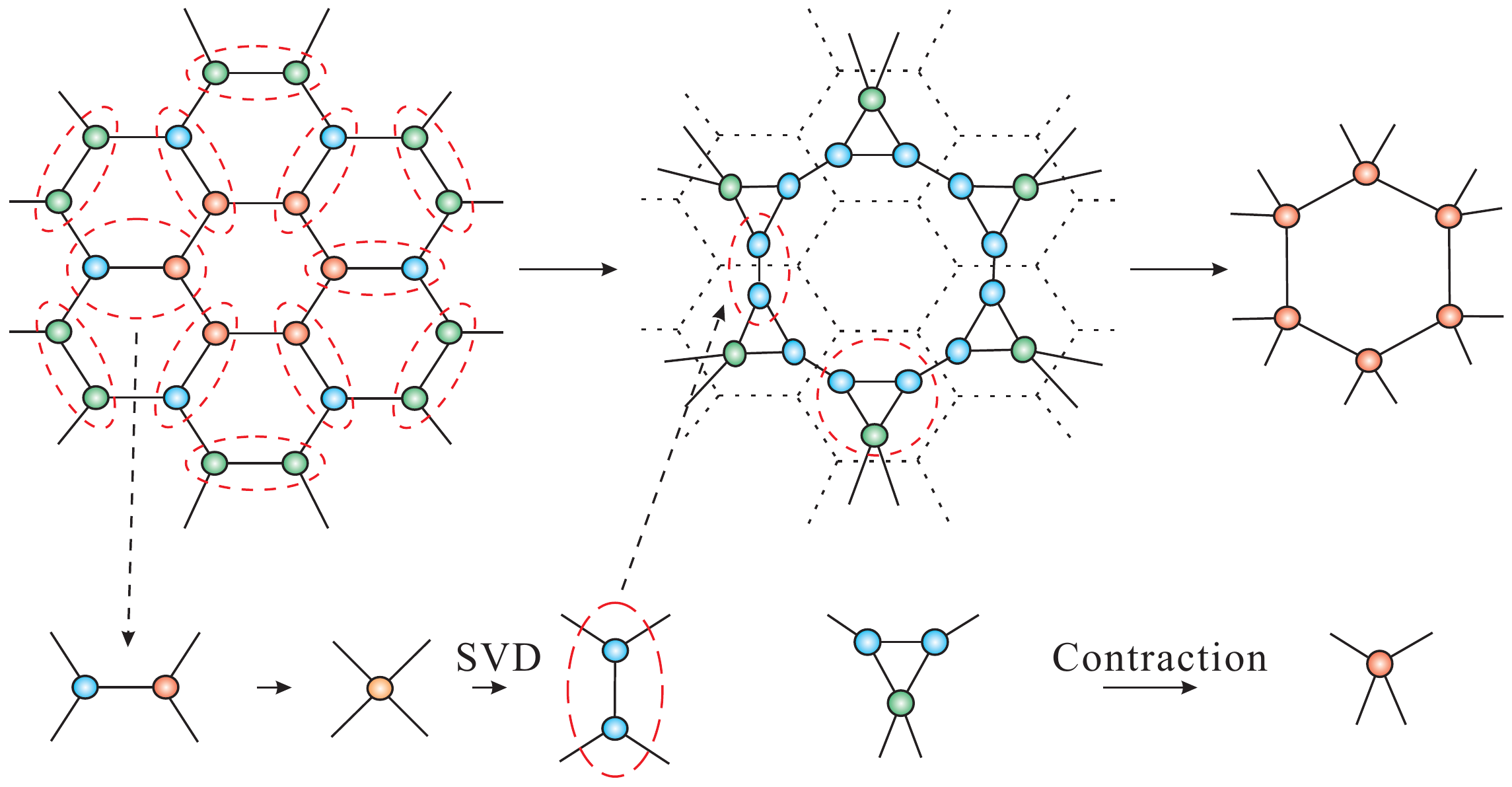}\\
(c)\\
\includegraphics[width=8.6cm]{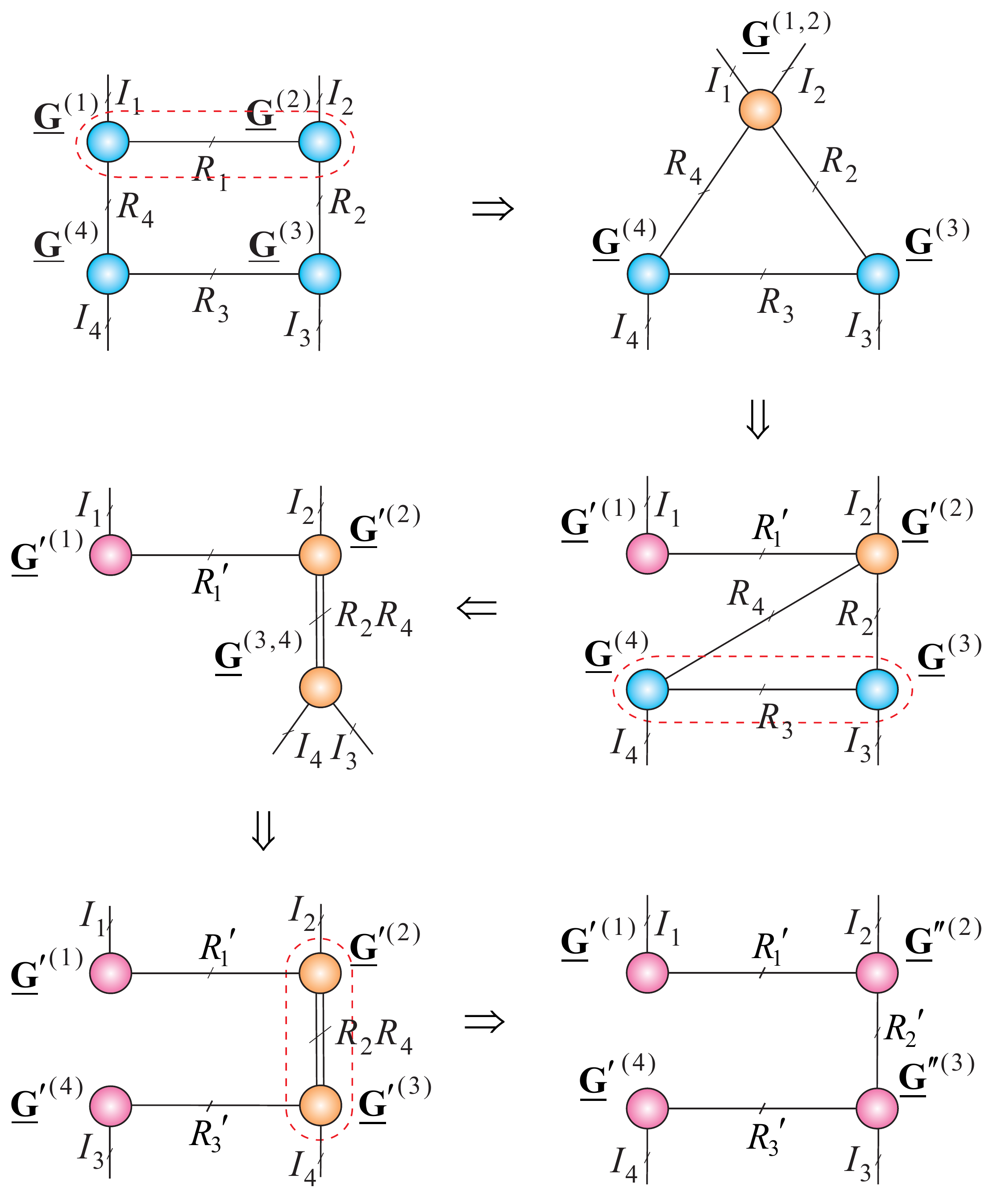}\\
\caption{Illustration of  basic transformation of tensors via: (a) Contraction, unfolding, matrix factorization (SVD) and reshaping of matrices back into tensors. (b) Transformation of Honey-Comb lattice into Tensor Chain (TC) via tensor contractions and the SVD. (c) Transformation of the TC,  i.e., a TT/MPS with periodic boundary conditions (PBC) to the standard TT, i.e., the MPS with  open boundary conditions (OBC). }
\label{Fig:TCTT}
\end{figure}

Performing contraction of core tensors iteratively for  tree--structured tensor networks
has usually a much smaller complexity than tensor networks containing many cycles.
One could transform a specific tensor network with cycles into a tree structure, perform stable computations{\footnote{The TT decomposition is stable in the sense that the best approximation of a data tensor with bounded TT-ranks always exist and a quasi-optimal approximation can be computed by a sequence of truncated SVDs of suitably reshaping matrices of cores \cite{Vidal03,OseledetsTT11}.}}, with it and re-transform it back to the original structure if necessary.
Furthermore, in the cases that we need to compare or analyze a  set  of blocks of tensor data,
  it is important that  such
tensors are represented by  the same or very similar structures  to analyze link or correlation between  them or detect common  cores or hidden components.
 Performing  such analysis with differently structured tensor networks is in general difficult or even impossible.

A Tensor network can be relatively easily  transformed  from one form to another one via tensor contractions, reshaping and basic matrix factorizations,  typically using SVD \cite{Vidal03,OseledetsTT11}.
The basic approach to modify tensor structure is to perform: sequential core contractions, unfolding contracting tensors into matrices, performing matrix factorizations (typically, SVD) and finally reshaping matrices back to  new core tensors.
These principles are illustrated graphically in Figs \ref{Fig:TCTT} (a), (b), (c).

For example, in Fig \ref{Fig:TCTT} (a) in the
first step we perform  a contraction of two core tensors $\underline \bG^{(1)} \in \Real^{I_1 \times I_2 \times R}$ and $\underline \bG^{(2)}  \in \Real^{R \times I_3 \times I_4}$, as:
\be
\underline \bG^{(1,2)} =\underline \bG^{(1)} \times_3^1 \; \underline \bG^{(2)} \in \Real^{I_1 \times I_2 \times I_3 \times I_4},
\ee
with entries  $g_{i_1,i_2,i_3,i_4}^{(1,2)} = \sum_{r=1}
^R g_{i_1,i_2,r}^{(1)} \; g_{r,i_3,i_4}^{(2)}$.
In the next step, we transform  the tensor
$\underline \bG^{(1,2)}$ into a matrix via unfolding and low-rank matrix factorization  via the SVD
\be
\bG^{(1,2)}_{\overline{i_1,i_4};\overline{i_2,i_3}} \cong \bU \mbi \Sigma \bV^T \in \Real^{I_1 I_4 \times I_2 I_3}.
\ee
In the last step, we reshape factor matrices $\bU \mbi \Sigma^{1/2} \in \Real^{I_1 I_4 \times R'}$ and $\bV \mbi \Sigma^{1/2} \in \Real^{R' \times I_2 I_3}$ back to new core tensors: $\underline \bG^{'(1)} \in \Real^{I_1 \times R' \times I_4 }$ and $\underline \bG^{'(2)} \in \Real^{I_2 \times I_3 \times R'}$.

The above procedure has been applied  in Fig. \ref{Fig:TCTT} (b)
to transform Honey-Comb lattice  into tensor chain (TC) along with tensor
contraction of three cores \cite{Zhao-Xie-RTNS10}.

In Fig. \ref{Fig:TCTT} (c) we  have illustrated how to convert tensor chain (TC) into TT/MPS with OBC, by contracting sequentially two core tensors, unfolding them,  applying SVD  and reshaping matrices back into core tensors \cite{Handschuh12}.
More precisely, in the
first step, we perform  a contraction of two tensors $\underline \bG^{(1)} \in \Real^{I_1 \times R_4 \times R_1}$ and $\underline \bG^{(2)}  \in \Real^{R_1 \times  R_2 \times I_2}$, as:
\be
\underline \bG^{(1,2)} = \underline \bG^{(1)} \times_3^1 \; \underline \bG^{(2)} \in \Real^{I_1 \times R_4 \times R_2 \times I_2},
\ee
with entries  $g_{i_1,r_4,r_2,i_2}^{(1,2)} = \sum_{r_1=1}
^{R_1} g_{i_1,r_4,r_1}^{(1)} \; g_{r_1,r_2,i_2}^{(2)}$.
In the next step, we can transform this tensor
$\underline \bG^{(1,2)}$ into a matrix  in order to perform the truncated SVD:
\be
\bG^{(1,2)}_{i_1 \, ; \, \overline{r_4,r_2,i_2}} \cong \bU \mbi \Sigma \bV^T \in \Real^{I_1  \times R_4 R_2 I_2}.
\ee
In the next step, we reshape orthogonal matrices $\bU \mbi \Sigma^{1/2} \in \Real^{I_1 \times R'_1}$ and $\bV \mbi \Sigma^{1/2} \in \Real^{R'_1 \times R_4 R_2 I_2}$ back to core tensors: $\underline \bG^{'(1)}=\bU \mbi \Sigma^{1/2} \in \Real^{1 \times I_1 \times R'_1}$ and $\underline \bG^{'(2)} \in \Real^{R'_1 \times R_4 \times  R_2 \times I_2}$.
The procedure is repeated again and again for different pair of cores as illustrated in the Fig. \ref{Fig:TCTT} (c).

\subsection{\bf Distributed (Concatenated) Representation of Tensors}

\begin{figure*}[p]
(a) Tensor Train (TT) model -- MPS/MPO with the Open Boundary Conditions (OBC)
\begin{center}
\includegraphics[width=12.8cm]{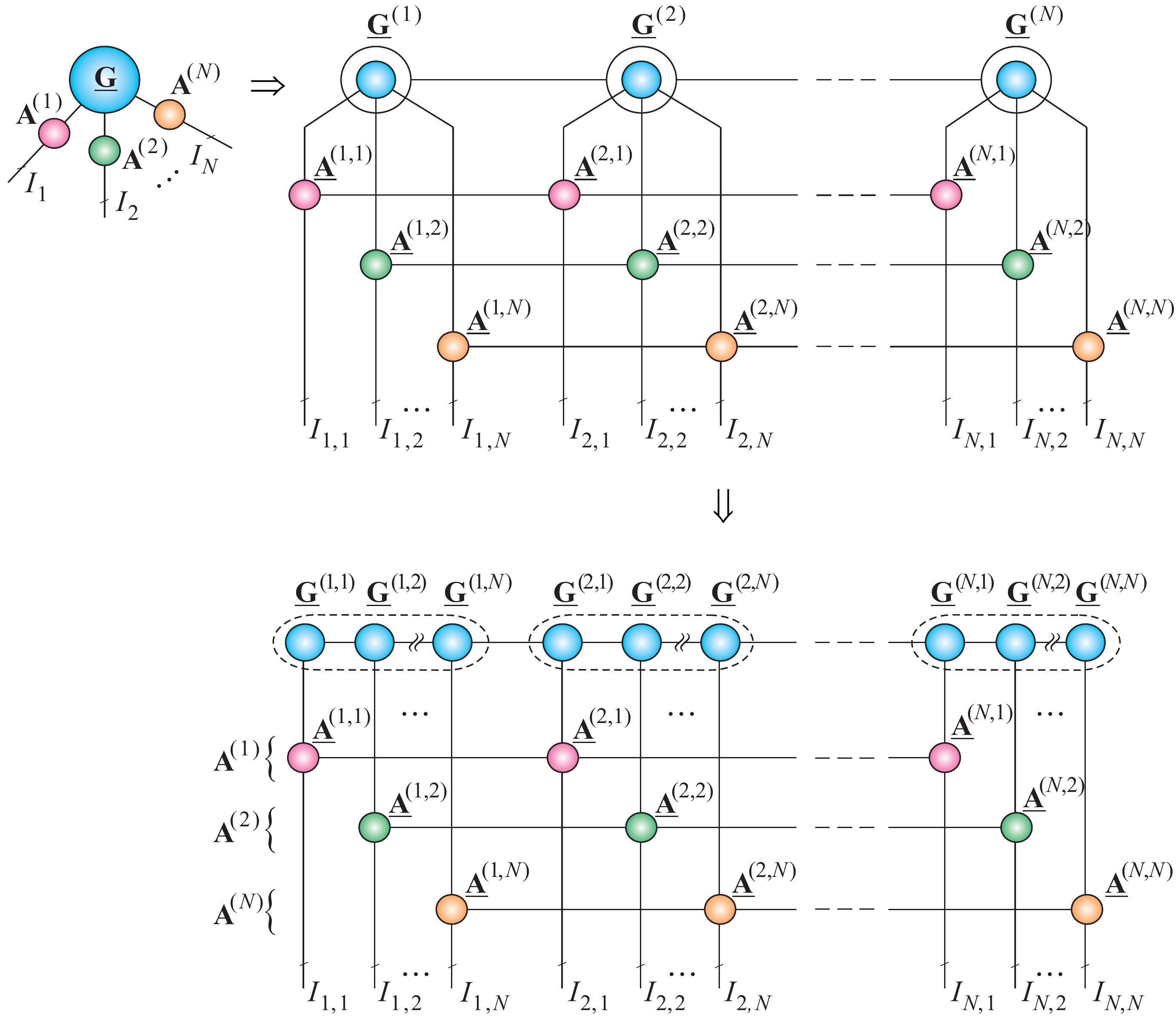}\\
\end{center}
(b) Tensor Chain (TC) model -- MPS/MPO with the Periodic Boundary Conditions (PBC)
\begin{center}
\includegraphics[width=9.6cm]{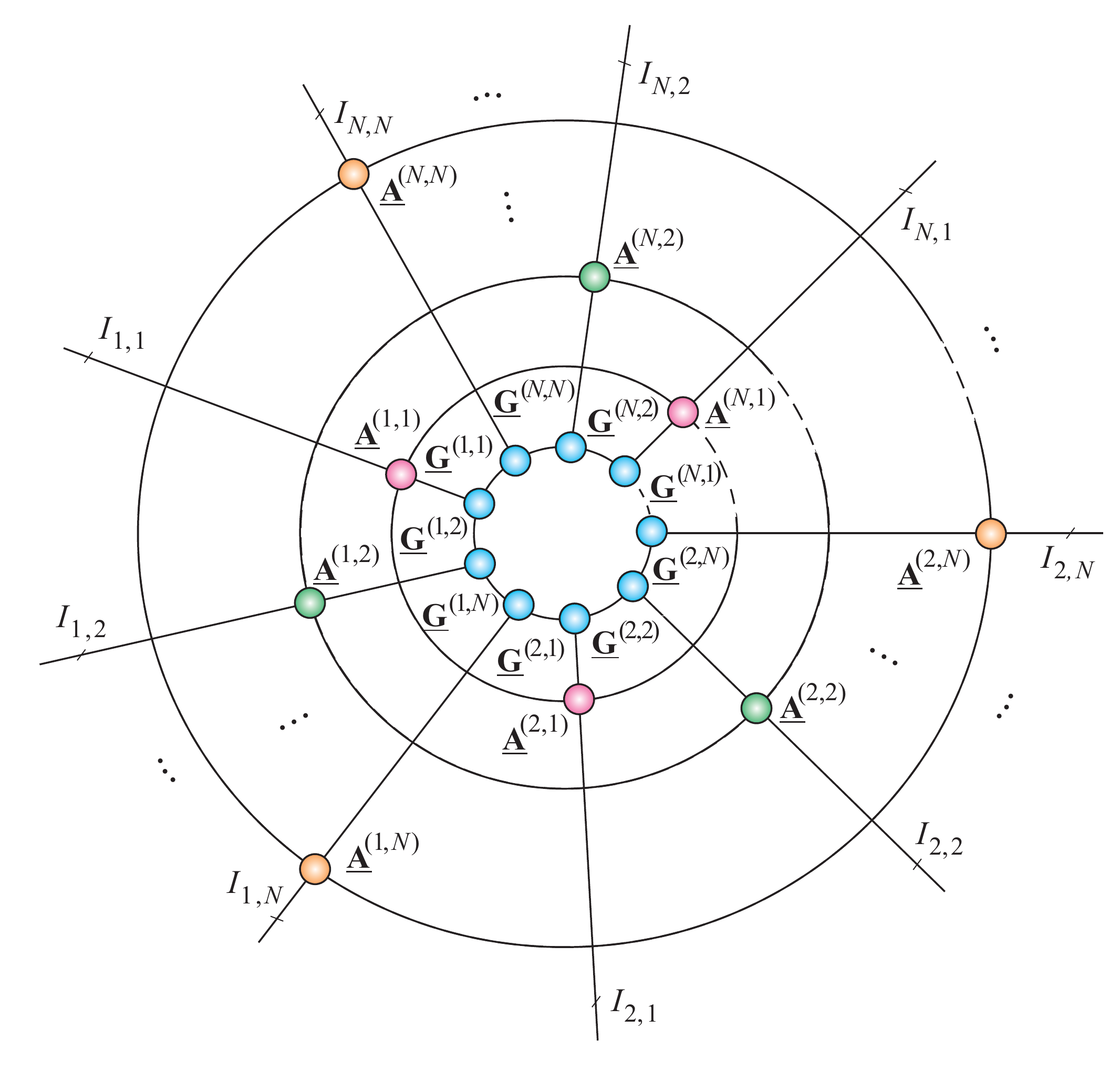}\\
\end{center}
\caption{New distributed models  of the Tucker-$N$ decomposition $\underline \bX =\underline \bG \times_1 \bA^{(1)} \times_2 \bA^{(2)} \cdots \times_N \bA^{(N)} \in \Real^{I_1 \times I_2 \times \cdots I_N}$, with $I_n=I_{1,n} I_{2,n} \cdots I_{N,n},\;$ $(n=1,2,\ldots,N)$.}
\label{Fig:Tucker-TT2}
\end{figure*}

\begin{figure} 
\begin{center}
\includegraphics[width=8.6cm,height=8.5cm]{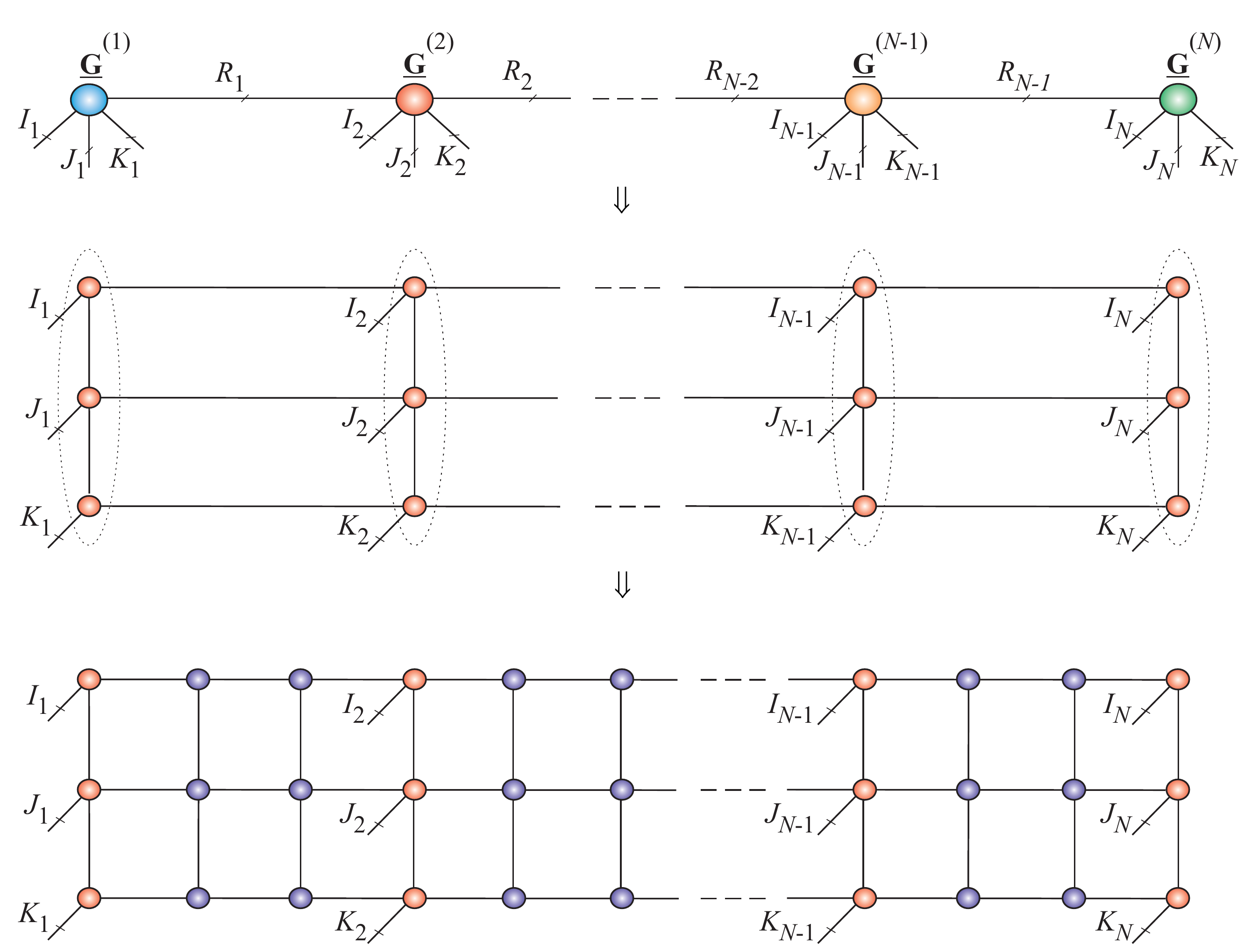}\\
\end{center}
\begin{center}
\includegraphics[width=8.1cm,height=6.9cm]{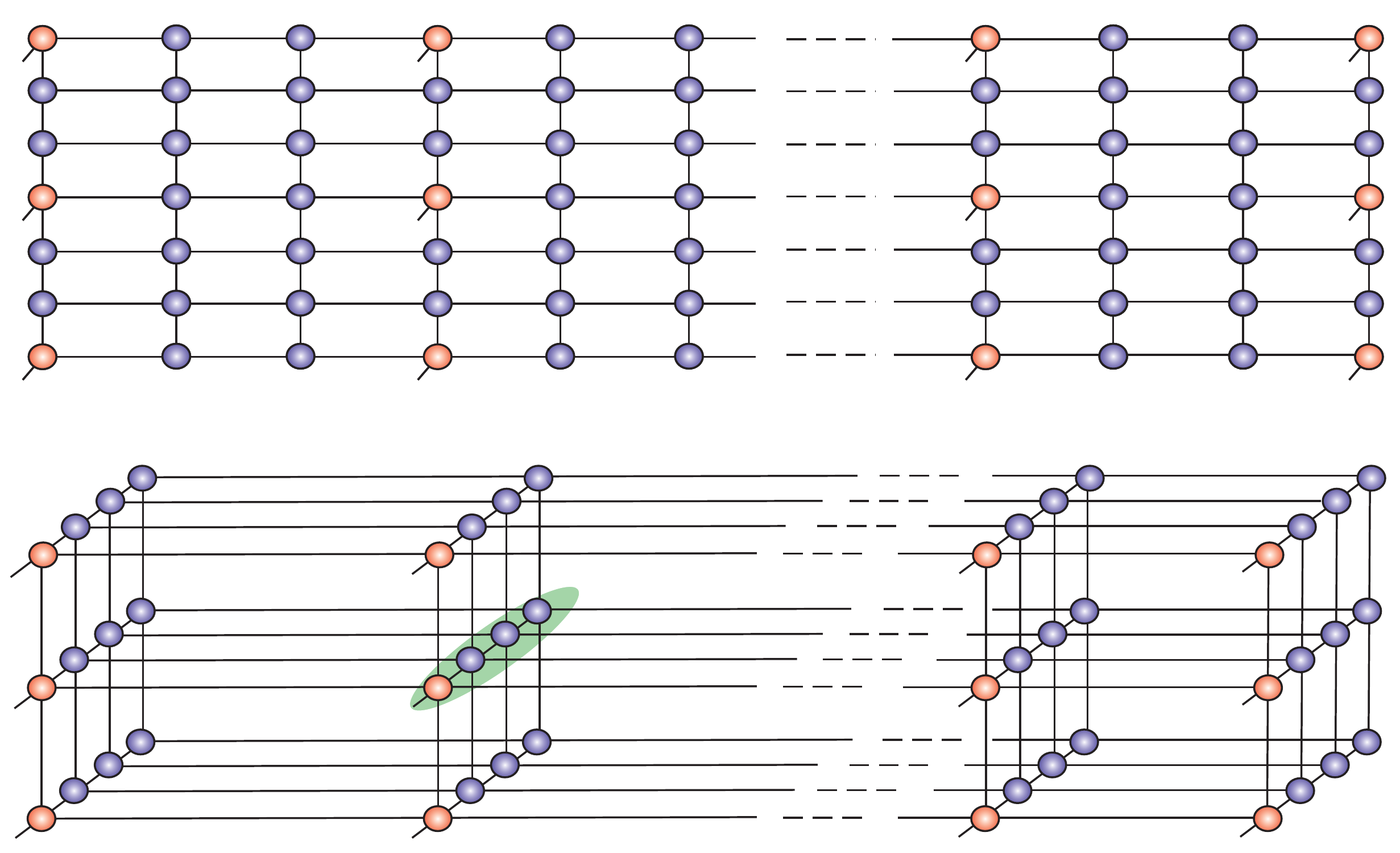}\\
\end{center}
\caption{Graphical  representation of TT via Tensor Product State (TPS) or equivalently PEPS for a large-scale data tensor and its transformation to  distributed 2D and 3D PEPS \cite{hubener2010concatenated}.}
\label{Fig:TTPEPS}
\end{figure}

A simple approach to reduce the size or rank of core tensors is to apply distributed tensor networks (DTNs), which consists of two kind of cores (nodes): internal nodes which has no free edges and external nodes which have free edges representing natural (physical) indices of a data tensor as illustrated in
Figs. \ref{Fig:COMB} and \ref{Fig:Tucker-TT2}. A simple idea  is that each of the core tensor in an original TN is itself repeatedly replaced by another TN (see Fig.  \ref{Fig:TTPEPS}),
resulting in another TN in which only some core tensors are associated with physical (natural) modes  of the original data tensor \cite{hubener2010concatenated}.

The  main advantage of DTNs  is that  the size of
each of the core tensors in the internal tensor network structure is usually much smaller
than the initial core tensor so consequently the total number of parameters can be  reduced
\cite{hubener2010concatenated}.
 However, it should be noted that the contraction of the resulting tensor network becomes more difficult
when compared to the initial tree structure. This is due to the fact that the distributed tensor
network contains loops.

Many algorithms applied to  tensor networks scale with the size $R_k$ or $I_k$ of the core tensors of the network. In spite of
the usually polynomial scaling of these algorithms, the computations quickly become intractable
for increasing $R_k$, so that a network containing core tensors with  small dimensions are favorable in
general. See as examples the distributed Tucker models shown in Fig. \ref{Fig:Tucker-TT2} (a) and (b).

\section{\bf Tensorization   -- Blessing of Dimensionality}
\label{sect:tensorization}

\begin{figure}[t]
\includegraphics[width=8.6cm]{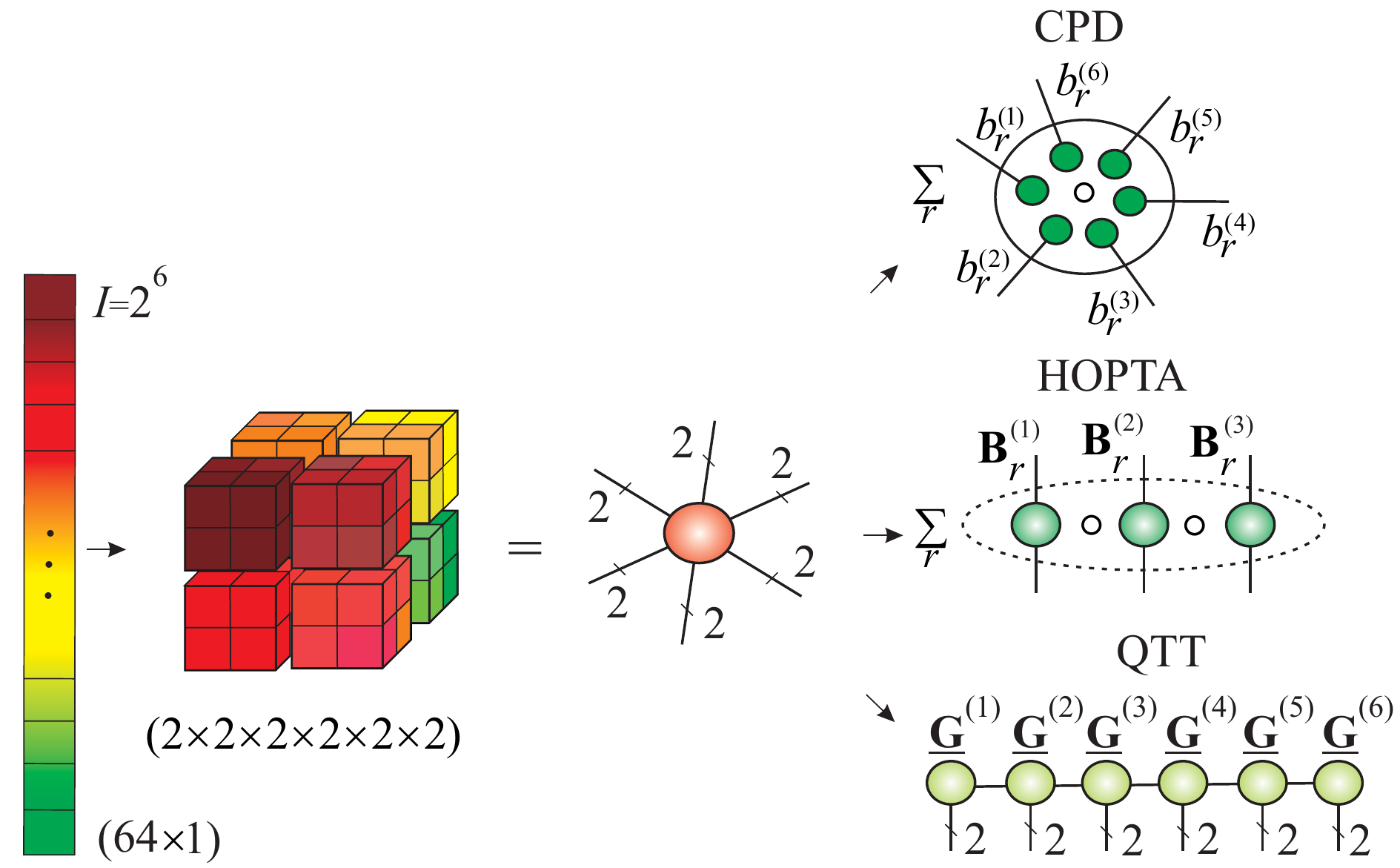}
\caption{The conceptual illustration of tensorization of a large-scale vector
into a higher-order quantized tensor.
In order to achieve super-compression we need to apply a suitable  tensor decomposition: e.g., CPD decomposition into rank-1 tensors $\underline \bX \cong \sum_{r=1}^R \bb_r^{(1)} \circ \bb_r^{(2)} \cdots \circ \bb_r^{(6)}$,  Hierarchical Outer Product Tensor Approximation (HOPTA) using rank-$q$  terms: $\underline \bX \cong \sum_{\tilde r=1}^{\tilde R} \bB_{\tilde r}^{(1)} \circ \bB_{\tilde r}^{(2)} \circ \bB_{\tilde r}^{(3)}$ or  quantized TT (QTT) using 3rd-order cores: $\underline \bX = \underline \bG^{(1)} \times^1_3 \underline \bG^{(2)} \times^1_3 \cdots \times^1_3  \underline \bG^{(6)}$.}
\label{Fig:Tensorization}
\end{figure}

The procedure of creating a higher-order tensor from lower-order original data is referred to as {\em tensorization}.  In other words,  lower-order data tensors can be  reshaped (reformatted) into high-order tensors. The purpose of a such tensorization or reshaping is to achieve a low-rank approximation with high level of compression.
For example, big vectors, matrices even low-order tensors can be easily  tensorized to very high-order tensors, then efficiently compressed by applying a suitable tensor network decomposition; this is the underlying principle for big data analysis  \cite{Cich-Lath,Cichocki-era,Oseledets10,Khoromskij-SC} (see also Figs. ~\ref{Fig:Tensorization} - \ref{Fig:TT=MPS=MPO}).

\subsection{\bf Curse of Dimensionality}

The term curse of dimensionality, in the context of tensors,  refers to the fact that the number of elements of an $N$th-order $(I \times I \times \cdots \times I)$ tensor, $I^N$, grows exponentially with the tensor order $N$.
%
%
For example, for the Tucker decomposition the number of entries of an original data tensor but also a core tensor  scales exponentially in the tensor order, for instance, the number of entries of an $N$th-order  $(R \times R \times \cdots \times R)$ core tensor is   $R^N$.

If all computations are performed on a CP tensor format and not on the raw data tensor itself, then instead of the original $I^N$  raw data entries, the number of parameters in a CP decomposition reduces to $N R I$, which scales linearly in $N$ and $I$. This effectively bypasses the curse of dimensionality, however the CP approximation may provide a poor fit to the data and may involve numerical problems, since existing CPD algorithms are not stable for high-order tensors.
  In this paper we exploit TT decompositions which are stable and robust with ability   to control an
  approximation error i.e., to achieve any desired accuracy of TT approximation  \cite{OseledetsTT11,Khoromskij-TT}.
The main idea  of using  low-rank tensor-structured approximations is to
reduce the complexity of computation and relax or avoid the curse of dimensionality.

 \subsection{\bf Quantized Tensor Networks}

 The curse of dimensionality can be  overcome relatively easily
 through quantized tensor networks, which represent a tensor of possibly very high-order as a set of sparsely interconnected of  low dimensions (typically, 3rd-order) cores \cite{Khoromskij-TT,Oseledets10}.
  The concept of quantized tensor networks was first proposed by
   Oseledets \cite{Oseledets10} and Khoromskij \cite{Khoromskij-SC}.

\begin{figure*}[p]
(a)\\
\begin{center}
\includegraphics[width=8.6cm]{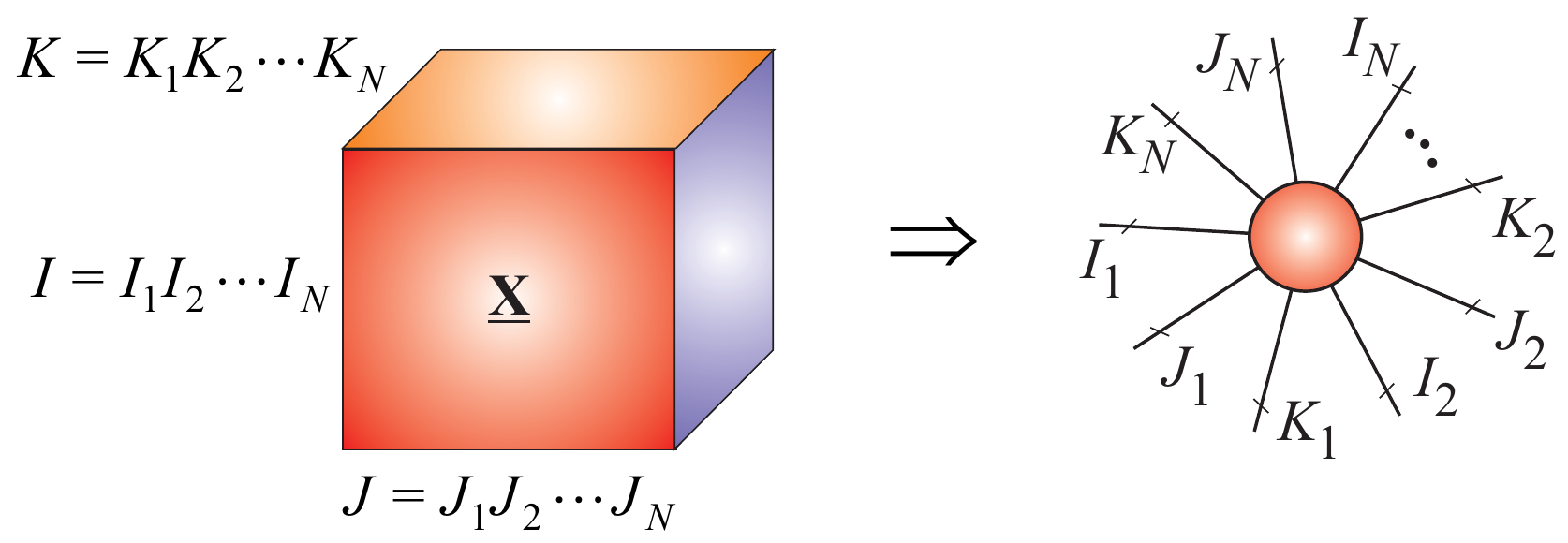}\\
\end{center}
(b)\\
\begin{center}
\includegraphics[width=17.5cm,height=12.1cm]{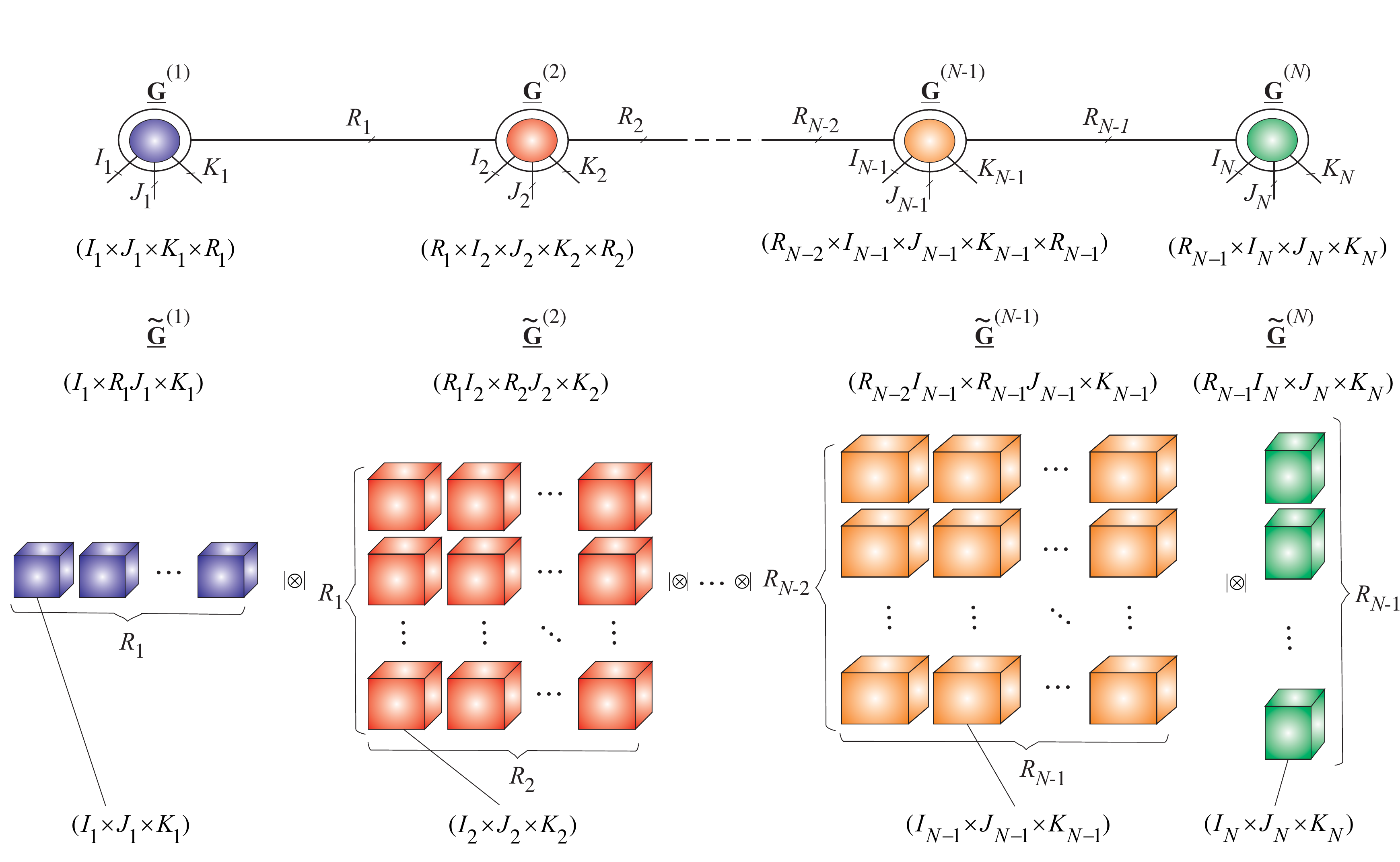}
\end{center}
\caption{(a) Example 1 of tensorization and decomposition of a large-scale 3rd-order tensor $\underline \bX \in \Real^{ I \times J \times K}$ into $3N$th-order tensor, assuming that $I =I_1 I_2 \cdots I_N$, $J =J_1 J_2 \cdots J_N $ and $K =K_1 K_2 \cdots K_N$. (b) Decomposition of the tensor via generalized Tensor Train
referred to as the Tensor Product State (TPS). The  data tensor can be expressed by the strong Kronecker product of block tensors as $\underline \bX \cong \underline {\widetilde \bG}^{(1)} \; |\otimes| \; \underline {\widetilde \bG}^{(2)} \;  |\otimes|  \cdots \otimes \; \underline {\widetilde \bG}^{(N)}  \in \Real^{I_1 \cdots I_N \times  J_1 \cdots J_N  \times  K_1 \cdots K_N}$,
where each block of the core $\underline {\widetilde \bG}^{(n)} \in \Real^{R_{n-1} I_{n} \times R_n J_{n} \times K_{n}}$ is a 3rd-order tensor of size $(I_{n} \times J_{n} \times K_{n})$, with $R_0=R_N=1$.
 The strong Kronecker product of two block cores  $\underline {\widetilde\bG}^{(n)} \in \Real^{R_{n-1} I_{n} \times R_n J_{n} \times K_{n}}$ and
 $\underline {\widetilde\bG}^{(n+1)} \in \Real^{R_n I_{n+1} \times R_{n+1} J_{n+1} \times K_{n+1}}$ is defined as the block tensor
$\underline \bC = \underline {\widetilde\bG}^{(n)} |\otimes|  \underline {\widetilde\bG}^{(n+1)} \in \Real^{R_{n-1}I_n I_{n+1} \times R_{n+1} J_n J_{n+1} \times K_n K_{n+1}}$, with 3rd-order tensor blocks $\underline \bC_{r_{n-1},r_{n+1}}=\sum_{r_n=1}^{R_n} \underline \bG^{(n)}_{r_{n-1},r_n} \otimes  \underline \bG^{(n+1)}_{r_n,r_{n+1}} \in \Real^{I_n I_{n+1} \times J_n J_{n+1} \times K_n K_{n+1}}$, where $\underline \bG^{(n)}_{r_{n-1},r_n} \in \Real^{I_{n} \times J_{n} \times K_{n}}$ and  $\underline\bG^{(n+1)}_{r_{n},r_{n+1}} \in \Real^{I_{n+1} \times J_{n+1} \times K_{n+1}}$ are block tensors of $\underline {\widetilde\bG}^{(n)}$ and $\underline {\widetilde\bG}^{(n+1)}$, respectively. In the special cases: $J=K=1$ and $K=1$ the model simplifies to standard tensor train models shown in Fig. \ref{Fig:TT=MPS=MPO} (a) and (b).} 
\label{Fig:TPS}
\end{figure*}

 For example,
the  quantization and tensorization  of a huge vector $\bx \in \Real^I$, $I = 2^K$ can be achieved through reshaping to give an
$(2 \times 2  \times \cdots \times 2)$  tensor $\underline \bX$ of order $K$,   as illustrated in Fig. \ref{Fig:Tensorization}. Such a quantized tensor $\underline \bX$ often admits low-rank matrix/tensor approximations, so that a good compression of a huge vector $\bx$ can be achieved by enforcing a maximum possible low-rank structure on the tensor $\bX$, thus admitting highly compressed  representation  via a tensor network. 

Even more generally,  an $N$th-order tensor $\underline \bX \in \Real^{I_1 \times \cdots \times I_N}$, with  $I_n=q^{K_n}$,  can be quantized  in all modes simultaneously  to yield a
$(q \times q \times \cdots  q)$  quantized tensor $\underline \bY$ of  higher-order, with small $q$.
%


 In practice, a fine ($q=2,3,4$ ) quantization is desirable to create as many virtual modes as possible, thus allowing us to implement an efficient low-rank tensor approximations.
For example, the binary encoding ($q=2$) reshapes an $N$th-order tensor with $(2^{K_1} \times 2^{K_2} \times \cdots \times 2^{K_N})$ elements into a tensor of order $(K_1+K_2 +\cdots + K_N)$, with the same number of elements.
In other words, the idea of the quantized tensor is quantization of the each $n$-th ``physical'' mode
(dimension) by replacing it with $K_n$ ``virtual'' modes, provided that the corresponding mode size $I_n$ are factorized as $I_n= I_{n,1} I_{n,2} \cdots I_{n, K_n}$. This corresponds to reshaping the $n$-th mode of size $I_n$ into $K_n$ modes of sizes $I_{n,1}, I_{n,2}, \ldots, I_{n,K_n}$.

In example shown in Fig. \ref{Fig:TPS}, the Tensor Train of huge 3rd-order tensor is expressed by the strong Kronecker products of block tensors with relatively small 3rd-order tensor blocks.
Since large-scale tensors cannot be loaded explicitly in main memory, they  usually  reside in distributed storage by splitting tensors to smaller blocks. Our approach is to apply tensor networks and represent big data by high-order tensors  not explicitly but in compressed TT formats.

The TT decomposition applied to quantized tensors is
referred to as the QTT; 
it was first introduced as a compression scheme for large-scale structured matrices, which admit low-rank TT approximation  \cite{Oseledets10}, and also developed  for  more general settings  \cite{Khoromskij-SC,QTT-Tucker,QTT-Laplace}.
The attractive property of QTT is that not only
its rank is typically small (below 10) but it is  almost  independent or at least uniformly bounded by data  size, providing a logarithmic (sub-linear) reduction of storage requirements: ${\cal{O}}(I^N) \rightarrow {\cal{O}}(N \log_q(I))$ -- so-called super-compression  \cite{Khoromskij-SC}.

Note also that, unlike in Tucker or CPD,  the TT decomposition relies on a
certain ordering of the modes so that reordering modes may
affect the numerical values of TT ranks significantly.

Quantization is  quite important for reducing the computational complexity
further, since  it allows the TT decomposition to resolve and
represent more structure in the data by splitting the ``virtual''
dimensions introduced by the quantization, as well as the
``physical'' ones. In practice it appears the most efficient to use
as fine a quantization  as possible (typically, with $q=2$) and to
generate as many virtual modes as possible.

A TT decomposition of the quantized vector is referred to as QTT
decomposition of the original vector; the ranks of this TT decomposition are called ranks
of the QTT decomposition of the original vector.


\section{\bf Mathematical and Graphical Representation of Tensor Trains}

In order to perform efficiently  various mathematical operations in the TT formats we need to represent TT decompositions in compact and easily understandable  mathematical and graphical representations \cite{Cichocki-era,NLee-Cich14}.

\subsection{\bf Vector TT/MPS Decomposition}

 \begin{figure*}[p]
 (a)\\
\begin{center}
\includegraphics[width=15.9cm,height=8.3cm]{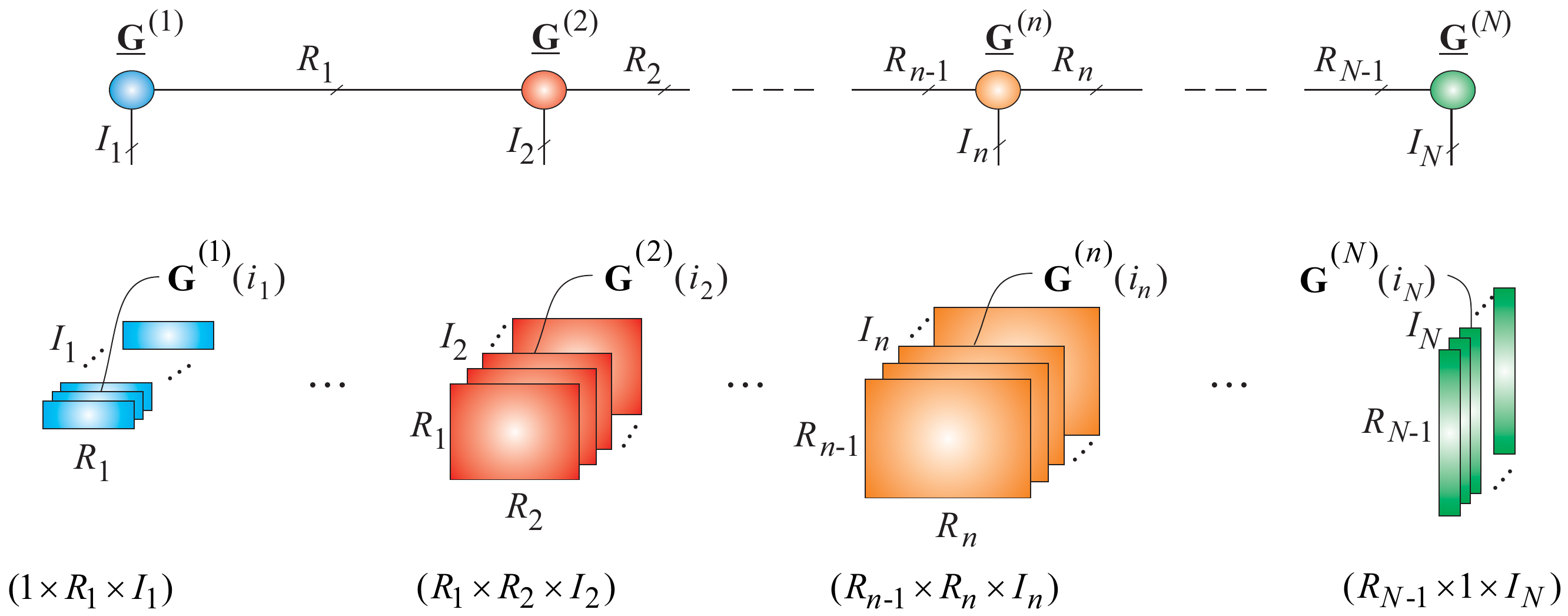}
\end{center}
(b)\\
\begin{center}
\includegraphics[width=15.9cm,height=8.1cm]{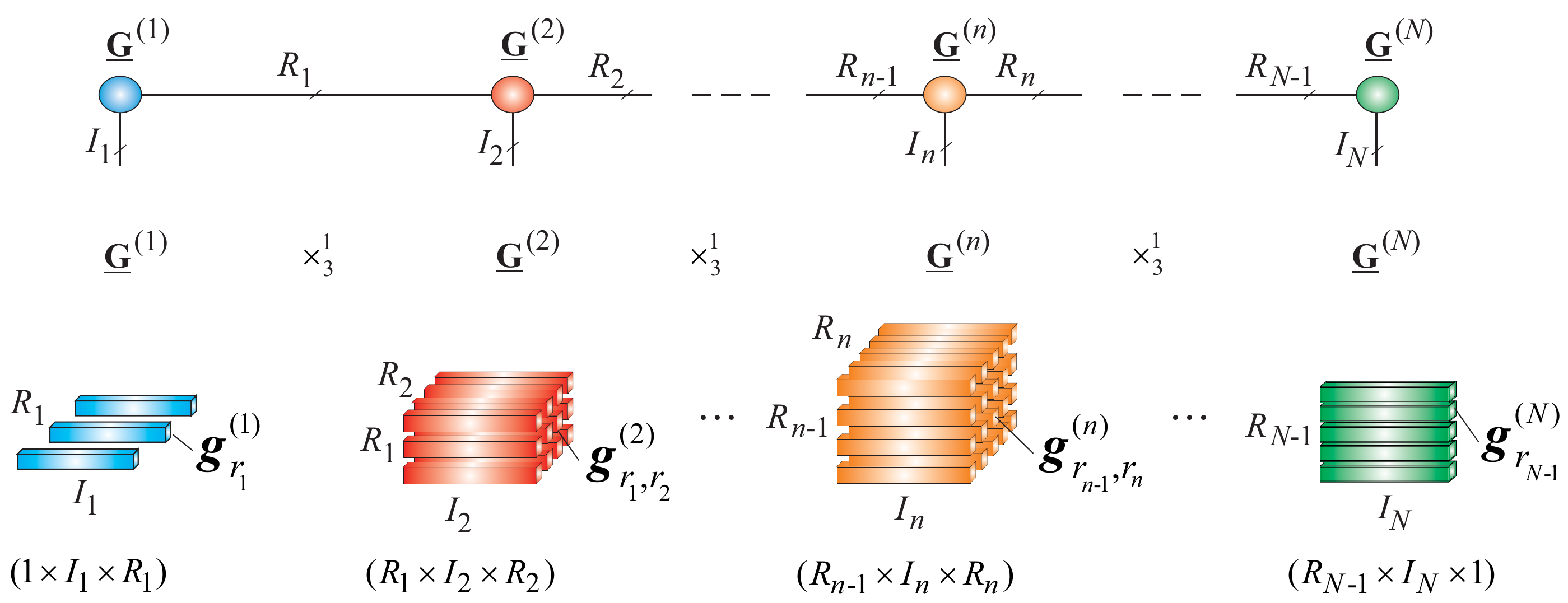}\\
\end{center}
\caption{Alternative representations of the tensor train decomposition (TT/MPS) for an $N$th-order tensor
$\underline \bX \in \Real^{I_1 \times I_2 \times I_3 \times \cdots \times I_N}$; (a)  Representation of the TT/MPS in a scalar form  via slice matrices as: $x_{i_1,i_2,\ldots,i_N} \cong
\bG^{(1)}(i_1) \;  \bG^{(2)}(i_2) \; \cdots \;  \bG^{(N)}(i_N)=
  \sum_{r_1=1}^{R_1} \sum_{r_2=1}^{R_2} \cdots \sum_{r_{N-1}=1}^{R_{N-1}}  \; g^{(1)}_{1,i_1,r_1} \; g^{(2)}_{r_1,i_2,r_2} \; g^{(3)}_{r_2,i_3,r_3} \; \cdots  g^{(N)}_{r_{N-1},i_N,1}$; (b) expressed by  the outer product of vectors (sum of rank-1 tensors) as: $\underline \bX \cong  \sum_{r_1=1}^{R_1} \sum_{r_2=1}^{R_2} \cdots \sum_{r_{N-1}=1}^{R_{N-1}} (\bg^{(1)}_{\,1,r_1} \; \circ \; \bg^{(2)}_{\,r_1, r_2} \; \circ \cdots \circ \bg^{(N-1)}_{\,r_{N-2},r_{N-1}} \; \circ \; \bg^{(N)}_{\,r_{N-1}, 1})$.  All vectors $\bg^{(n)}_{r_{n-1}r_n} \in \Real^{I_n}$ are considered to be the column vectors.}
\label{Fig:TTouter}
\end{figure*}

The (vector) tensor train (TT/MPS)  for $N$th-order data tensor $\underline \bX \in \Real^{I_1 \times I_2  \times \cdots \times I_N}$ can be described  in  the standard (tedious and rather complicated) scalar form as \cite{Oseledets_TyrTT11,OseledetsTT11}:
\be
x_{i_1,i_2,\ldots,i_N} &\cong& \sum_{r_1,r_2,\ldots,r_{N-1}=1}^{R_1,R_2, \ldots,
 R_{N-1}} g^{(1)}_{i_1,r_1} \; g^{(2)}_{r_1,i_2,r_2}  \cdots g^{(n)}_{r_{N-1},i_N}
\notag \\
\ee
or equivalently by using slice representations (see Fig. \ref{Fig:TTouter} (a)):
\be
x_{i_1,i_2,\ldots,i_N} \cong \bG^{(1)}(i_1) \; \bG^{(2)}(i_2)
 \cdots  \bG^{(N)}(i_N),
\ee
%
where slice matrices  are defined as
\begin{equation} \bG^{(n)}(i_n)= \bG^{(n)}(:, i_n, :) \in \Real^{R_{n-1} \times R_n}, \notag
\end{equation}
 i.e., $\bG^{(n)}(i_n)$ is an $i_n$th  lateral slice of the  core
 $\underline \bG^{(n)} \in \Real^{R_{n-1} \times I_n \times R_n}$ for $n=1,2,\ldots,N$, with $R_0=R_N=1$.

 However we can use several more convenient compact mathematical forms
  as follows (see Figs.  \ref{Fig:TTouter} and \ref{Fig:TT=MPS=MPO}(a) and Table \ref{table:MPS-MPO}):

\begin{enumerate}

\item  In a  tensor form using multilinear products of cores:
\be
\underline \bX &\cong& \underline \bG^{(1)}\times_3^1 \underline \bG^{(2)} \times_3^1
\cdots \times_3^1 \underline \bG^{(N-1)} \times_3^1 \underline \bG^{(N)} \nonumber \\
  &=& \llbracket \underline \bG^{(1)}, \underline \bG^{(2)},  \ldots, \underline \bG^{(N-1)}, \underline \bG^{(N)} \rrbracket,
\ee
where 3rd-order cores{\footnote{Note that the cores  $\underline \bG^{(1)}$ and $\underline \bG^{(N)}$ are now two-dimensional arrays (matrices), but to apply uniform representation, we assume
that 2nd-order cores are represented also as 3rd-order cores of mode sizes $1 \times I_1 \times R_1$ and $R_{N-1} \times I_N \times 1$, respectively.}} are defined as $\underline \bG^{(n)} \in \Real^{R_{n-1} \times I_n \times R_n}$  for $n=1,2,,\ldots,N$ with $R_0=R_N=1$ (see Fig. \ref{Fig:TT=MPS=MPO}(a)).

\item  In tensor/vector form  expressed as summation of rank-1  tensors, by using outer products of fibers  (see Fig. \ref{Fig:TTouter} (b)):
\be
\underline \bX &\cong& \sum_{r_1,r_2,\ldots,r_{N-1}=1}^{R_1,R_2, \ldots,
 R_{N-1}} \bg^{(1)}_{1,r_1} \; \circ \; \bg^{(2)}_{r_1, r_2} \; \circ \; \cdots
\circ \; \bg^{(N)}_{r_{N-1},1}, \notag \\
\label{TT-outerprod}
\ee
where
$\bg^{(n)}_{r_{n-1},r_{n}} = \underline \bG^{(n)}(r_{n-1}, \, :, \, r_n) \in \Real^{I_n}$ are mode-2 fibers, i.e., column vectors of  matrices $\bG^{(n)}_{(2)}=[\bg^{(n)}_{1,1}, \; \bg^{(n)}_{2,1},\ldots,\bg^{(n)}_{R_{n-1},1},\bg^{(n)}_{1,2},\ldots,\bg^{(n)}_{R_{n-1},R_n}] \in \Real^{I_n  \times R_{n-1} R_n}$ ($n=1,2,\ldots, N$), with $R_0=R_N=1$ or equivalently in the vector form using
the Kronecker products
\begin{equation}
\bx
\cong \sum_{r_1,r_2,\ldots,r_{N-1}=1}^{R_1,R_2, \ldots,
 R_{N-1}} \bg^{(1)}_{1,r_1} \; \otimes \; \bg^{(2)}_{r_1, r_2} \; \otimes \; \cdots
\otimes\; \bg^{(N)}_{r_{N-1},1},\\
\label{TT-Kron1}
\end{equation}
where  the vector is defined as $\bx= \bx_{\overline{i_1, i_2, \ldots, i_N}} =\mbox{vec} (\bX) \in \Real^{I_1 I_2 \cdots I_N}$.


\item  In the vector form expressed by the strong Kronecker products of block matrices \cite{KazeevT13,Cichocki-era}
(see Fig. \ref{Fig:TT=MPS=MPO}(a)):
\begin{equation}
\bx_{{\overline{i_1, i_2, \ldots, i_N}}} \cong  \widetilde \bG^{(1)} \; |\otimes| \;  \widetilde \bG^{(2)} \; |\otimes| \cdots |\otimes| \widetilde \bG^{(N)}, \\
\end{equation}
where   the cores $\underline \bG^{(n)} \in \Real^{R_{n-1} \times I_n \times R_n}$  are represented by block matrices $\widetilde \bG^{(n)} = (\bG^{(n)}_{(3)})^T \in \Real^{R_{n-1} I_n \times R_n}$ for $n=1,2,\ldots,N$, with blocks $\bg^{(n)}_{r_{n-1},r_n} \in \Real^{I_n \times 1}$, $R_0=R_N=1$, and the symbol  $|\otimes|$ denotes the strong Kronecker product.

\end{enumerate}

\begin{table*}[ht!]
\caption{Equivalent forms of the Tensor Trains (TT): MPS and MPO (with open boundary  conditions) representation of an $N$th-order tensor $\underline \bX \in \Real^{I_1 \times I_1 \times \cdots \times I_N}$ and a $2N$th-order tensor $\underline \bY \in \Real^{I_1 \times J_1 \times I_2 \times J_2 \cdots \times  I_N \times J_N }$, respectively. It is assumed that the TT rank is $\{R_1,R_2,\ldots,R_{N-1}\}$, with $R_0=R_N=1$.}
 { \shadingbox{
    \begin{tabular*}{1.01\textwidth}[t]{@{\extracolsep{\fill}}l@{\hspace{1em}}l} \hline  &  \\
{ \raisebox{0mm}[0mm][4mm]{\hspace{2.9cm}TT/MPS}}
&  \hspace{2.9cm}TT/MPO \\ \hline &   \\

 \multicolumn{2}{c}{Tensor Representations: Multilinear Products (tensor contractions)}\\
  & \\
   $ \underline \bX = \underline \bG^{(1)} \times_3^1 \; \underline \bG^{(2)} \times_3^1 \; \cdots \times_3^1 \; \underline \bG^{(N-1)} \times_3^1 \; \underline \bG^{(N)} $  &  \hspace{2em}
 $ \underline \bY = \underline \bG^{(1)} \times_4^1 \; \underline \bG^{(2)} \times_4^1  \cdots \times_4^1 \; \underline \bG^{(N-1)} \times_4^1 \; \underline \bG^{(N)}$
 \\ & \\
 $\underline \bG^{(n)} \in \Real^{ R_{n-1} \times I_n \times R_n}$, $\; (n=1,2,\ldots,N)$ &  \hspace{2em}
 $\underline \bG^{(n)} \in \Real^{ R_{n-1} \times I_n \times J_n \times R_n}$
 \\  &   \\ \hline & \\
  \multicolumn{2}{c}{Tensor Representations: Outer Products}\\
  & \\
   $ \underline \bX  = \displaystyle{\sum_{r_1,r_2,\ldots,r_{N-1}=1}^{R_1,R_2, \ldots,
 R_{N-1}}} \;\; \bg^{(1)}_{\;1,r_1} \; \circ \; \bg^{(2)}_{\;r_1, r_2}  \circ \cdots
 \circ \; \bg^{(N-1)}_{\;r_{N-2}, r_{N-1}} \; \circ \; \bg^{(N)}_{\;r_{N-1},1} $  &  \hspace{-0em}
 $  \underline \bY  = \displaystyle{\sum_{r_1,r_2,\ldots,r_{N-1}=1}^{R_1,R_2, \ldots,
 R_{N-1}}} \;\; \bG^{(1)}_{\;1,r_1} \; \circ \; \bG^{(2)}_{\;r_1, r_2}  \circ \cdots
 \circ \; \bG^{(N-1)}_{\;r_{N-2}, r_{N-1}} \; \circ \; \bG^{(N)}_{\;r_{N-1},1}$
   \\ & \\
$\bg^{(n)}_{\;r_{n-1}, r_n} \in \Real^{I_n}$ blocks of a matrix $\widetilde\bG^{(n)} =(\bG^{(n)}_{(3)})^T \in \Real^{ R_{n-1}  I_n \times  R_n}$  &  \hspace{1em}
 $\bG^{(n)}_{\;r_{n-1}, r_n} \in \Real^{I_n \times J_n} $ blocks of a matrix $\widetilde \bG^{(n)} \in \Real^{ R_{n-1} I_n  \times R_n J_n}$
 \\  &   \\ \hline & \\
  \multicolumn{2}{c}{Vector/Matrix Representations: Strong Kronecker Products}
   \\ & \\
   $ \bx_{\overline{i_1 \cdots i_N}} =  \widetilde \bG^{(1)} \; |\otimes| \;  \widetilde \bG^{(2)} \; |\otimes| \cdots |\otimes|\;  \widetilde \bG^{(N)} \in \Real^{I_1 I_2 \cdots I_N}$  &  \hspace{-3em}
 $ \bY_{\overline{i_1 \cdots i_N}; \; \overline{j_1 \cdots j_N}} =  \widetilde \bG^{(1)} \; |\otimes| \; \widetilde \bG^{(2)} \; |\otimes| \cdots \;|\otimes|\; \widetilde\bG^{(N)} \in \Real^{I_1  \cdots I_N \; \times \; J_1  \cdots J_N}$
 %
  \\ & \\
 $\widetilde\bG^{(n)} \in \Real^{ R_{n-1} I_n \times  R_n}$ a block matrix with blocks $\bg^{(n)}_{r_{n-1}, r_n} \in \Real^{I_n}$ &  \hspace{-2em}
 $\widetilde \bG^{(n)} \in \Real^{ R_{n-1} I_n  \times R_n J_n}$ a block matrix with blocks $\bG^{(n)}_{r_{n-1}, r_n} \in \Real^{I_n \times J_n}$
  \\   &  \\ \hline & \\
    \multicolumn{2}{c}{Scalar (standard) Representations}\\
  & \\
   $ x_{\;i_1,i_2,\ldots,i_N}  = \displaystyle{\sum_{r_1,r_2,\ldots,r_{N-1}=1}^{R_1,R_2, \ldots,
 R_{N-1}}} \;\;  g^{(1)}_{\;1, i_1, r_1} \; g^{(2)}_{\;r_1, i_2, r_2} \;
 g^{(3)}_{\;r_{2}, i_{3}, r_{3}}  \cdots
 g^{(N)}_{\;r_{N-1}, i_N,1}$  &  \hspace{0em}
 $y_{i_1,j_1, i_2, j_2, \ldots, i_N, j_N}  = \displaystyle{\sum_{r_1,r_2,\ldots,r_{N-1}=1}^{R_1,R_2, \ldots,
 R_{N-1}}} \;\; g^{(1)}_{\;1,i_1, j_1, r_1} \; g^{(2)}_{\;r_1, i_2, j_2, r_2} \cdots
  g^{(N)}_{\;r_{N-1}, i_{N}, j_{N}, 1}$
  \\ & \\
$g^{(n)}_{\;r_{n-1}, i_n, r_n}$ entries of a 3rd-order core $\underline \bG^{(n)} \in \Real^{ R_{n-1} \times I_n \times R_n}$ 
&
\hspace{0em}
 $g^{(n)}_{\;r_{n-1}, i_n, j_n, r_n}$ entries of a 4th-order core $\underline \bG^{(n)} \in \Real^{ R_{n-1} \times I_n \times J_n \times R_n}$
  \\  &   \\ \hline & \\
  \multicolumn{2}{c}{Slice Representations}\\
  & \\
  $ x_{\;i_1,i_2,\ldots,i_N} =  \bG^{(1)}(i_1) \; \bG^{(2)}(i_2) \cdots \bG^{(N-1)}(i_{N-1}) \; \bG^{(N)}(i_N)$  &  \hspace{0em}
 $y_{i_1, j_1, i_2, j_2,\ldots, i_N, j_N}  =  \bG^{(1)}(i_1, j_1) \; \bG^{(2)}(i_2, j_2) \cdots \bG^{(N)}(i_N, j_N)$
  \\ & \\
$\bG^{(n)}(i_n) \in \Real^{R_{n-1} \times R_n}$  lateral slices of cores $\underline \bG^{(n)} \in \Real^{ R_{n-1} \times I_n \times R_n}$
&
\hspace{0em}
$\bG^{(n)}(i_n, j_n) \in \Real^{R_{n-1} \times R_n}$ slices of cores $\underline \bG^{(n)} \in \Real^{ R_{n-1} \times I_n \times J_n \times R_n}$
\\  &  \\ \hline  
    \end{tabular*}
    }}
\label{table:MPS-MPO}
\end{table*}

In general, the strong Kronecker product of two block matrices (e.g., unfolding cores) \cite{Seberry94,Kazeev2013LRT,Cichocki-era}:
\be
{\bA} = \begin{bmatrix}\bA_{1,1}&\cdots&  \bA_{1,R_2}\\
				\vdots & \ddots & \vdots \\
				\bA_{R_{1},1}&\cdots&  \bA_{R_{1},R_2}
		\end{bmatrix} \in \Real^{R_{1} I \times R_2 J}\notag
\ee
 and\\
\be
{\bB} &=& \begin{bmatrix} \ \bB_{1,1}&\cdots& \bB_{1,R_{3}}\\
				\vdots & \ddots & \vdots \\
				 \bB_{R_{2},1}&\cdots& \bB_{R_{2},R_{3}}
		\end{bmatrix} \in \Real^{R_2 K \times R_{3} L}, \notag
 \ee
  is defined as a block matrix
\be
\bC =  {\bA}\; |\otimes| \;   {\bB} \in \Real^{R_1 I K \times R_3 J L},
\ee
with blocks $\bC_{r_{1},r_{3}}=\sum_{r_2=1}^{R_2}  \bA_{r_{1},r_2} \otimes   \bB_{r_2,r_{3}} \in \Real^{I K \times K L}$,    where $\bA_{r_{1},r_2} \in \Real^{I \times J}$ and  $\bB_{r_{2},r_{3}} \in \Real^{K \times L}$ are block matrices of ${\bA}$ and $ {\bB}$, respectively.


 The strong Kronecker product representation of a TT is
 probably the most comprehensive  and useful form for displaying
a tensor train since it allows us to perform all  operations by using
compact block matrices.


\subsection{\bf Matrix  TT (MPO) Decomposition}
\begin{figure*}[p]
(a)
\begin{center}
\includegraphics[width=14.99cm,height=10.3cm]{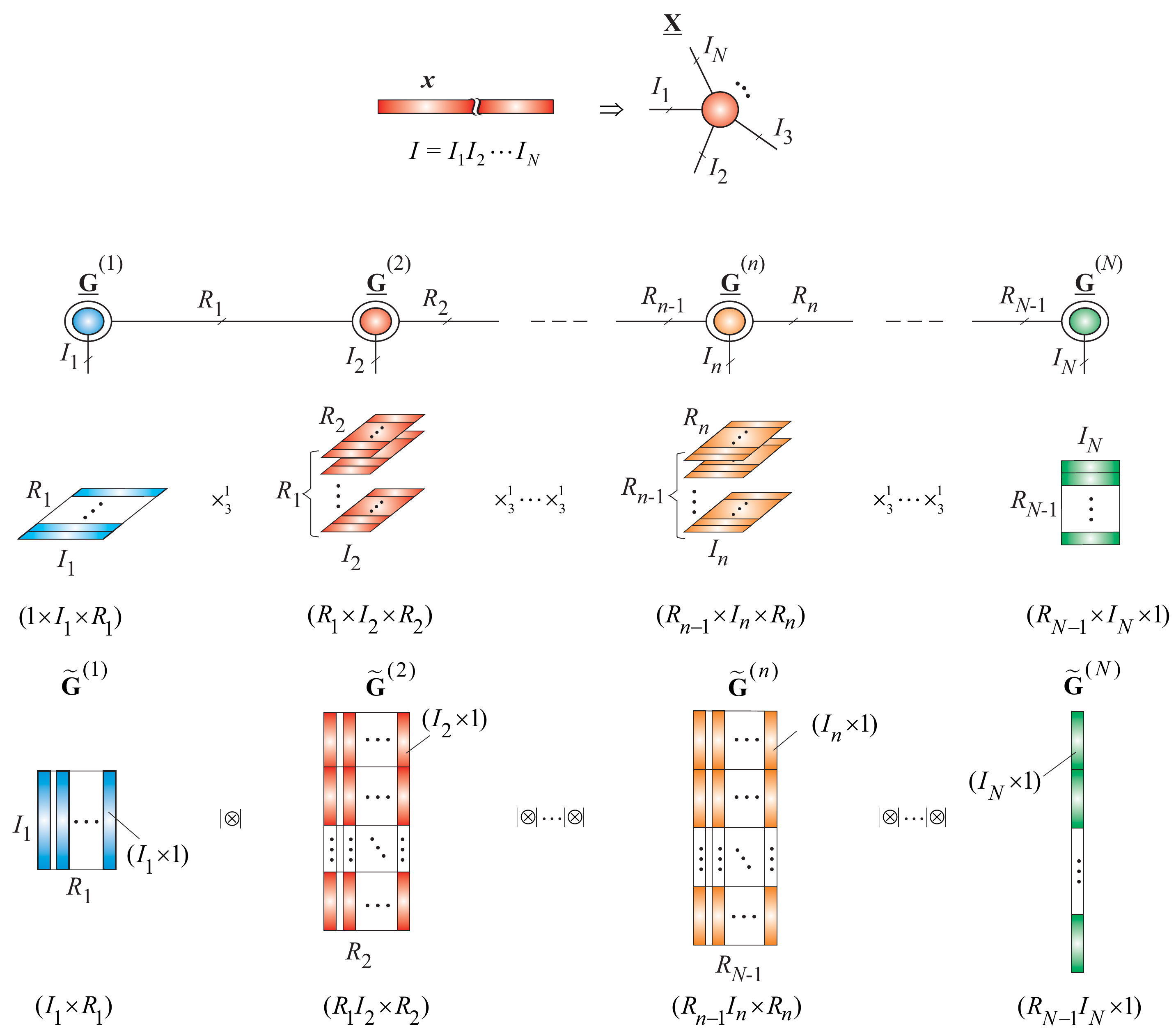}
\end{center}
(b)
\begin{center}
\includegraphics[width=5.3cm]{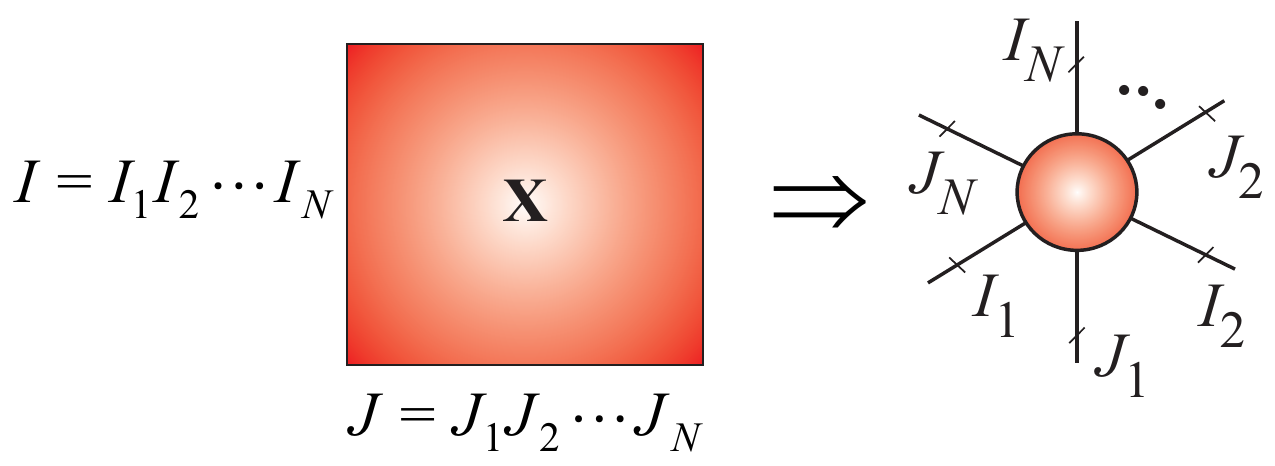}\\
\end{center}
\begin{center}
\includegraphics[width=15.8cm,height=8.9cm]{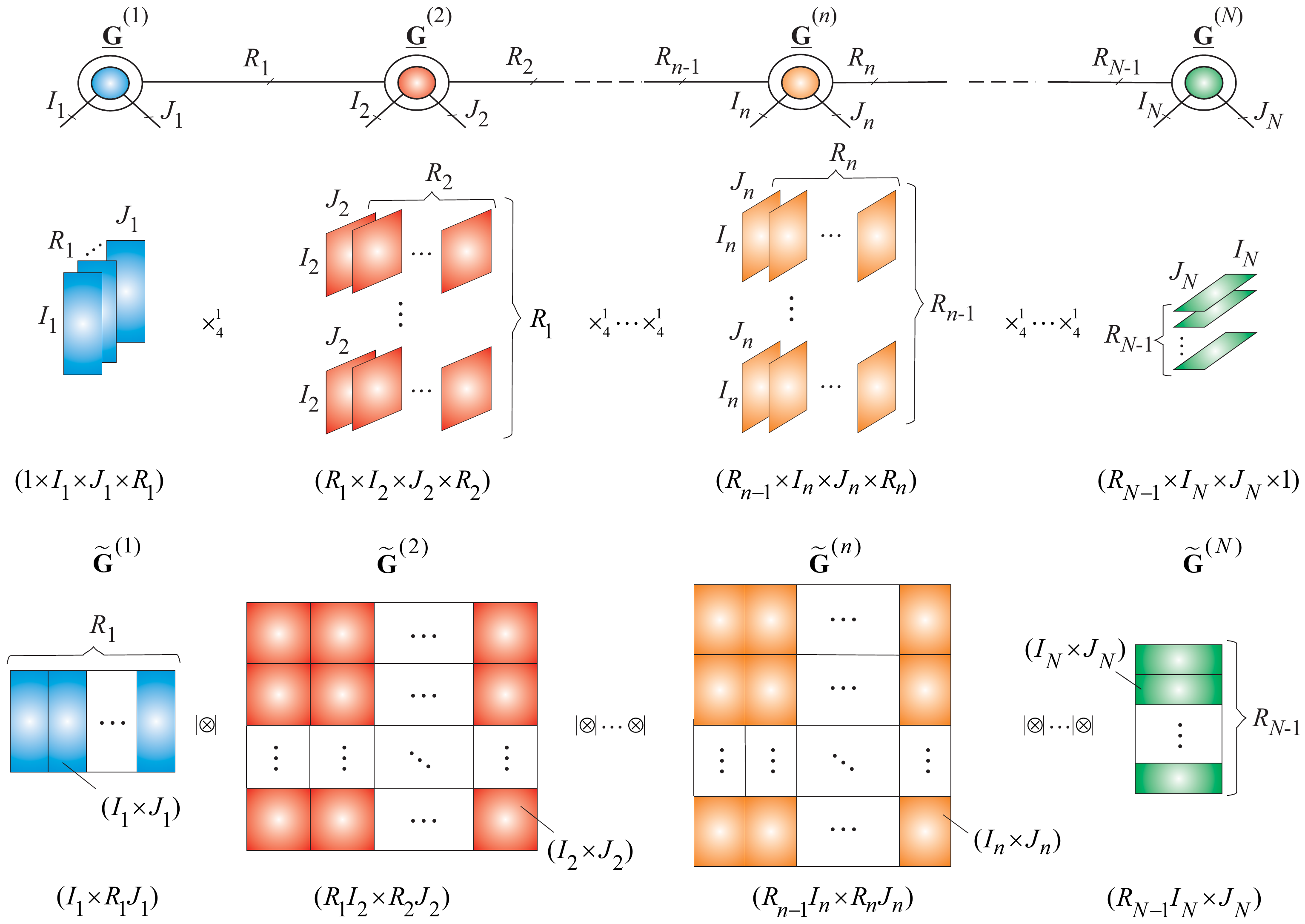}
\end{center}
\caption{Example 2 of tensorization and TT decompositions of a huge vector (a) and a  matrix (b). TT networks are represented via strong Kronecker products of block matrices.}
\label{Fig:TT=MPS=MPO}
\end{figure*}

In a similar way, we can represent a large scale matrix $\bX \in \Real^{I \times J}$, as a $2N$th-order tensor $\underline \bX \in \Real^{I_1 \times J_1 \times I_2 \times J_2 \cdots I_N \times J_N}$ with $I=I_1 I_2 \cdots I_N$ and $J=J_1 J_2 \cdots J_N$ (see Fig. \ref{Fig:TT=MPS=MPO} (b)).
This leads to an important model: the matrix TT, called also MPO (Matrix Product Operator with open boundary conditions) that  consists of chain (train) of 3rd-order and 4th-order cores{\footnote{Note that the cores  $\underline \bG^{(1)}$ and $\underline \bG^{(N)}$
 are now three-dimensional arrays, however to apply uniform representation, we assume
that 3rd-order cores are considered  also as 4th-order cores of mode sizes: $1 \times I_1 \times J_1 \times R_1$ and $R_{N-1} \times I_N \times J_N \times 1$, respectively. }} as illustrated in Fig. \ref{Fig:TT=MPS=MPO} (b). Note that the 3rd-order core tensors can be represented  as a block row and column  vectors in which each element (block) is a matrix (a lateral slice) of the cores, while  4th-order  core tensor can be interpreted  equivalently  as a  block matrix as illustrated in Fig. \ref{Fig:TT=MPS=MPO} (b).

 Since  $\bX$ is usually a full rank matrix the straightforward $2N$ dimensional exact TT  decomposition is inefficient, as it has the rank $R_n =I^N$ in the middle of a chain. Therefore, the matrix TT/MPO decompositions employ the index permutation as illustrated in Fig. \ref{Fig:TT=MPS=MPO} (b),
  and can be described  in a  scalar form as:
 \be
x_{i_1,j_1, \ldots,i_{N},j_N} &\cong& \sum_{r_1=1}^{R_1} \sum_{r_2=1}^{R_2}\cdots \sum_{r_{N-1}=1}^{R_{N-1}} g^{(1)}_{\,1,i_1,j_1,r_1} \; g^{(2)}_{r_1,i_2,j_2,r_2}  \notag \\
&\cdots& g^{(N-1)}_{r_{N-2},i_{N-1},j_{N-1},r_{N-1}} \; g^{(N)}_{r_{N-1},i_{N},j_{N},1}.
\ee
or equivalently using slice representation
\begin{equation}
x_{i_1,j_1, \ldots,i_{N},j_N} \cong  \bG^{(1)}(i_1, j_1) \; \bG^{(2)}(i_2, j_2) \cdots \bG^{(N)}(i_N, j_N), \\
\end{equation}
where $\bG^{(n)}(i_n, j_n)\cong\underline \bG^{(n)}(:, \,i_n,j_n,\,:)$ are slices of the cores
$\underline \bG^{(n)} \in \Real^{R_{n-1} \times I_n \times J_n \times R_n}$

However, the TT/MPO model  for an
$2N$th-order tensor $\underline \bX \in \Real^{I_1 \times J_1 \times \cdots I_N \times J_N}$ can be described mathematically and graphically, in  more elegant global  and  compact  forms{\footnote{i.e., not for each individual entry of a tensor.}}(see also Table \ref{table:MPS-MPO} for detailed and comparative descriptions):

A) In the tensor compact form using multilinear products
\be
\underline \bX & \cong & \underline \bG^{(1)} \; \times_4^1 \; \underline \bG^{(2)} \; \times_4^1 \; \cdots \times_4^1 \; \underline \bG^{(N)} \nonumber \\
  &=& \llbracket \underline \bG^{(1)}, \underline \bG^{(2)}, \dots, \underline \bG^{(N)} \rrbracket,
\ee
where the cores are defined as $\underline \bG^{(n)} \in \Real^{R_{n-1} \times I_{n} \times J_{n} \times R_{n}}$, with $R_0=R_N=1$, ($n=1,2,\ldots,N$).


B) In the block matrix form using the strong Kronecker products:
\begin{equation}
\bX \cong  \widetilde \bG^{(1)} \; |\otimes| \; \widetilde \bG^{(2)} \; |\otimes| \cdots \;|\otimes|\; \widetilde\bG^{(N)},  \\
\end{equation}
where $\bX = \bX_{(\overline{i_1,i_2,\ldots,i_N}\;;\;\overline{j_1,j_2,\ldots,j_N})} \in \Real^{I_1 I_2 \cdots I_{N} \times J_1 J_2 \cdots J_{N}}$ is unfolding matrix of $\underline \bX$  and $\widetilde \bG^{(n)} \in \Real^{R_{n-1} I_{n} \times R_{n} J_{n}}$ are block  matrices with blocks  $\bG^{(n)}_{r_{n-1},r_n} \in \Real^{I_{n} \times J_{n}}$ and the number of blocks $R_{n-1}\times R_{n}$.
 In the special case, when ranks of the TT/MPO $R_n=1, \; \forall n$ the strong Kronecker products simplify to the standard Kronecker products.


\section{\bf Basic Operations in TT Formats}

\begin{figure}[ht!]
(a)\\
\includegraphics[width=8.99cm]{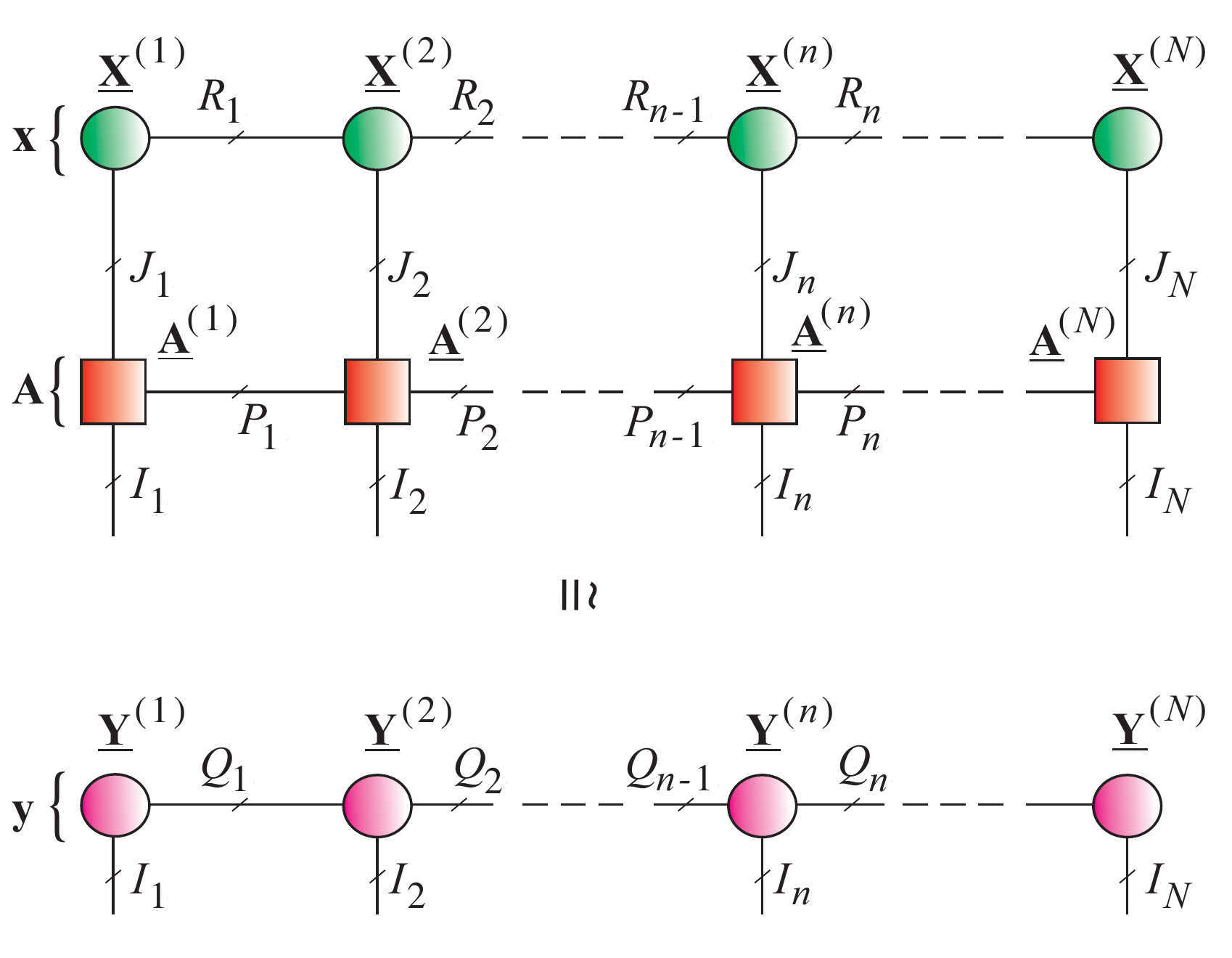}\\
(b)\\
\includegraphics[width=8.99cm]{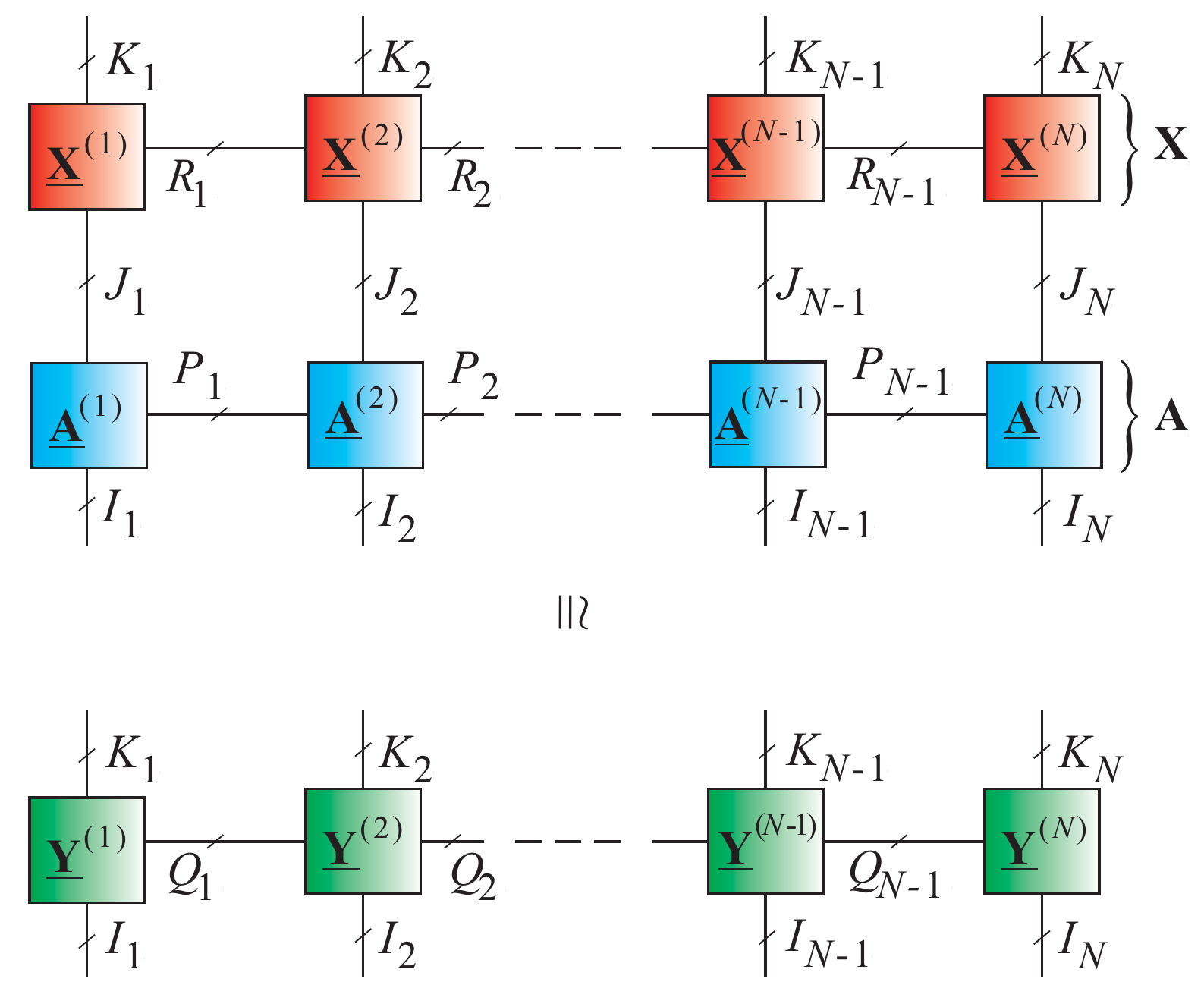}
\caption{Distributed representation of matrix equations (a) $\bA \bx =\by$ and (b)
$\bA \bX =\bY$ in TT formats.}
\label{Fig:YAX}
\end{figure}

Using the compact representations  of the TT/MPS and TT/MPO decompositions described in the previous section, we can perform easily basic mathematical operations (e.g., matrix by vector and matrix by matrix multiplications) using  block matrices.
For example, the  large-scale  matrix equation
\be
\bA \bx = \by,
\label{Axy}
\ee
where $\bA \in \Real^{I \times J}$, $\;\;\bx \in \Real^J$ and $\by \in \Real^I$  can be represented
in TT format (after performing suitable tensorization of the matrix and vectors), as shown in Fig \ref{Fig:YAX} (a), with $I=I_1 I_2 \cdots I_N$ and $J=J_1 J_2 \cdots J_N$,
and the cores defined as
\be
&&\underline \bA^{(n)} \in \Real^{P_{n-1} \times I_n \times J_n \times P_n} \notag \\
&&\underline \bX^{(n)} \in \Real^{R_{n-1} \times J_n \times R_n} \notag \\
&&\underline \bY^{(n)} \in \Real^{Q_{n-1} \times I_n  \times Q_n}. \notag
\ee
By representing the entries of the matrix $\bA$ and vectors $\bx$ and $\by$ by outer products as
\be
 \underline \bA  &=& \displaystyle{\sum_{p_1,p_2,\ldots,p_{N-1}=1}^{P_1,P_2, \ldots,
 P_{N-1}}} \;\; \bA^{(1)}_{\;1,p_1} \; \circ \; \bA^{(2)}_{\;p_1, p_2}  \circ \cdots
 \circ \; \bA^{(N)}_{\;p_{N-1},1} \notag\\
\underline \bX  &=& \displaystyle{\sum_{r_1,r_2,\ldots,r_{N-1}=1}^{R_1,R_2, \ldots,
 R_{N-1}}} \;\; \bx^{(1)}_{\;r_1} \; \circ \; \bx^{(2)}_{\;r_1, r_2}  \circ \cdots
 \circ \; \bx^{(N)}_{\;r_{N-1}}\\
 \underline \bY & =& \displaystyle{\sum_{q_1,q_2,\ldots,q_{N-1}=1}^{Q_1,Q_2, \ldots,
 Q_{N-1}}} \;\; \by^{(1)}_{\;q_1} \; \circ \; \by^{(2)}_{\;q_1, q_2}  \circ \cdots
 \circ \; \by^{(N)}_{\;q_{N-1}}, \notag
 \label{Ax=y-outer}
\ee
we can establish the following formulas:
\be
 \by^{(n)}_{q_{n-1},q_n} &=& \by^{(n)}_{\overline{r_{n-1}p_{n-1}},\;\overline{r_n p_n}} \\
&=& \bA^{(n)}_{p_{n-1},p_n}  \; \bx^{(n)}_{r_{n-1},r_n}
 \in \Real^{I_n},  \nonumber
\ee
with $Q_n=P_n R_n$ for $n=1,2,\ldots,N$.

On the other hand, by representing the  matrix $\bA$ and vectors $\bx$, $\by$ via the strong Kronecker products:
\be
\bA &=& \tilde \bA^{(1)} |\otimes| \tilde \bA^{(2)} |\otimes| \cdots |\otimes| \tilde \bA^{(N)} \notag\\
\bx &=& \tilde \bX^{(1)} |\otimes| \tilde \bX^{(2)} |\otimes| \cdots |\otimes| \tilde \bX^{(N)}\\
\by &=& \tilde \bY^{(1)} |\otimes| \tilde \bY^{(2)} |\otimes| \cdots |\otimes| \tilde \bY^{(N)}, \notag
\ee
with $\tilde \bA^{(n)} \in \Real^{P_{n-1}  I_n \times J_n P_n}, \;\;
\tilde \bX^{(n)} \in \Real^{R_{n-1}  J_n \times R_n}$ and $
\tilde \bY^{(n)} \in \Real^{Q_{n-1} I_n  \times Q_n}$,
we can easily establish a simple  relationship
\be
\tilde \bY^{(n)} &=&\tilde \bA^{(n)} |\cdot| \tilde \bX^{(n)} \in
\Real^{R_{n-1} P_{n-1} I_n \times R_{n} P_{n}}, \notag\\
&& n=1,2,\ldots,N,
\label{ACproduct}
\ee
where operator $|\cdot|$ means the AC product of two block matrices.

\begin{figure}[ht]
\begin{center}
\includegraphics[width=8.79cm]{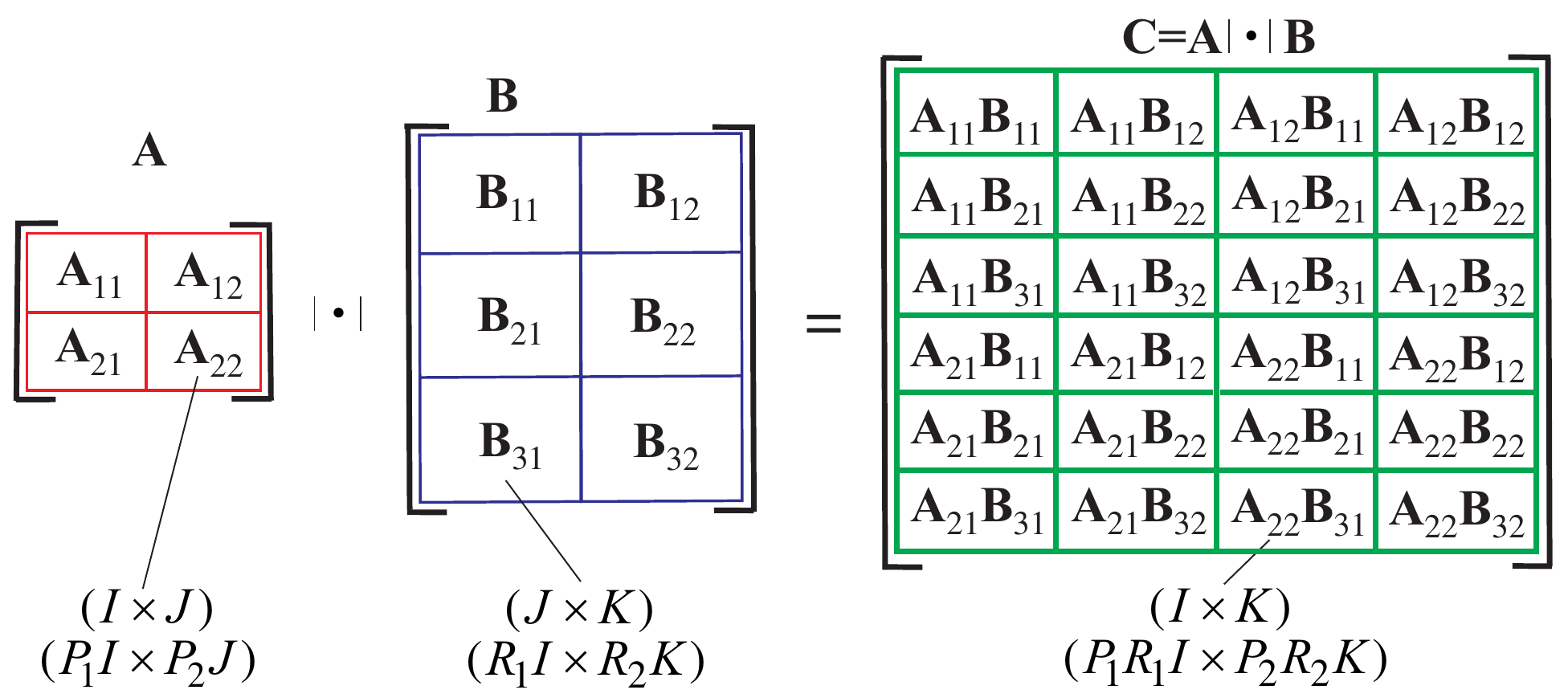}
\end{center}
\caption{Graphical illustration of the AC product for two block matrices.}
\label{Fig:SKPAC-product}
\end{figure}

In general, the  AC product of a block  matrix $\bA^{(n)} \in \Real^{P_{n-1} I_n \times P_n J_n}$
 (with blocks $\bA^{(n)}_{p_{n-1},p_n} \in \Real^{I_n \times J_n}$) and a block matrix $\bB^{(n)} \in \Real^{R_{n-1} J_n \times R_n K_n}$
 (with  blocks $\bB^{(n)}_{r_{n-1},r_n} \in \Real^{J_n \times K_n}$)  is defined as a block matrix
 $\bC^{(n)}= \bA^{(n)} |\cdot| \bB^{(n)} \in \Real^{Q_{n-1} I_n \times Q_n K_n}$
 (with blocks $\bC^{(n)}_{q_{n-1},q_n} = \bA^{(n)}_{p_{n-1},p_n} \bB^{(n)}_{r_{n-1},r_n} \in \Real^{I_n \times K_n}$ as illustrated in Fig. \ref{Fig:SKPAC-product}. 

The AC product of two block matrices is similar to the  Tracy-Singh
product but  the Kronecker product for block matrices is replaced by  the ordinary products
matrix-by-matrix.


In a similar way, we  can  represent in TT format a matrix equation
 \be
\bY \cong \bA \bX=\bA \bB^T,
\label{AXY}
\ee
where $\bA \in \Real^{I \times J}$, $\;\;\bX = \bB^T \in \Real^{J \times K}$ and $\bY \in \Real^{I \times K}$  as shown in Fig \ref{Fig:YAX} (b), with $I=I_1 I_2 \cdots I_N$, $J=J_1 J_2 \cdots J_N$ and
$K=K_1 K_2 \cdots K_N$
with the cores defined as
\be
&&\underline \bA^{(n)} \in \Real^{P_{n-1} \times I_n \times P_n \times J_n} \notag \\
&&\underline \bX^{(n)} \in \Real^{R_{n-1} \times J_n \times R_n \times K_n} \notag \\
&&\underline \bY^{(n)} \in \Real^{Q_{n-1} \times I_n \times Q_n  \times K_n}. \notag
\ee
%
%

It can be proved that by assuming that matrices: $\bA \in \Real^{I \times J}$ and $\bX \in \Real^{J \times K}$ are represented in TT formats and expressed via the strong Kronecker product of block matrices as:
$\bA =  \tilde \bA^{(1)} |\otimes| \tilde \bA^{(2)} |\otimes| \cdots |\otimes| \tilde \bA^{(N)}$ and
$\bX = \tilde \bX^{(1)} |\otimes|  \tilde \bX^{(2)} |\otimes| \cdots |\otimes| \tilde \bX^{(N)}$, with $\tilde \bA^{(n)} \in \Real^{P_{n-1} I_n \times J_n P_n}$ and $\tilde \bX^{(n)} \in \Real^{R_{n-1} J_n \times K_n R_n}$, respectively,
  then the matrix $\bY = \bA \bX$  can be  expressed in TT format via the strong Kronecker products:
  $\bY =  \tilde \bY^{(1)} |\otimes| \tilde \bY^{(2)} |\otimes| \cdots |\otimes| \tilde \bY^{(N)}$, where
  $\tilde \bY^{(n)} = \tilde \bA^{(n)} |\cdot| \tilde \bX^{(n)} \in \Real^{Q_{n-1} I_n \times K_n Q_n}$, $(n=1,2, \ldots, N)$, with blocks $\tilde \bY_{q_{n-1},q_n}^{(n)} = \tilde \bA^{(n)}_{p_{n-1},p_n}   \tilde \bX^{(n)}_{r_{n-1},r_n} $, where $Q_n=R_n P_n,\;\; \forall n$.\\


The above operation assumes precise contraction of cores. However,  an exact contraction of core tensors for very large scale data is impossible, and the choice of the approximating procedure determines the
efficiency and accuracy of  algorithms implemented for specific computational or optimization problems \cite{ContracPEPS14}.
In other words, contraction operations as matrix-by-vector or matrix-by-matrix products TT ranks grows and the TT ranks could become excessively large and therefore truncation  (called also recompression) or low-rank matrix approximations are needed. In the truncation procedure (usually, performed via QR/SVD or CUR) the core tensors $\underline \bG^{(n)}$ are approximated by other core tensors with minimal possible TT-ranks with desired prescribed accuracy \cite{OseledetsTT11}.

\section{\bf Tensor Train (TT/MPS) Splitting}

\begin{figure}[t]
(a)\\
\includegraphics[width=8.99cm]{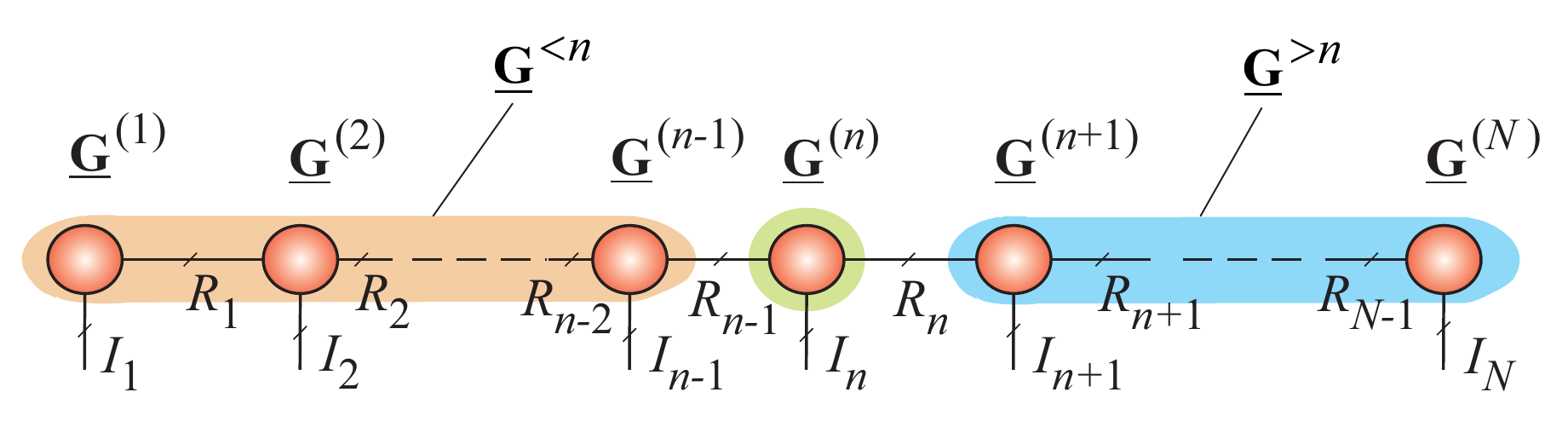}\\
(b)\\
\includegraphics[width=9.1cm]{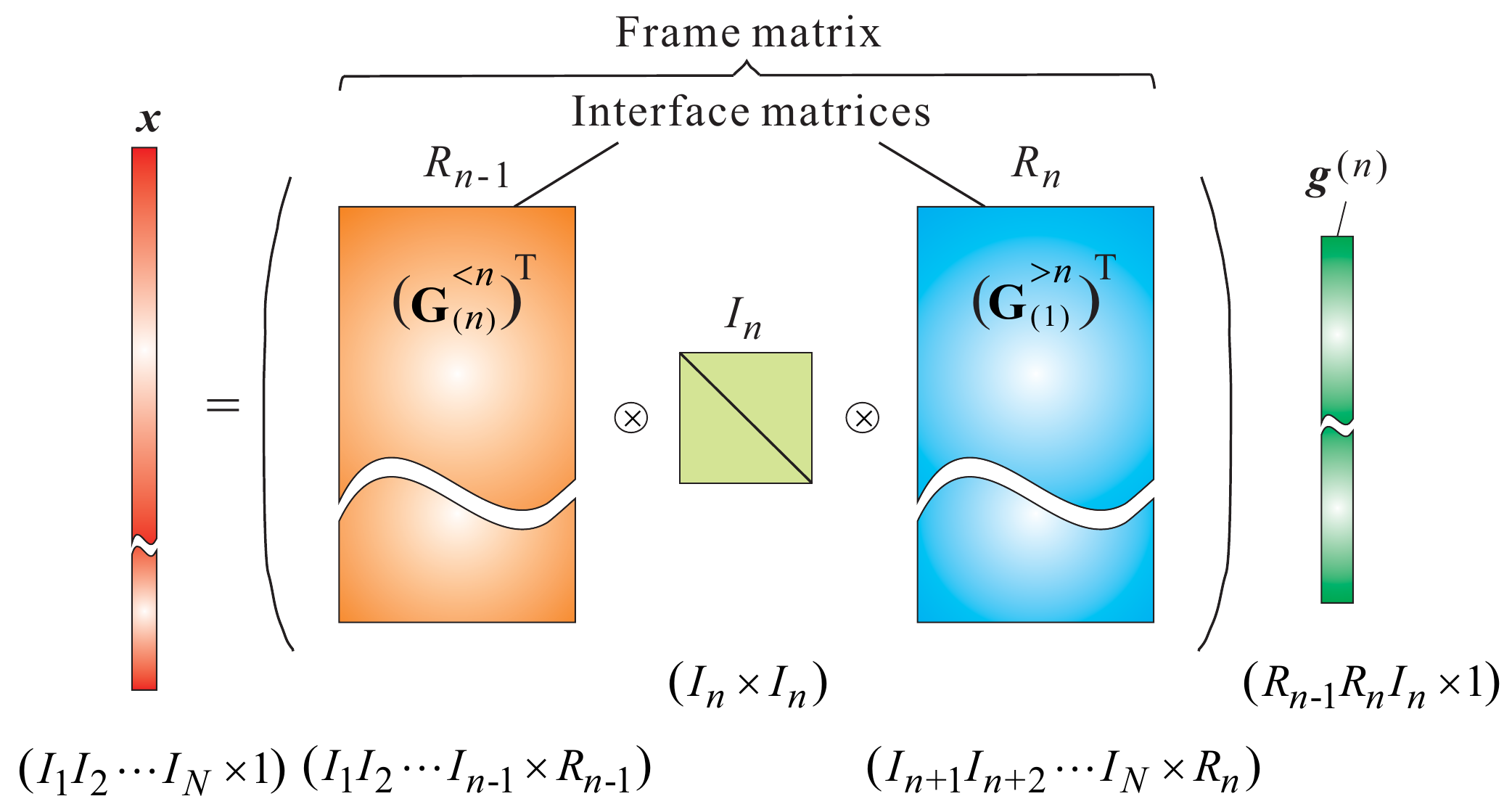}
\caption{Extraction of a single core. (a) Graphical representation and notations of vector tensor train and subtrains. (b) Graphical illustration of the equation expressed via interface matrices or a frame matrix (see Eqs. (\ref{frame-eq}) and (\ref{frame-matrices})).}
\label{Fig:TTsplitt}
\end{figure}

In practical applications it is very useful and efficient  to divide a  TT decomposition, representing a tensor $\underline \bX = \llbracket \underline \bG^{(1)}, \underline \bG^{(2)}, \ldots, \underline  \bG^{(N)}\rrbracket \in \Real^{I_1 \times I_2 \times \cdots \times I_N}$, into subtrains
as illustrated in Fig. \ref{Fig:TTsplitt}.

\subsection{\bf Extraction of a single core}
For this purpose, we define subtrains as follows
\begin{equation}
\underline \bG^{<n}=\llbracket \underline \bG^{(1)}, \underline \bG^{(2)}, \ldots, \underline \bG^{(n-1)}\rrbracket \in \Real^{I_1 \times I_2 \times \cdots \times I_{n-1} \times R_{n-1}} \\ \label{TT-splitt1}
\end{equation}
\begin{equation}
\underline \bG^{>n} = \llbracket \underline \bG^{(n+1)}, \underline \bG^{(n+2)}, \ldots, \underline \bG^{(N)}\rrbracket \in \Real^{R_n \times I_{n+1}  \times \cdots \times I_{N}} \\
\label{TT-splitt2}
\end{equation}
with corresponding unfolding matrices called interface matrices:
\be
&&\bG^{<n}_{(n)} \in \Real^{R_{n-1} \times I_1 I_2  \cdots  I_{n-1}}  \\
&&\bG^{>n}_{(1)}  \in \Real^{R_n \times I_{n+1}   \cdots I_{N}}
\label{interface-matrices}
\ee

Using basic multilinear algebra, we  can construct a set of linear  equations referred to as
the frame equation:
\be
\bx = \bG_{\neq\,n} \; \bg^{(n)}, \qquad n=1,2,\ldots, N,
\label{frame-eq}
\ee
where $\bx= vec(\underline \bX) \in \Real^{I_1 I_2 \cdots I_N}$, $\;\;\bg^{(n)} =vec(\underline \bG^{(n)}) \in \Real^{R_{n-1} I_n R_n}$ and
a tall-and-skinny matrix, called the frame matrix,  formulated as
\begin{equation}
 \bG_{\neq \, n} = (\bG^{<n}_{(n)})^T \otimes \bI_{I_n} \otimes (\bG^{>n}_{(1)})^T \in \Real^{I_1 I_2  \cdots  I_{N} \times R_{n-1} I_n R_n}.   \\
\label{frame-matrices}
\end{equation}
The frame and interface matrices help to show a very important property of TT,
 -- TT is linear with respect to each core $\bG^{(n)}$ in the vectorized form \cite{OseledetsTT11}.\\

\subsection{\bf Extraction of two cores for two--sided DMRG}

\begin{figure}[t]
(a)\\
\includegraphics[width=8.99cm]{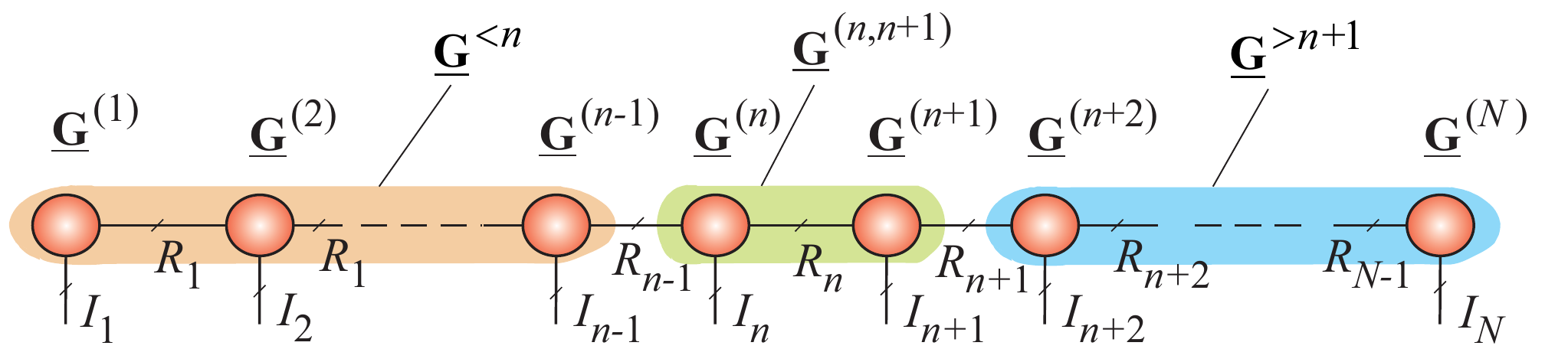}\\
(b)\\
\includegraphics[width=9.1cm]{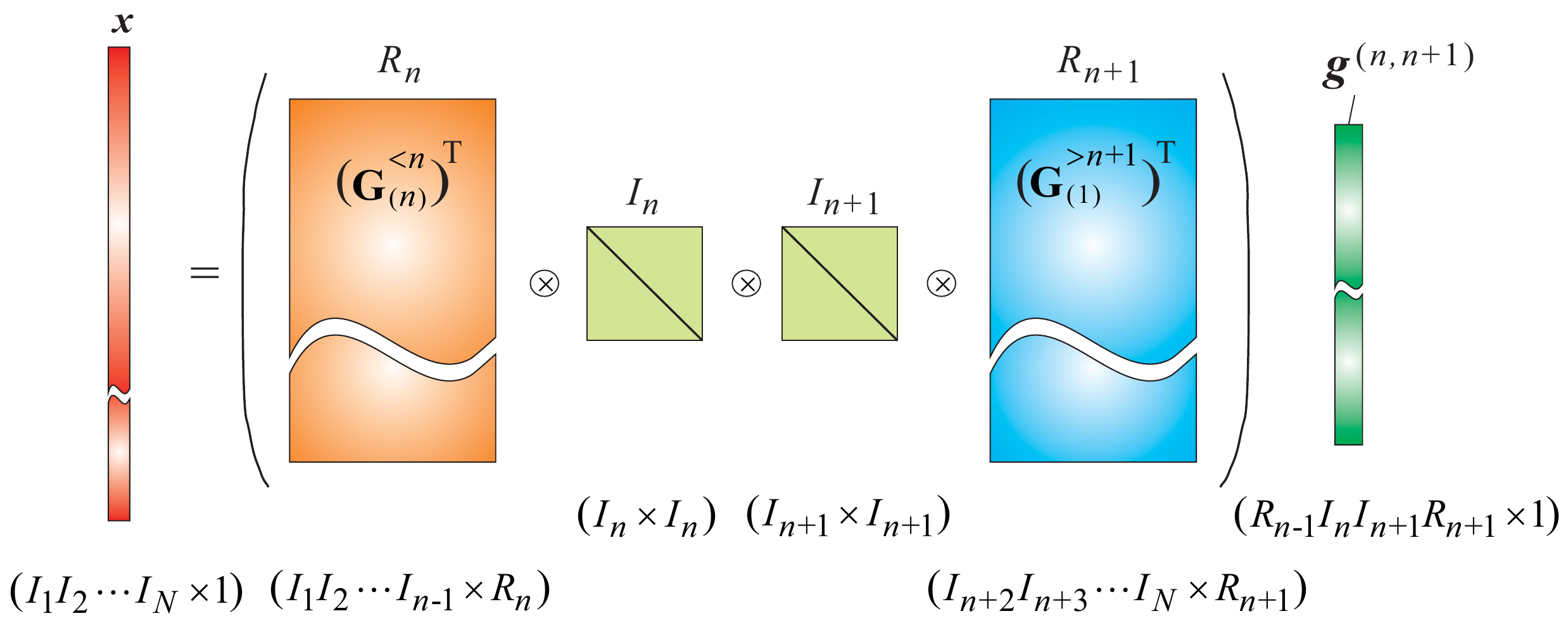}
\caption{Extraction of two cores: (a) Graphical representation of the tensor train and subtrains. (b) Graphical illustration of the frame equation (see  Eqs. (\ref{frame-eq2}) and (\ref{frame-matrices2})).}
\label{Fig:TTsplitt2}
\end{figure}

In a similar way, we can formulate equations for the 2-sided DMRG, where we extract block of two consecutive cores (see Fig. \ref{Fig:TTsplitt2}):
\begin{equation}
\bx =\bG_{\neq \, n,n+1} \; \bg^{(n,n+1)}, \;\; n=1,2,\ldots, N-1,
\\
\label{frame-eq2}
\end{equation}
where the frame (tall-and-skinny) matrix is formulated as
\be
\bG_{\neq \, n,n+1} &=& (\bG^{<n}_{(n)})^T \; \otimes \; \bI_{I_n} \; \otimes \bI_{I_{n+1}}  \otimes \; (\bG^{>n+1}_{(1)})^T  \notag \\
&\in& \Real^{I_1 I_2  \cdots  I_{N} \times R_{n-1} I_n  I_{n+1} R_{n+1}}
\label{frame-matrices2}
\ee
and  $\bg^{(n,n+1)} =vec[\bG^{(n)\;T}_{(3)} \bG^{(n+1)}_{(1)}]=vec(\underline \bG^{(n,n+1)}) \in \Real^{R_{n-1} I_n I_{n+1} R_{n+1}}$ for $n=1,2,\ldots,N-1$.

Simple matrix manipulations give the following useful relationships \cite{NLee-Cich14}:
\be
\bX_{(n)} &=& \bG^{(n)}_{(n)} \; ( \bG^{<n}_{(n)} \otimes \bG_{(1)}^{>n}) ,\\
\bG_{\neq \;n+1} &=& \bG_{\neq \;n,n+1} \; (\bI_{R_{n+1} I_{n+1}} \otimes \; \bG^{(n)\;T}_{(3)}).
\label{frame-matrices12}
\ee
If  cores  are normalized in a such way that all cores to the left of the currently considered
(optimized) core $\underline \bG^{(n)}$ are left-orthogonal:
 \be \bG_{(3)}^{(k)} \bG_{(3)}^{(k)\;T} =\bI_{R_k}, \quad k<n,
 \ee
 and all cores to the right of the $\underline \bG^{(n)}$ are right-orthogonal, i.e.:
\be
\bG_{(1)}^{(p)} \bG_{(1)}^{(p)\; T} =\bI_{R_{p-1}}, \quad p>n,
\ee
 then the frames matrices have orthogonal columns \cite{NRG-DMRG12,dolgovEIG2013,KressnerEIG2014}:
\be
\bG_{\neq\, n}^T \; \bG_{\neq \; n}&=&\bI_{R_{n-1} I_n R_n}, \\
\bG_{\neq \, n,n+1}^T \; \bG_{\neq\,n,n+1}&=&\bI_{R_{n-1} I_n I_{n+1} R_{n+1}}.
\label{orth-frames}
\ee
Orthogonalization of cores is usually performed by the QR/SVD algorithm \cite{OseledetsTT11}
(see  Section \ref{sec:Applic} for detail).

\section{\bf Application of TT Decompositions to Large-Scale Optimization Problems}
\label{sec:Applic}

For extremely large-scale problems, due to curse of dimensionality, most  computation and optimization problems (such as solving eigenvalue problems, SVD,  sparse PCA, Canonical Correlation Analysis (CCA),   system of linear  equations) are  intractable when using standard numerical methods.

Our goal and objective is to seek for alternative solutions for specific optimization problems in approximative  tensor compressed  formats.
 The key idea  discussed in this section is to represent huge data in TT formats
  and to apply some kind of  separation of variables \cite{dolgovEIG2013,KressnerEIG2014,Holtz-TT-2012}.
In other words, we approximate involved vectors and matrices by suitable TT networks and  convert  a large-scale specific optimization problem  into a set of much smaller optimization problems.

We next illustrate this approach  by considering several fundamental optimization problems
for very large-scale data.

\subsection{\bf Computing a Few Extreme Eigenvalues and  Eigenvectors  for Symmetric EVD
in TT Format}

In many applications we need to compute  extreme (minimum or maximum) eigenvalues
and corresponding eigenvectors of a huge structured symmetric matrix.
The basic problem we try to solve is the standard symmetric
eigenvalue decomposition (EVD), which can be formulated as
\be
\bA\; \bx_k = \lambda_k \bx_k,   \qquad  (k=1,2,\ldots,K),
\label{eq:eigp} \ee
where $\bx_k \in \Real^I$ are the orthonormal eigenvectors, $\lambda_k$
are the corresponding eigenvalues of a symmetric matrix  $\bA \in \Real^{I \times I}$ (e.g., a positive-definite covariance matrix of zero-mean signals $\by(t)$).
Note that (\ref{eq:eigp}) can be
written in the matrix form as
\be
\bX^T \, \bA \, \bX = \mbi \Lambda_K,
\ee
 where $\mbi\Lambda_K$ is the diagonal matrix
of K smallest or largest eigenvalues (ranked in  ascending  or descending order, respectively).\\
%

\subsubsection{\bf Tensor Network for Computing Single Eigenvalue and Corresponding  Eigenvector}

Many iterative algorithms for  extreme eigenvalue and the corresponding eignevector
exploit the Rayleigh quotient (RQ) of the  symmetric  matrix
 as a cost function.
The Rayleigh quotient $R(\bx)$ is defined for $\bx \neq \0$, as
\be
J(\bx) = R(\bx, \bA) = \frac{\bx^T \bA \bx}{\bx^T \bx}= \frac{\langle \bA \bx, \bx\rangle}{\langle\bx,\bx\rangle},
\ee
where
\be
\lambda_{max} = \max R( \bx, \bA ),\quad
 \lambda_{min} = \min  R( \bx, \bA ),
\ee
where $\lambda_{max}$ and $\lambda_{min}$ denote respectively  largest and  smallest
eigenvalue of the  matrix $\bA$.
More generally, the critical points and critical values of
$R(\bx, \bA)$ are the eigenvectors and eigenvalues of $\bA$.

\begin{figure}[t]
\begin{center}
\includegraphics[width=8.99cm]{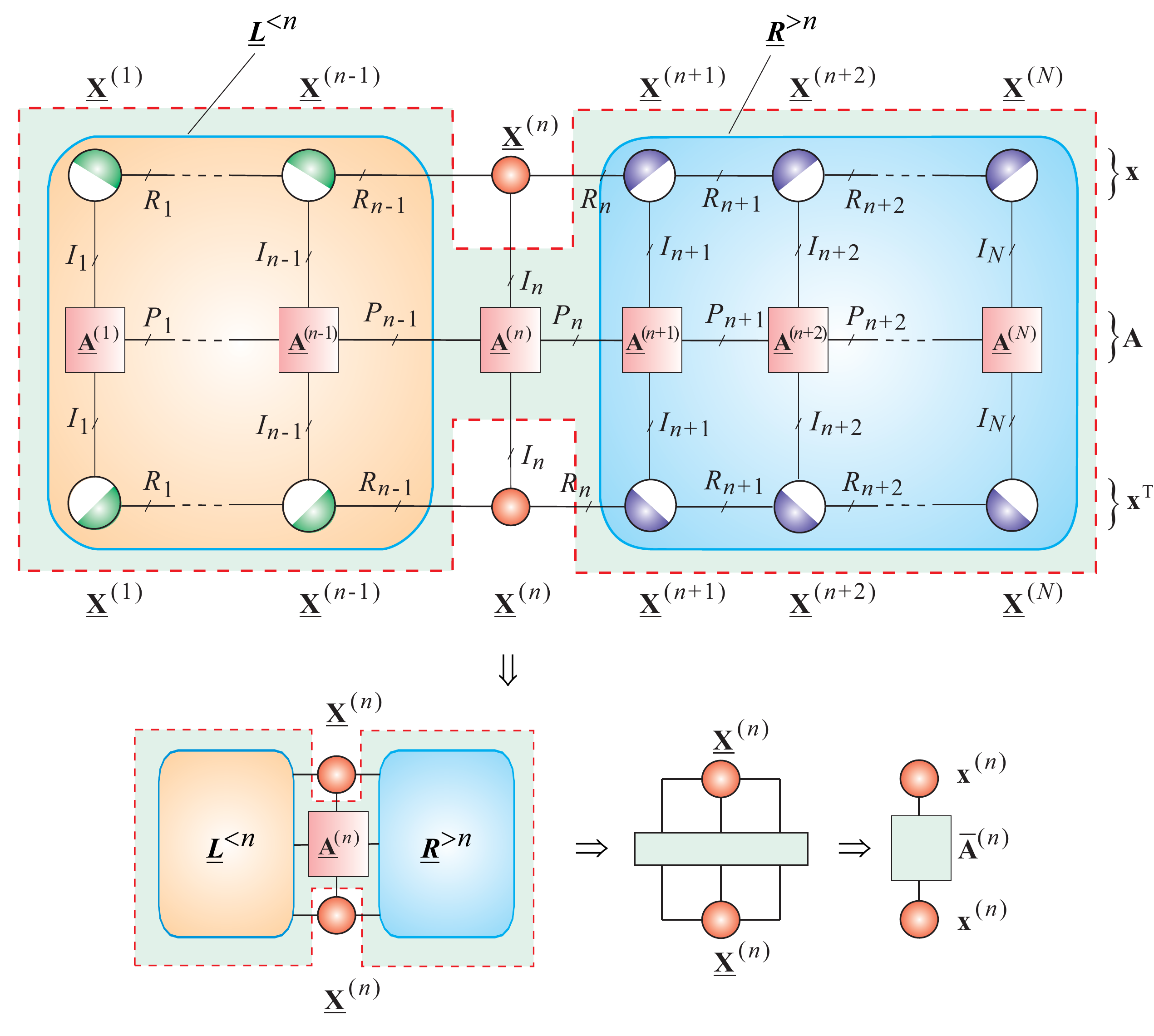}
\end{center}
%
\caption{Computation of a single extreme eigenvalue and the corresponding eigenvector
$\bx \in \Real^I$ in the TT format for a symmetric  matrix $\bA \in \Real^{I \times I}$. The frame matrix  maps a TT core  into a large vector.  The  tensor network corresponds to the Rayleigh quotient,  with the matrix $\bA$ and vectors $\bx \in \Real^I$ given in the tensor train format with distributed  indices $I=I_1 I_2 \cdots I_N$. The cores included in the  shaded areas form the matrix $\bar \bA^{(n)}$ (the effective Hamiltonian), which can computed by sequential core contractions.}
\label{Fig:EIG1}
\end{figure}

If the matrix $\bA$  admits low-rank TT approximation,
 we can convert large-scale problems into smaller optimization problems
by representing the eigenvector $\bx$ and the matrix $\bA$ in TT (MPO/MPS) formats (see also Fig. \ref{Fig:EIG1}) as:
\be
\underline \bA &=& \llbracket \underline\bA^{(1)}, \underline\bA^{(2)}, \ldots, \underline\bA^{(N)} \rrbracket \notag \\
\underline \bX &=& \llbracket \underline\bX^{(1)}, \underline\bX^{(2)}, \ldots,\underline \bX^{(N)} \rrbracket
\ee
and by computing iteratively the frame equation $\bx =\bX_{\neq n} \; \bx^{(n)}$, $(n=1,2,\ldots,N)$, with
the frame  matrices:
\be
\bX_{\neq n} = (\bX^{<n}_{(n)})^T \otimes \bI_{I_n} \otimes (\bX^{>n}_{(1)})^T \in \Real^{I_1 I_2  \cdots  I_{N} \times R_{n-1} I_n R_n}. \notag
\label{frame-matricesX}
\ee
Assuming that cores $\underline \bX^{(n)}$ are constrained to be left and right orthogonal, we can minimize (or maximize) the RQ as follows:
\be
\min_{\bx} J(\bx) &=& \min_{\bx^{(n)}} J(\bX_{\neq n} \bx^{(n)}) \\
&=& \min_{\bx^{(n)}} \frac{\langle \bar \bA^{(n)} \bx^{(n)}, \bx^{(n)}\rangle}{\langle\bx^{(n)},\bx^{(n)}\rangle}, \quad n=1,2,\ldots, N, \nonumber
\label{RQframe}
\ee
where $\bx^{(n)} =vec(\underline\bX^{(n)}) \in \Real^{R_{n-1} I_n R_n}$ and the matrix $\bar \bA$,
often called  the effective Hamiltonian, and can be expressed as
\be
\bar \bA^{(n)} = (\bX_{\neq n})^T \bA \bX_{\neq n} \in \Real^{R_{n-1}I_n R_n \times R_{n-1}I_n R_n}
\label{frameAxn}
\ee
for $n=1,2,\ldots,N$.

\begin{figure}[ht]
(a)
\begin{center}
 \includegraphics[width=8.2cm]{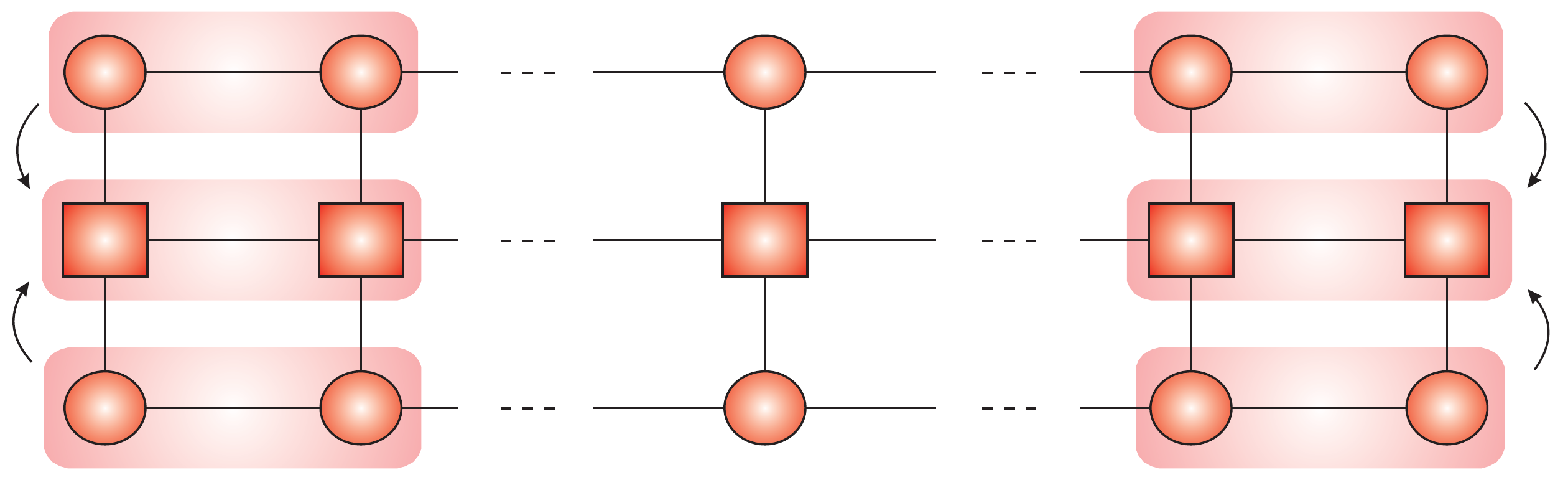}\\
 \end{center}
(b)
\begin{center}
  \includegraphics[width=8.2cm]{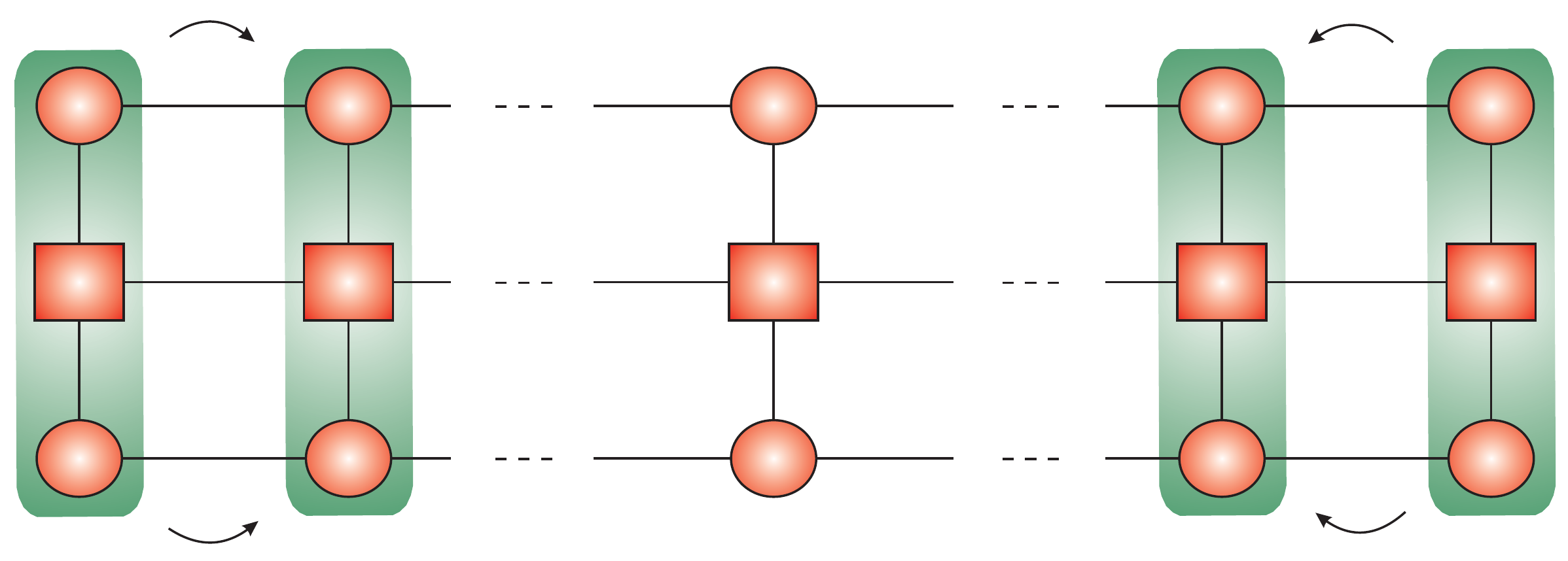}
  \end{center}
\caption{(a) Non-optimal (inefficient) and (b) optimal (efficient) contraction of the
TT (MPS/MPO) network.}
\label{Fig:TTcontrac}
\end{figure}

Note that the matrices $\bar \bA^{(n)}$  are usually much smaller than the original matrix $\bA$ if the TT rank is relatively small, then, the large-scale optimization problem can be converted into a much smaller
set of  EVDs, i.e., by solving the  set of equations:
\be
\bar \bA^{(n)} \bx^{(n)} = \lambda \bx^{(n)}, \quad n=1,2,\ldots,N.
\ee
In practice, we never compute the matrices $\bar \bA^{(n)}$ directly by Eq. (\ref{frameAxn}),
 but iteratively via optimized  and approximative contraction of  cores of the tensor network as  shown
  in Fig. \ref{Fig:TTcontrac}.

This is achieved by sweeping through the tensor network in a recursive  forward and backward manner and forth through each node.
An initial guess for all cores $\underline \bX^{(n)}$ is first made, and  then we sweep
through the set of the cores with the index $n$,  keeping all other cores fixed and choosing the   $\underline \bX^{(n)}$,
such that the cost  function  gradually  decreases. By repeating such sweeps (from the left to the right and from the right to the left)
through the tensor network  several times that usually  leads to a converged  approximation.
Note that this sweeping process works in a similar fashion as a self-consistent
recursive loops, where we iteratively and gradually improve the solution.

In order to efficiently estimate the matrix $\bar \bA^{(n)}$, we need to compute
 blocks $\underline \bL^{<n}$ and  $\underline \bR^{>n}$ (see Fig. \ref{Fig:EIG1}).
However, $\underline \bL^{<n}$ and  $\underline \bR^{>n}$ can be built iteratively in order to best  reuse available information;
this involves an optimal arrangement of a tensor  network contraction.
In a practical implementation of the algorithm the full
network contraction is never carried out globally, but we rather look at
blocks $\underline \bL^{n}$ and $\underline \bR^{>n}$ that are growing and shrinking in size sweeping along the tensor network \cite{schollwock11-DMRG}.
In other words, the construction of blocks $\underline \bL^{<n}$ and  $\underline \bR^{>n}$
is an iterative process  in a way that directly matches block growth and shrinkage.
If we sweep through the chain from right to left or vice-versa we can build up $\underline \bL^{<n}$
and  $\underline \bR^{>n}$ iteratively from the previous steps, which is the most efficient way \cite{schollwock11-DMRG}.
Furthermore, we can exploit left- and right-orthogonalization of the cores in order to simplify the tensor
contraction process \cite{dolgovEIG2013,KressnerEIG2014,NLee-Cich14}.

As in any  iterative optimization  based on gradient descent the cost function can only decrease, however we have no guarantee
that a global minimum is achieved. Moreover, for some initial conditions the iteration process can be slow.
To alleviate these problems we can exploit double site DMRG,
 in which we join two neighboring factors (cores), optimize the resulting ``supernode''  called also
 ``super-core'' or ``super-block'', and
split again the result into separated factors by  low-rank matrix factorizations
\cite{dolgovEIG2013,KressnerEIG2014,schollwock11-DMRG,Kressner-Uschmajew2014}.\\

\subsubsection{\bf Tensor Network for Computing Several Extreme Eigenvalues and
Corresponding Eigenvectors for Symmetric Eigenvalue Problem}

In a more general case, in order  to compute a  few, say $K$ eigenvectors corresponding to $K$ algebraically smallest eigenvalues for
a symmetric  matrix $\bA \in \Real^{I \times I}$, we
can employ the following trace minimization problem with orthogonality constraints
\be
\min_{\bX} \tr(\bX^T \bA \bX),\qquad \mbox{s.t.} \;\; \bX^T \bX =\bI_K,
\label{Trace1}
\ee
where $\bX=[\bx_1,\bx_2,\ldots,\bx_K] \in \Real^{J \times K}$, which is equivalent to the following unconstrained problem
\be
\min_{\bX} \{ \tr(\bX^T \bA \bX) +  \alpha ||\bX^T \bX -\bI_K||^2_F \},
\label{Trace1u}
\ee
where the penalty parameter $\alpha >0$ takes suitable finite value \cite{Wen2013trace}.

\begin{figure}[t]
(a) Block tensor train with left- and right-orthogonal cores  \\
\includegraphics[width=8.59cm]{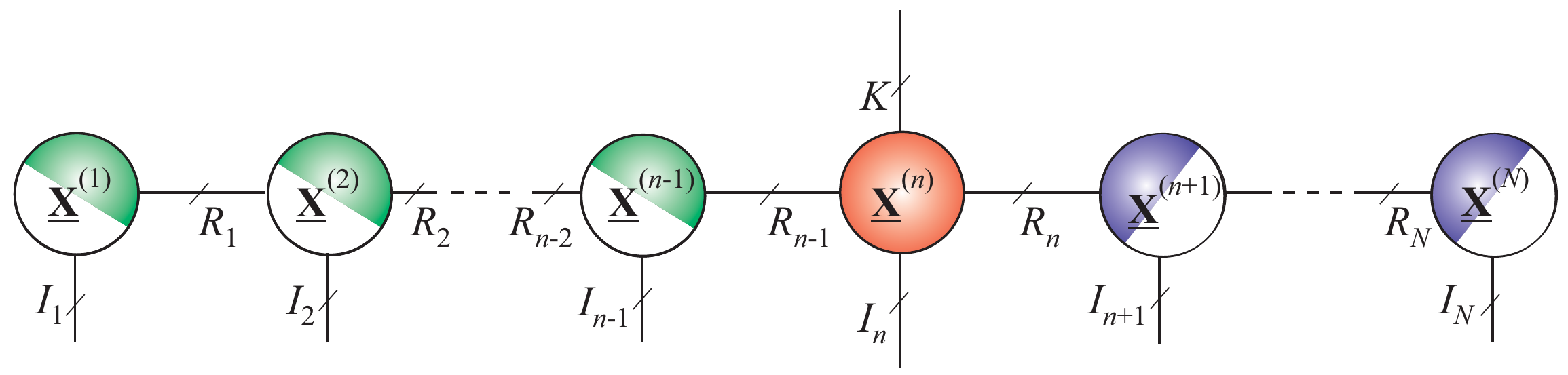}\\
(b)  Tensor network corresponding to the optimization problem (\ref{Trace1}) \\
\includegraphics[width=8.99cm,height=10.01cm]{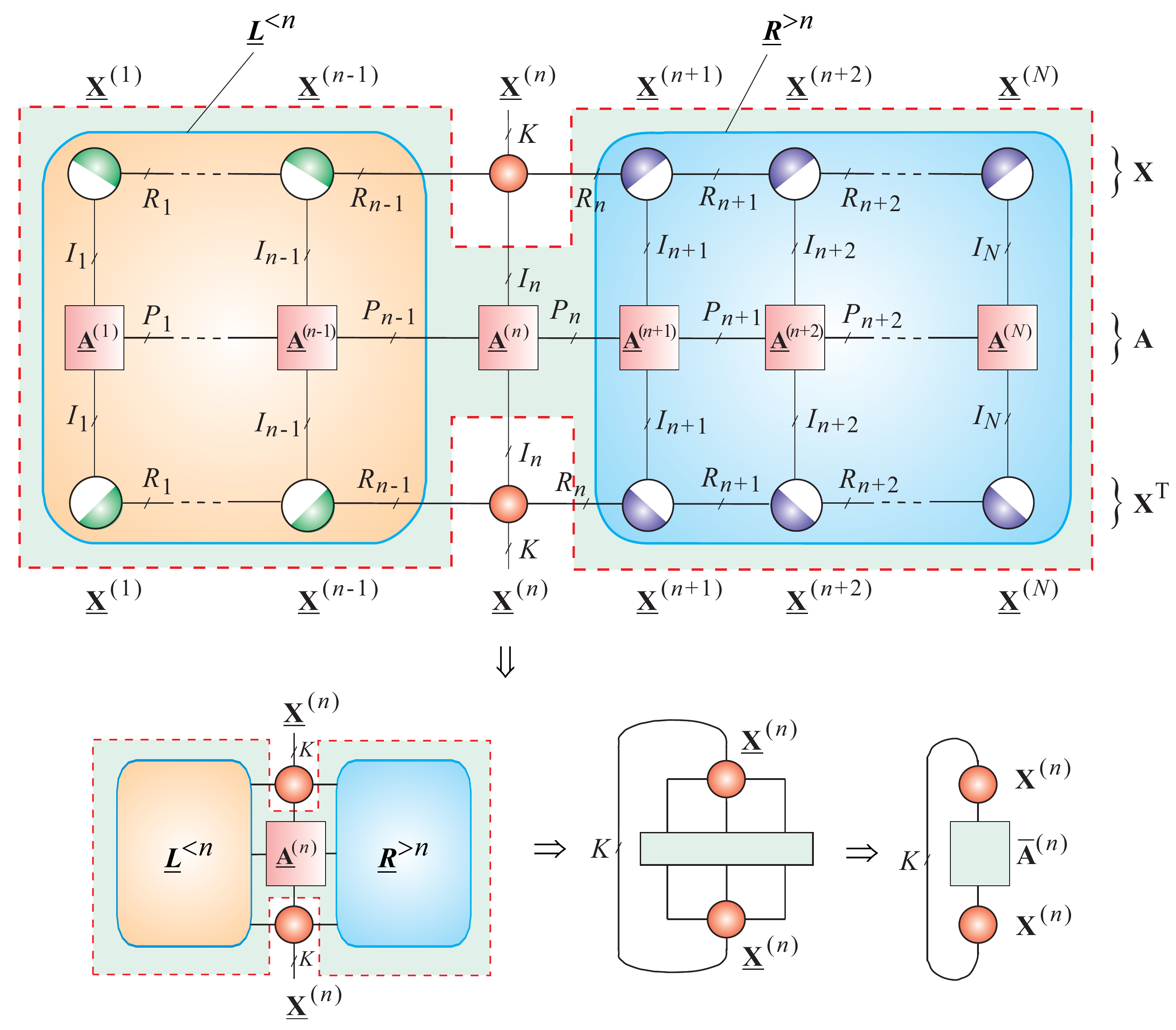}
%
\caption{Computation of $K$  eigenvectors corresponding to $K$ extreme eigenvalues in TT format for the symmetric matrix $\bA \in \Real^{I \times I}$ and the orthogonal matrix $\bX \in \Real^{I\times K}$ given in the distributed block tensor train formats. The extreme eigenvalues are computed as $\mbi \Lambda = \bX^{(n) \;T} \bar \bA^{(n)} \bX^{(n)}$.}
\label{Fig:EIG1K}
\end{figure}

 When computing $K>1$  eigenvectors, we need to work with $K$ vectors $\bx_k$ in parallel. Instead of representing each vector individually in the TT format, we can represent them jointly
 in a block TT format{\footnote{Instead of the block TT format for distributed matrix representations, we can use alternative models, see  Fig. \ref{Fig:TPO-traces}.}} introduced by   Dolgov et al. \cite{dolgovEIG2013} (see also Pi\ifmmode \check{z}\else \v{z}\fi{}orn, I. and Verstraete \cite{NRG-DMRG12} and  Kressner et al. \cite{KressnerEIG2014}). 
In the block TT all cores are 3rd-order tensors,  except one which is 4th-order tensor,
  where additional physical index $K$ represents the number of vectors as shown
  in Fig. \ref{Fig:EIG1K} (a).
It should be noted that the  position of such 4th-order core  $\underline \bG^{(n)}$,
which carries the index $K$  is not fixed;  we will move it back and forth
 from position 1 to $N$   during the sequential optimization \cite{dolgovEIG2013,Holtz-TT-2012}.

 If the block TT model is used to represent  $K$ orthogonal vectors, then
 the matrix frame equation takes the slightly modified form:
\be
\bX= \bX_{\neq n} \bX^{(n)} \in \Real^{I \times K}, \quad n=1,2,\ldots,N,
\ee
where  $\bX^{(n)} \in \Real^{R_{n-1} I_n R_n \times K}$.

Hence, we can express the trace in (\ref{Trace1}) as follows:
\be
\tr(\bX^T \bA \bX) &=& \tr((\bX_{\neq n} \bX^{(n)})^T \bA \bX_{\neq n} \bX^{(n)}) \notag \\
&=&\tr(( \bX^{(n)})^T [\bX_{\neq n}^T \bA \bX_{\neq n}] \bX^{(n)})\notag \\
&=& \tr(\bX^{(n)\,T}  \bar \bA^{(n)} \bX^{(n)}),
\ee
where $\bar \bA^{(n)} = \bX_{\neq n}^T \bA \bX_{\neq n} $.

Assuming that frame matrices have orthogonal columns,
  we can convert the optimization  problem (\ref{Trace1}) into a set of linked
   optimization problems:
\be
\min_{\bX^{(n)}} \tr( \bX^{(n) \,T} \bar \bA^{(n)} \bX^{(n)}),\quad \mbox{s.t.} \;\; \bX^{(n) \;T} \bX^{(n)}=\bI_K
\ee
for $n=1,2,\ldots,N$, where $\bar \bA^{(n)}$ is computed
iteratively by  tensor network contraction shown in Fig. \ref{Fig:EIG1K} (b).
In other words, the above problem is solved iteratively via optimized iterative contraction of the tensor network. This means that an active block (core) is sequentially selected in an iterative  manner
 for $n=1,2,\ldots,N$ by sweeping from left to right and back from right to left and so on until convergence \cite{dolgovEIG2013,KressnerEIG2014}.

It should be noted that the global orthogonality constraint $\bX^T \bX =\bI_K$ is equivalent to the set of local orthogonality constraints $(\bX^{(n)})^T \bX^{(n)}=\bI_K,\; \forall n$, since due to left and right orthogonality of the cores, we can write:
\be
\bX^T \bX &=& \bX^{(n)\,T } \bX^T_{\neq n}  \bX_{\neq n} \bX^{(n)} \\
&=& \bX^{(n)\,T} \bX^{(n)}, \;\; \forall n. \nonumber
\ee

\subsection{\bf Tensor Networks for Tracking a Few Extreme Singular Values and Vectors for SVD and Sparse PCA}

\begin{figure*}[p]
(a)
\vspace{-0.6cm}
\begin{center}
\includegraphics[width=12.29cm]{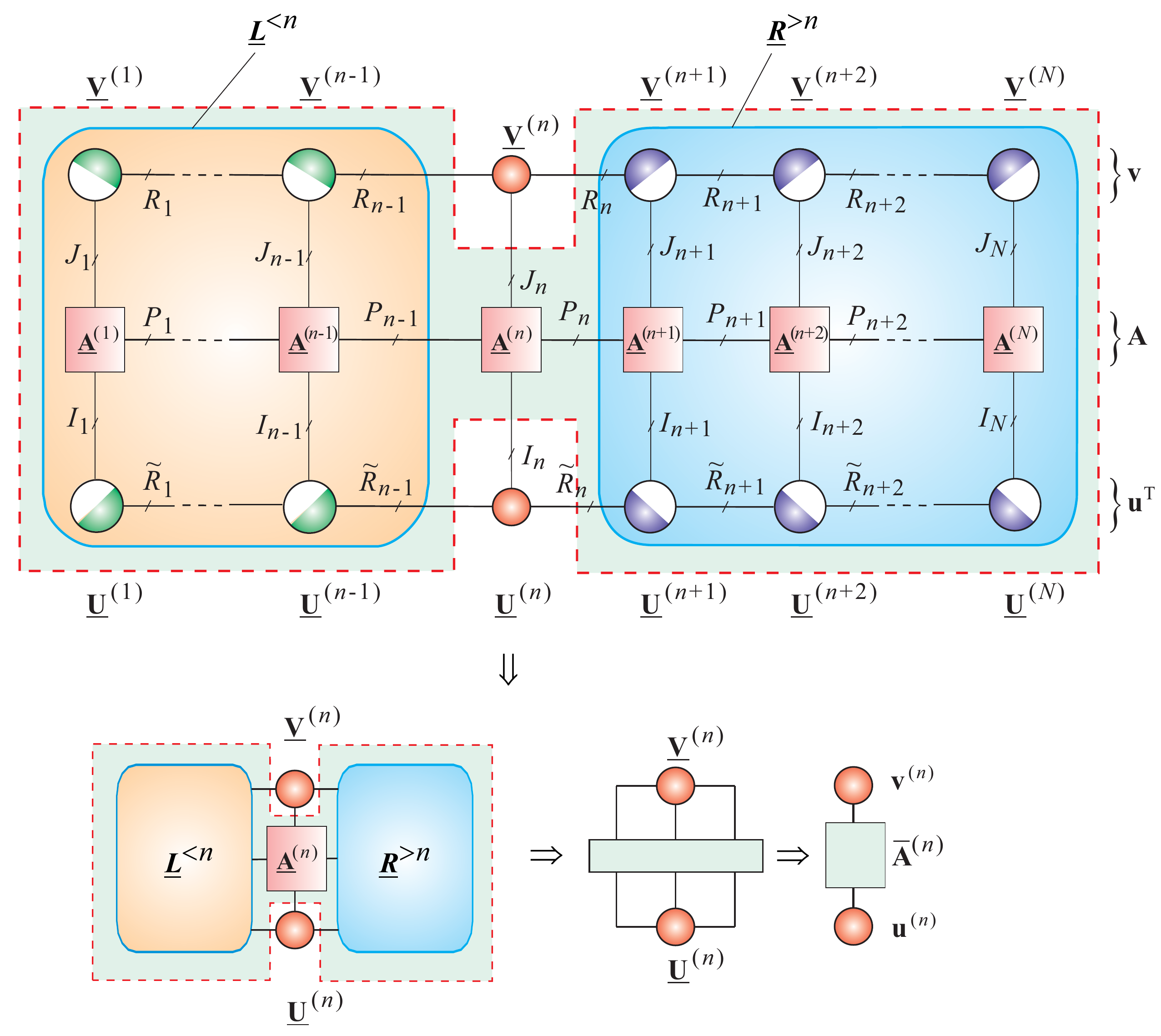}\\
\end{center}
(b)
\vspace{-0.3cm}
\begin{center}
\includegraphics[width=12.29cm]{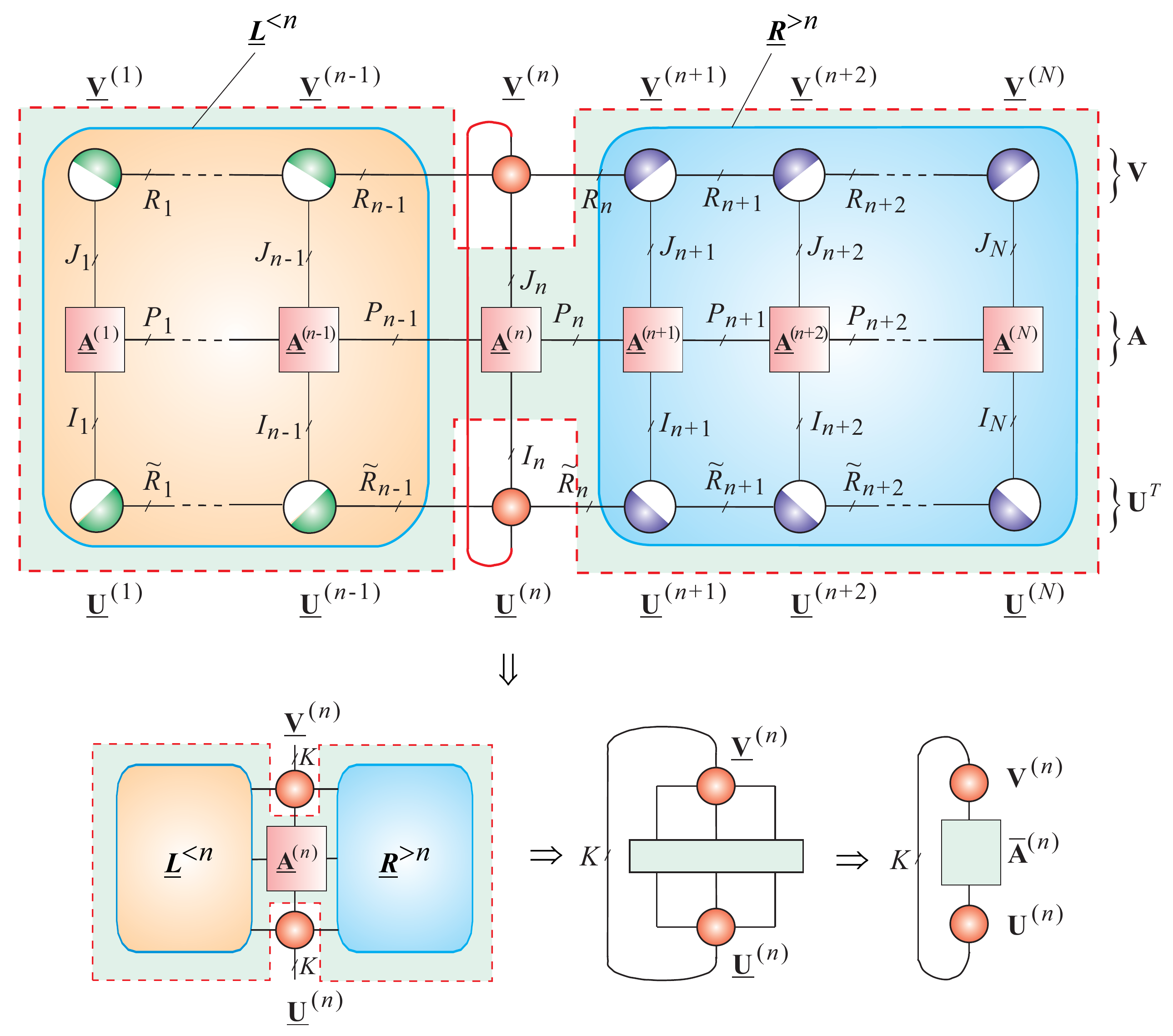}
\end{center}
\caption{(a) Tensor network for computing the SVD singular eigenvectors corresponding to a largest singular value. (b) Tensor network for computing $K$  left- and right-eignevectors corresponding to the $K$ largest singular values via maximization of the trace $\tr(\bU^T \bA \bV)$, subject to orthogonality constraints $\bU^T \bU =\bI_k$ and $\bV^T\bV=\bI_K$ \cite{Lee-Cich-TTSVD}. The singular values are computed as $\mbi \Sigma =\bU^{(n) \;T} \bar \bA^{(n)} \bV^{(n)}$.}
\label{Fig:SVD1}
\end{figure*}

The computation of the largest singular value and the corresponding left- and right eigenvector can be performed via the following optimization problem
\be
\max_{\bu,\bv} \{ \bu^T \bA  \bv \}, \;\; \mbox{s.t.} \;\; ||\bu||_2^2 = 1,\;  ||\bv||_2^2 = 1,
\label{SVD}
\ee
where   $\bA \in \Real^{I \times J}$ is arbitrary data matrix that admits low-rank TT decomposition.
Using TT decomposition of vectors $\bu \in \Real^I$, $\bv\in \Real^J$ and  the data matrix $\bA$ and assuming that cores $\underline \bU^{(n)}$ and $\underline \bV^{(n)}$ are  kept left- and right- orthogonal, the optimization problem (\ref{SVD}) can be converted into a set of usually much smaller scale optimization problems as follows (see Fig. \ref{Fig:SVD1} (a)):
\be
&& \max_{\bu^{(n)}, \bv^{(n)}} \{ (\bu^{(n)})^T \bar \bA^{(n)} \bv^{(n)} \}, \\
&& \mbox{s.t.} \quad ||\bu^{(n)}||^2_2=1,\;\; ||\bv^{(n)}||^2_2=1, \; \forall n, \nonumber
\label{SVDframe}
\ee
where $\bu^{(n)} =vec(\underline\bU^{(n)}) \in \Real^{\widetilde R_{n-1} I_n \widetilde R_n}$ and $\bv^{(n)} =vec(\underline \bV^{(n)}) \in \Real^{R_{n-1} J_n R_n}$ and
\be
\bar \bA^{(n)} = (\bU_{\neq n})^T \bA \bV_{\neq n} \in \Real^{\widetilde R_{n-1} I_n \widetilde R_n \times R_{n-1} J_n R_n}
\ee
for $n=1,2,\ldots,N$.

Note that taking into account that the frame matrices $\bU_{\neq n} \in \Real^{I_1 I_2 \cdots I_N \times \widetilde R_{n-1} I_n \widetilde R_n}$ and
$\bV_{\neq n} \in \Real^{J_1 J_2 \cdots J_N \times  R_{n-1} I_n  R_n}$ are orthogonal and
\be
\bu = \bU_{\neq n} \bu^{(n)}, \;\; \bv = \bV_{\neq n} \bv^{(n)},\;\; \;\; \forall n,
\ee
 we can easily check that
$||\bu||_2 =||\bu^{(n)}||_2$ and
$||\bv||_2 =||\bv^{(n)}||_2, \;\; \forall n$.

 An alternative approach to compute SVD for several maximal singular values and the corresponding left- and right- orthogonal eigenvectors, is to convert the SVD to the problem of  symmetric EVD by applying the following basic relationships.
It is
evident, that from the SVD of the matrix $\bA = \bU \, \mbi {\Sigma} \,
\bV^T \in \Real^{I \times J}$, where $\mbi {\Sigma}_1 = \diag\{ \sigma_1, \ldots, \sigma_R \}$,
we have
\be
\bA \bA^T &=& \bU \mbi {\Sigma}^2_{1} \bU^T, \\
\bA^T \bA &=& \bV \, \mbi {\Sigma}^2_{2} \, \bV^T,
\ee
where $\mbi {\Sigma}_1 = \diag\{ \sigma_1, \ldots, \sigma_I \}$
and $\mbi {\Sigma}_2 = \diag\{ \sigma_1, \ldots, \sigma_I \}$.
This means that the singular values of $\bA \in \Real^{I \times J}$ are the positive
square roots of the eigenvalues of $\bA^T \bA$ and the  eigenvectors $\bU$ of
$\bA \bA^T$ are the left singular vectors of $\bA$.
Note that if $R < I$, the matrix $\bA \bA^T$ will contain at least $I-R$
additional eigenvalues that are not included as singular values of $\bA$.

\begin{figure}[t]
(a)\\
\includegraphics[width=8.99cm]{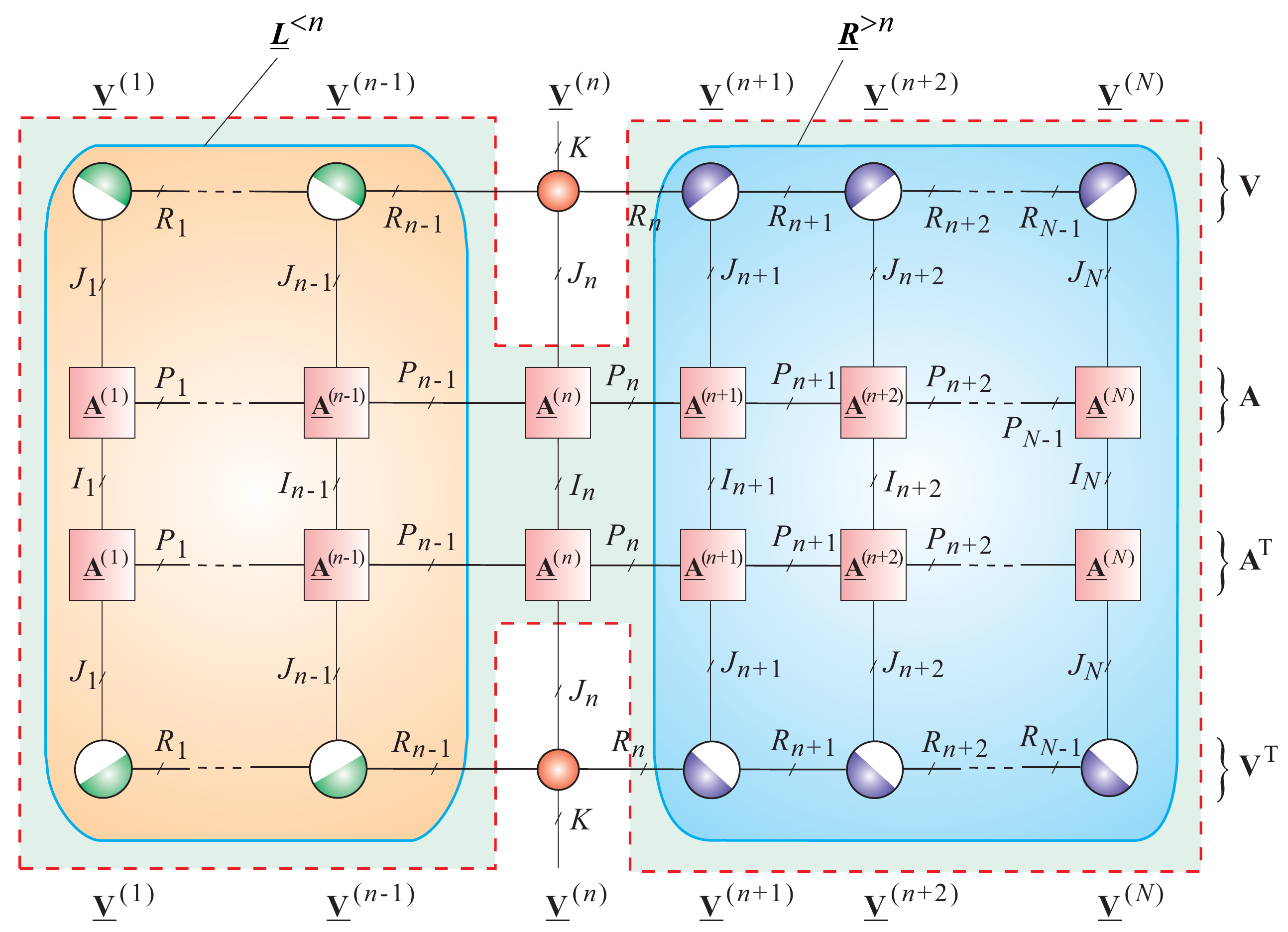}\\
\vspace{0.1cm}
(b)
\begin{center}
\includegraphics[width=8.69cm]{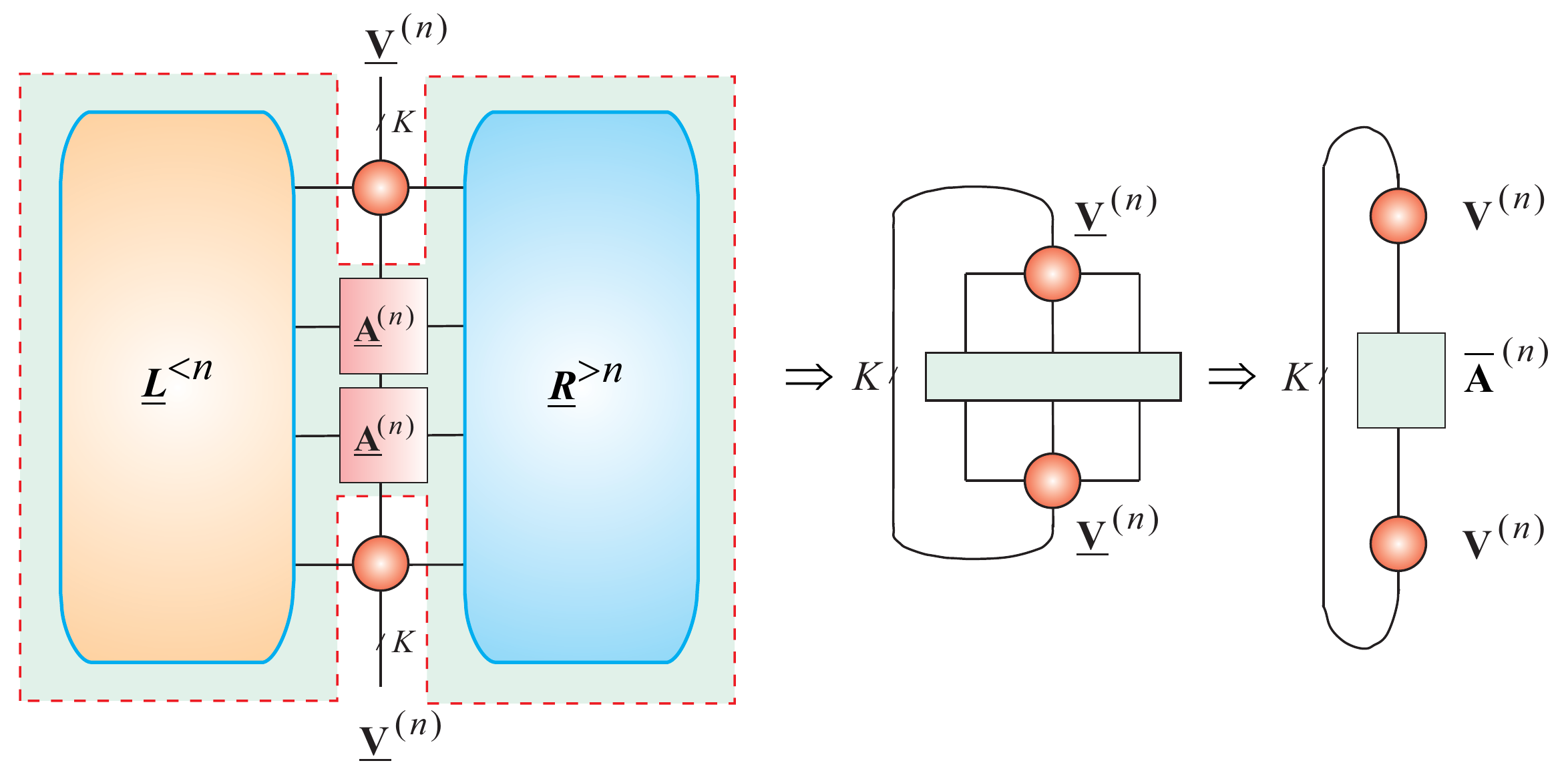}
\end{center}
\caption{Computation  of $K$  right eigenvectors corresponding to $K$ smallest singular values of the SVD in TT formats.}
\label{Fig:SVDK}
\end{figure}

Hence, in order to compute approximately $K$ smallest singular values and the corresponding right-eigenvectors,
we can employ formally the following optimization problem:
\be
&&\min_{\bV \in \Real^{I \times K}} \tr(\bV^T \bA^T \bA  \bV), \\
&& \mbox{s.t.}  \quad \bV^T \ \bV = \bI_K. \notag
 \ee
 The SVD problem for large--scale structured matrices that admit low-rank TT approximations
  can be solved iteratively in TT formats by the following set of smaller
   optimization (symmetric EVD) problems:
 \be
&&\max_{\bV^{(n)}} \tr((\bV^{(n)})^T [\bV_{\neq n}^T \bA^T \bA \bV_{\neq n}] \bV^{(n)}), \\
&& \mbox{s.t.}  \quad (\bV^{(n)})^T  \bV^{(n)}= \bI_K, \;\; n=1,2,\ldots,N, \notag
 \ee
where $\bV^{(n)} \in \Real^{R_{n-1} J_n R_n \times K}$ and
 \be
\bar \bA^{(n)}=\bV_{\neq n}^T \bA^T \bA  \bV_{\neq n} \in \Real^{R_{n-1} I_n R_n \times R_{n-1} I_n R_n}
  \ee
 for $n=1,2,\ldots,N$  are computed sequentially via tensor network contractions as illustrated
 in Fig. \ref{Fig:SVDK}.

The challenge is how to  extend and/or modify the above described approaches to the following large-scale  optimization problems for  structured matrices, if we need impose additional constraints such as sparsity, nonnegativity, orthogonality  or local smoothness:

\begin{itemize}

\item Sparse Principal Component Analysis (SPCA) using the Penalized Matrix Decomposition (PMD)   \cite{SPCA-Witten,Witten-PHD} 
\be
&&\max_{\bu,\bv} \{ \bu^T \bA  \bv \}, \quad \mbox{s.t.} \quad ||\bu||_2^2 \leq 1,\;\;  ||\bv||_2^2 \leq 1, \notag\\
&& P(\bv) \leq c_1,
\label{PMD-PCA}
\ee
where the positive parameter $c_1$ controls sparsity level and the convex penalty function $P(\bv)$ can take a variety of forms. Useful examples are \cite{Witten-PHD}:
\be
P(\bv)&=& ||\bv||_1 =\sum_{i=1}^I |v_i| \quad \mbox{(Lasso)}, \nonumber\\
P(\bv)&=& ||\bv||_0 = \sum_{i=1}^I |\sign (v_i)|,  \\
P(\bv) &=& \sum_{i=1}^I |v_i| + \lambda \sum_{i=2}^I |v_i-v_{i-1}|.  \nonumber
\label{PvPMD}
\ee

\item SPCA via regularized SVD (sPCA-rSVD)  \cite{Shen-Huang-SPCA-08},
\cite{SPCA-GPM}
\be
&& \max_{\bu,\bv} \; \{ \bu^T \bA  \bv - \alpha P(\bv)\} \notag \\
 &&\mbox{s.t.} \quad ||\bu||_2^2 \leq 1,\;\;  ||\bv||_2^2 \leq 1,
\label{rSVD}
\ee

\item Two-way functional PCA/SVD   \cite{Huang-2way} 
\be
&& \max_{\bu,\bv} \{ \bu^T \bA  \bv- \frac{\alpha}{2} P_1(\bu) P_2(\bv) \}, \notag \\
 &&\mbox{s.t.} \quad ||\bu||_2^2 \leq 1,\;\;  ||\bv||_2^2 \leq 1.
\label{FPCA}
\ee

\item Sparse SVD  \cite{Lee-2010-bicluster} 
\begin{equation}
\max_{\bu, \bv} \{ \bu^T \bA  \bv- \frac{1}{2} \bu^T \bu \bv^T \bv - \frac{\alpha_1}{2} P_1(\bu) - \frac{\alpha_2}{2} P_2(\bv) \}.\\
\label{SSVD}
\end{equation}

\item Generalized SPCA   \cite{Allen-HOPCA} 
\be
&&\max_{\bu, \bv} \{ \bu^T \bQ \bA \bR \bv - \frac{\alpha_1}{2} P_1(\bu) - \frac{\alpha_2}{2} P_2(\bv) \}, \notag \\
&& \mbox{s.t.} \quad \bu^T \bQ \bu  \leq 1,\;\;  \bv^T \bR \bv  \leq 1,
\label{GPCA}
\ee
where  $\bQ \in \Real^{T \times T}$ and $\bR \in \Real^{I \times I}$ are symmetric
positive-definite matrices.

\item Generalized nonnegative SPCA   \cite{Allen-NGPCA}
\be
&& \max_{\bu,\bv} \{ \bu^T \bA \bR \bv - \alpha ||\bv||_1 \}, \notag\\
&& \mbox{s.t.} \quad \bu^T \bu  \leq 1,\;\;  \bv^T \bR \bv  \leq 1,\;\; \bv \geq 0.
\label{NPCA}
\ee
%
\end{itemize}

\begin{table*}[ht]
\caption{Cost functions and constraints used in classical feature extraction
(dimension reduction) methods that
can be formulated  as generalized eigenvalue problem (\ref{GEVDK11}).
The objective is to find  an (orthogonal) matrix $\bV$, assuming that data matrices $\bX$, $\bW, \bD, \bH$ are known. The symmetric matrix $\bA$ can take different forms: $\bA= {\bf I}-\frac{1}{N} {\bf 11}^{T}, \; \bA=\bD-\bW,   \; \bA=(\bI - \bW^T)(\bI -\bW), \; \bA= \bI-\bH,$ depending on  method (for more detail see \cite{Kokiopoulou11}).}
\centering
\renewcommand{\arraystretch}{1.4}
{\shadingbox{
\begin{tabular}{l|c|c}
\hline
Method & Cost Function (min) & Constraints \\
\hline
Principal Component Analysis/\\/Multi-Dimensional Scaling (PCA/MDS) & {$\tr[-{\bf V}^{T}{\bf X}({\bf I}-\frac{1}{N} {\bf 11}^{T}){\bf X}^{T}{\bf V}]$} & ${\bf V}^{T}{\bf V}={\bf I}$ \\
 & & \\
Locally Preserving Projection (LPP) & {$\tr[{\bf V}^{T}{\bf X}({\bf D}-{\bf W}){\bf X}^{T}{\bf V}]$} & ${\bf V}^{T}{\bf XDX}^{T}{\bf V}={\bf I}$ \\
 & & \\
Orthogonal LPP (OLPP) & {$\tr[{\bf V}^{T}{\bf X}({\bf {D-W}}){\bf X}^{T}{\bf V}]$} & ${\bf V}^{T}{\bf V}={\bf I}$ \\
 & & \\
Neighborhood Preserving Projection (NPP) & {$\tr[{\bf V}^{T}{\bf X}({\bf I}-{\bf W}^{T})({\bf I}-{\bf W}){\bf X}^{T}{\bf V}]$} & ${\bf V}^{T}{\bf XX}^{T}{\bf V}={\bf I}$ \\
 & & \\
Orthogonal NPP (ONPP) & {$\tr[{\bf V}^{T}{\bf X}({\bf I}-{\bf W}^{T})({\bf I}-{\bf W}){\bf X}^{T}{\bf V}]$} & ${\bf V}^{T}{\bf V}={\bf I}$ \\
 & & \\
Linear Discriminant Analysis (LDA) & {$\tr[{\bf V}^{T}{\bf X}({\bf I}-{\bf H}){\bf X}^{T}{\bf V}]$} & ${\bf V}^{T}{\bf XX}^{T}{\bf V}={\bf I}$ \\
 & & \\
Spectral Clustering (Ratio Cut) & {$\tr[{\bf V}^{T}({\bf D}-{\bf W}){\bf V}]$}
& ${\bf V}^{T}{\bf V}={\bf I}$ \\
Spectral Clustering (Normalized Cut) & {$\tr[{\bf V}^{T}({\bf D}-{\bf W})
{\bf V}]$} & ${\bf V}^{T}{\bf D V}={\bf I}$ \\
\hline
\end{tabular}
}}
\label{Table:LPP}
\end{table*}

\begin{figure*}[ht!]
\centering
\includegraphics[width=16.99cm]{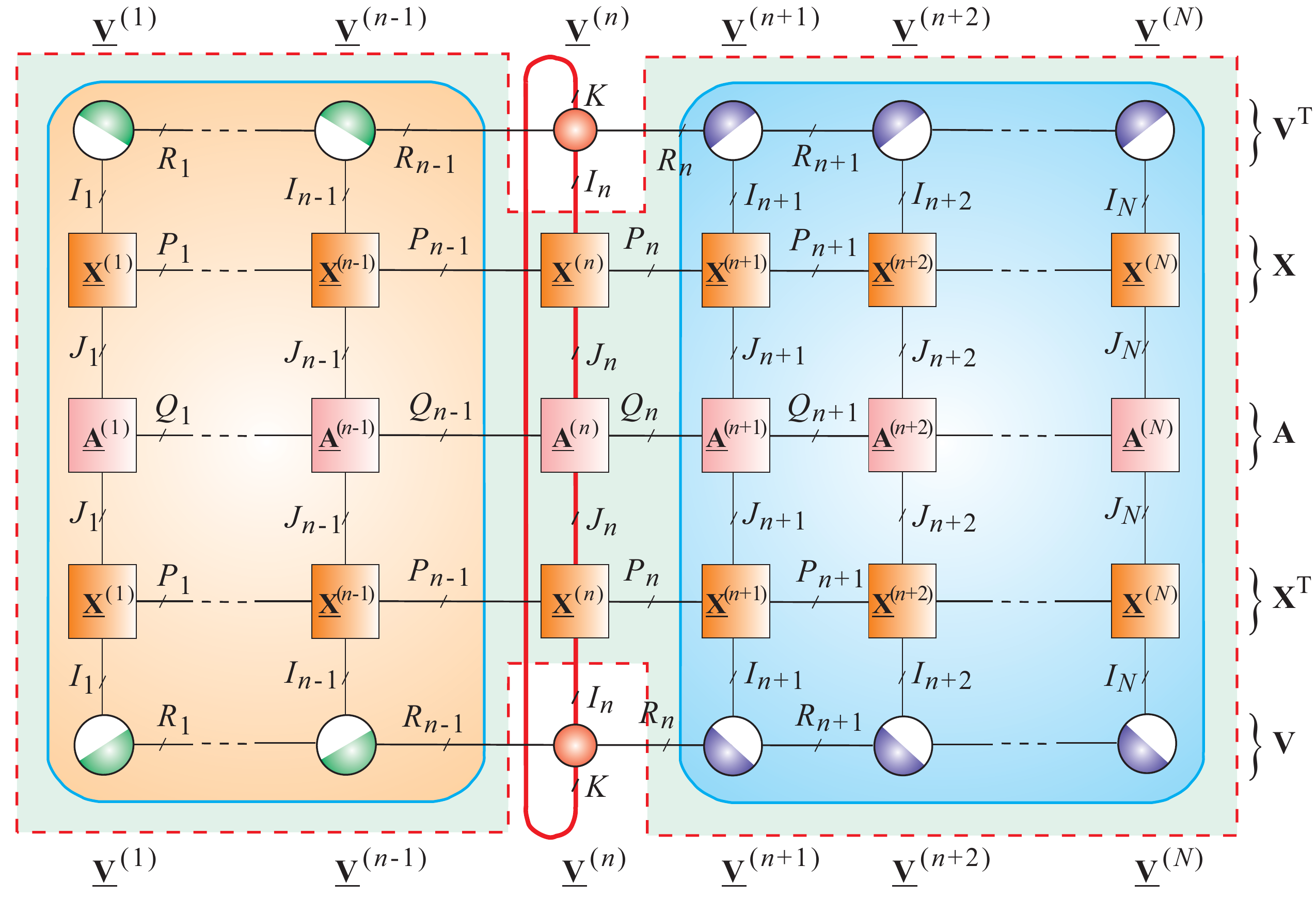}
\caption{Tensor network for computation of $K$  eigenvectors corresponding to the $K$ extreme eigenvalues in TT formats for the generalized eigenvalue problem (\ref{GEVDK11}).}
\label{Fig:GEVDK}
\end{figure*}


\subsection{\bf Generalized Eigenvalue Problems in TT formats}

In many practical applications, especially in dimension reduction and classification problems
 (e.g.,  in PCA/MDS, LPP, ONPP, LDA -- see Table \ref{Table:LPP} for more detail), we need to minimize the following trace  optimization problem formulated as a generalized eigenvalue problem (GEVD) \cite{Kokiopoulou11}:
 \be
\min_{\bV \in \Real^{I \times K}} \tr(\bV^T \bX \bA \bX^T \bV), \quad \mbox{s.t.}  \quad \bV^T \bB \bV = \bI_K,
\label{GEVDK11}
 \ee
%
 where it is assumed that the structured data matrices: $\bX \in \Real^{I \times J}$, symmetric matrix  $\bA \in \Real^{J \times J}$, and symmetric positive-definite matrix  $\bB \in \Real^{I \times I}$ are known.

 The problem is equivalent to the unconstrained optimization problem
  \be
\min_{\bV \in \Real^{I \times K}} \{ \tr(\bV^T \bX \bA \bX^T \bV)+ \alpha \; || \bV^T \bB \bV - \bI_K||^2_F.
\label{GEVDK11u}
 \ee
 Note that by changing of the variable $\bW=\bB^{1/2} \bV$ the GEVD can be converted to the standard symmetric
 EVD problem
  \be
\min_{\bW \in \Real^{I \times K}} \tr(\bW^T \bB^{-1/2} \bX \bA \bX^T \bB^{-1/2} \bW), \;\; \mbox{s.t.}  \;\; \bW^T  \bW = \bI_K. \nonumber
\label{EVDK12}
 \ee

 The objective is to estimate the  matrix $\bV \in \Real^{I \times K}$ in a TT format, assuming that large-scale matrices $\bX$ and $\bA$ ($\bW, \bD, \bH)$ are known and admit low-rank TT approximations.
  The problem for  structured matrices that admit low-rank TT approximations
  can be solved iteratively:
 \be
&&\min_{\bV^{(n)}} \tr((\bV^{(n)})^T \;[\bV_{\neq n}^T \bX \bA \bX^T \bV_{\neq n}]\; \bV^{(n)}), \\
&& \mbox{s.t.}  \quad (\bV^{(n)})^T \;[\bV_{\neq n}^T \bB \bV_{\neq n}]\; \bV^{(n)}= \bI_K, \notag
\label{TTGEVDK}
 \ee
where the relatively low-dimension matrices:
 \be
  \bar \bA^{(n)}= [\bV_{\neq n}^T \bX \bA \bX^T \bV_{\neq n}]\in \Real^{R_{n-1} I_n R_n \times R_{n-1} I_n R_n}
  \ee
  and
  \be
  \bar \bB^{(n)}= [\bV_{\neq n}^T \bB \bV_{\neq n}] \in \Real^{R_{n-1} I_n R_n \times R_{n-1} I_n R_n}
   \ee
can be  computed sequentially for $n=1,2,\ldots,N$ via tensor network contractions shown in Fig. \ref{Fig:GEVDK}.



\subsection{\bf Canonical Correlation Analysis in TT Format}

The Canonical Correlation Analysis (CCA), introduced by Hotelling, 
 can be considered as a generalization  of PCA and it is
 a classical method for
determining the relationship between two sets of variables.
Given two  zero-mean (i.e., centered) data sets $\bX \in \Real^{I \times J}$ and
$\bY \in \Real^{L \times J}$   on the same set of $J$ observations, CCA seeks linear
combinations of the variables in $\bX$ and the variables in $\bY$
that are maximally mutually correlated with each other.
Formally, the classical CCA computes two projection vectors $\bw_x =\bw^{(1)}_x  \in \Real^I$ and
$\bw_y =\bw^{(1)}_y \in \Real^L$ such that the correlation coefficient
\be
\rho =\frac{\bw^T_x \bX \bY^T  \bw_y}{\sqrt{(\bw^T_x \bX \bX^T \bw_x)(\bw^T_y \bY \bY^T \bw_y)}}
\ee
is maximized.

In a similar way, we  can formulate  kernel CCA by replacing inner product matrices by kernel matrices:
\be
\rho = \max_{\mbi \alpha_x, \mbi\alpha_y} \frac{\mbi \alpha^T_x \bK_x \bK_y \mbi \alpha_y}{\sqrt{(\mbi \alpha^T_x \bK_x \bK_x \mbi \alpha_x)(\mbi \alpha^T_y \bK_y \bK_y \mbi \alpha_y)}},
\ee
where $\bK_x \in \Real^{J \times J}$ and $\bK_y \in \Real^{J \times J} $ are suitably designed kernel matrices. 
The above optimization problem can be reformulated as a generalized eigenvalue decomposition (GEVD).

Since $\rho$ is invariant to the scaling of the vectors $\bw_x$ and $\bw_y$, the standard CCA can be equivalently formulated as the following constrained optimization problem:
\be
&&\max_{\bw_x,\bw_y} \{ \bw_x^T \bX  \bY^T \bw_y  \} \\
 && \mbox{s.t.} \quad \bw_x^T \bX  \bX^T \bw_x=\bw_y^T \bY^T  \bY \bw_y=1.
\label{CCA}
\ee
We will refer to $\bt_1=\bX^T \bw_x$ and $\bu_1=\bY^T \bw_y$ as the canonical variables.

\begin{figure}[t]
\centering
\includegraphics[width=8.99cm,height=7.1cm]{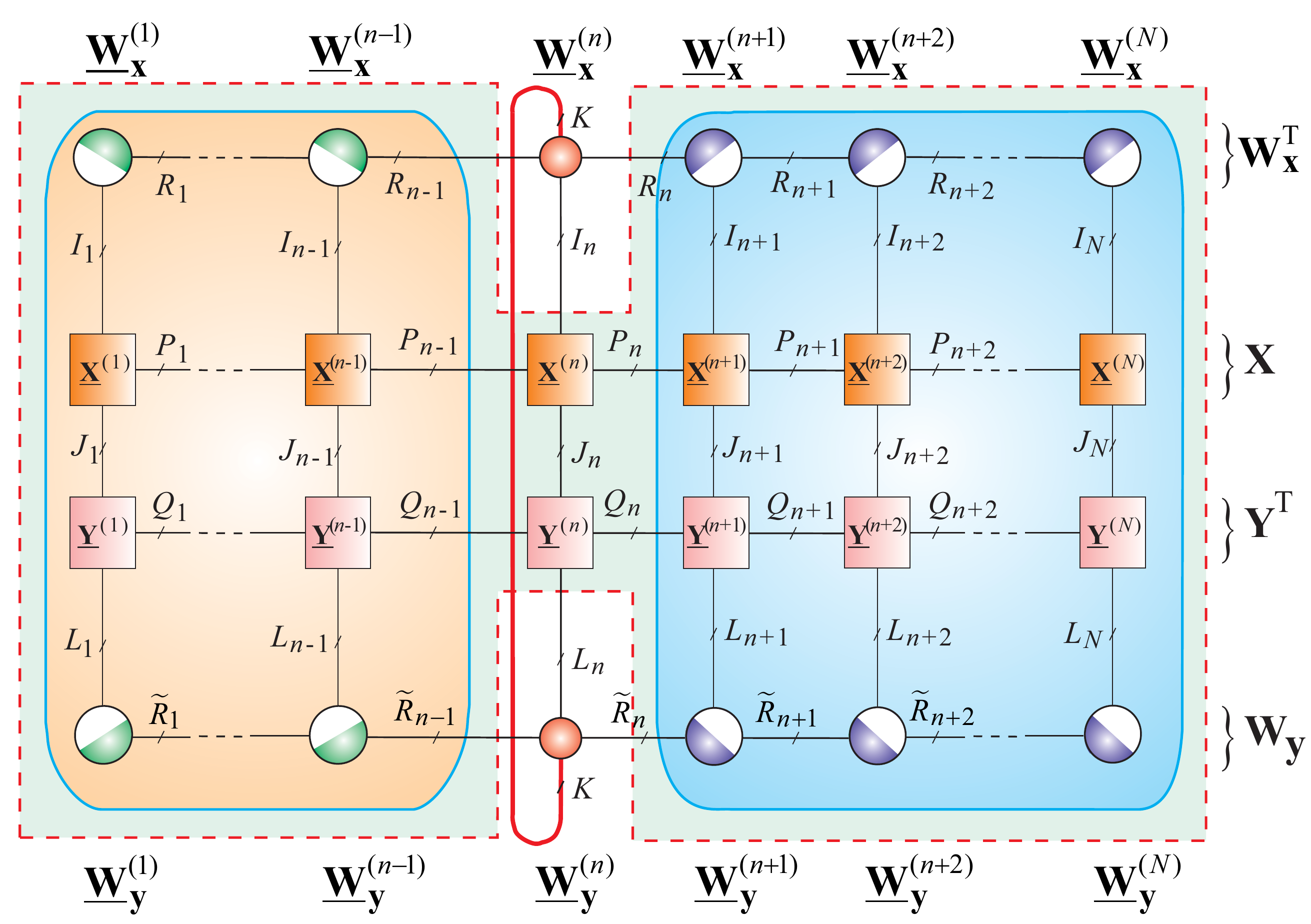}
\caption{Tensor network for computation of multiple (sparse) CCA.}
\label{Fig:CCAL}
\end{figure}

For  sparse CCA, we usually assume that the columns of $\bX$ and $\bY$ have been standardized to have zero
mean and standard deviation one.
 The cross product matrices $ \bX \bX^T$ and  $ \bY \bY^T$ are often approximated by identity matrices, and consequently the constraints $\bw_x^T \bX \bX^T \bw_x \leq 1$ and $\bw_y^T \bY \bY^T \bw_y \leq 1$ can
 be simplified as  $||\bw_x||_2^2 \leq 1$ and $||\bw_y||_2^2 \leq 1$, respectively under some conditions \cite{Witten-PHD}.
Hence, in order to compute sparse CCA we must impose  suitable sparsity constraints on the canonical vectors, for example, by applying the PMD approach \cite{SPCA-Witten,Witten-PHD}:
\be
&&\max_{\bw_x,\bw_y} \{\bw_x^T \bX \bY^T \bw_y \} \\
 &&\mbox{s.t.} \quad ||\bw_x||_2^2 \leq 1, \;\; ||\bw_y||_2^2 \leq 1, \notag \\
&&  P_1(\bw_x)\leq c_1, \;\; P_2(\bw_y)\leq c_2, \notag
\label{PMD-CCA}
\ee
where $P_1$ and $P_2$ are convex penalty functions and positive parameters $c_1,c_2$ control sparsity level. Since $P_1$ and $P_2$ are generally chosen to yield sparse projection vectors  $\bw_x$ and $\bw_y$, we call this criterion the sparse CCA (see Eqs. (\ref{PvPMD})).



In order to compute multiple  canonical vectors for the standard CCA, we  can formulate the following
optimization problem:
\be
\max_{\bW_x,\bW_y}&\{tr( \bW^T_x \bX \bY^T \bW_y)\}, \\
\mbox{s.t.}  &\bW^T_x \bX \bX^T \bW_x =\bI_K \notag \\
& \bW^T_y \bY \bY^T \bW_y =\bI_K, \notag
\ee
where $\bW_x =[\bw_x^{(1)},\bw_x^{(2)},\ldots,\bw_x^{(K)}] \in \Real^{I \times K}$ and
$\bW_y =[\bw_y^{(1)},\bw_y^{(2)},\ldots,\bw_y^{(K)}] \in \Real^{L \times K}$.


This optimization scheme in a TT format is illustrated in Fig. \ref{Fig:CCAL} and performs
iteratively the following set of optimization problems:
\be
&& \max_{\bW^{(n)}_x,\bW^{(n)}_y} \tr( (\bW_x^{(n)})^T  \;[\bW^T_{x,\;\neq n}\bX \bY^T \bW_{y, \;\neq n} ]\; \bW^{(n)}_y), \nonumber \\
 &&\mbox{s.t.} \; \; (\bW_x^{(n)})^T \;[\bW^T_{x,\; \neq n} \bX \bX^T \bW_{x,\;\neq n}] \; \bW^{(n)}_x =\bI_K \nonumber \\
&& \qquad (\bW_y^{(n)})^T  \;[ \bW^T_{y, \; \neq n} \bY \bY^T  \bW_{y,\; \neq n}]\; \bW^{(n)}_y =\bI_K, \nonumber \\
&&\qquad n=1,2,\ldots,N.
\ee

 Note that for large-scale sparse CCA  the cross product matrices $ \bX \bX^T$ and  $ \bY \bY^T$ can be approximated by identity matrices, and consequently the  above constraints can  be simplified \cite{SPCA-Witten}.

\subsection{\bf Solving  Large-Scale Systems of Linear Equations}

\begin{figure}[t]
\centering
\includegraphics[width=8.9cm]{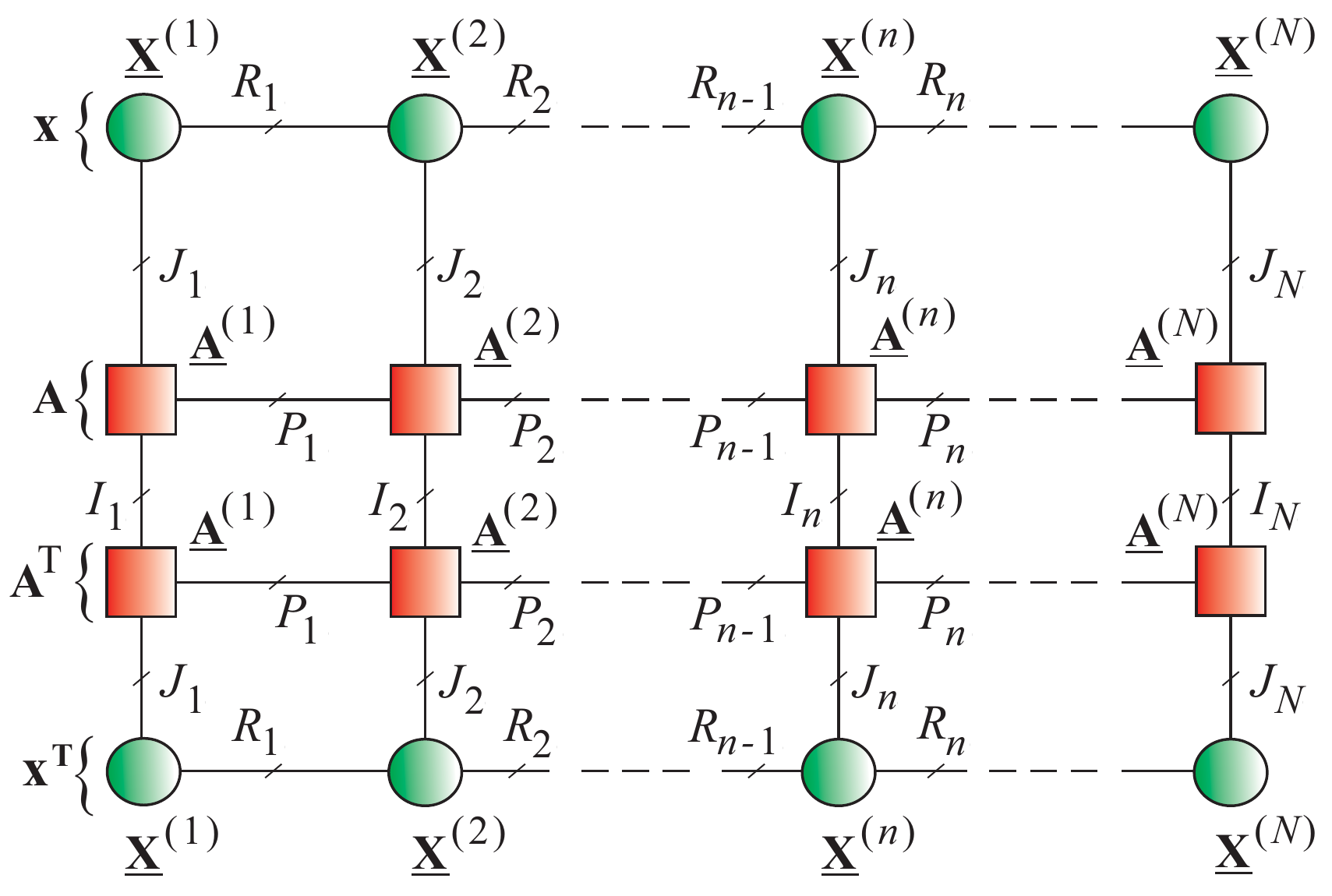}\\
\vspace{0.5cm}
\includegraphics[width=8.9cm]{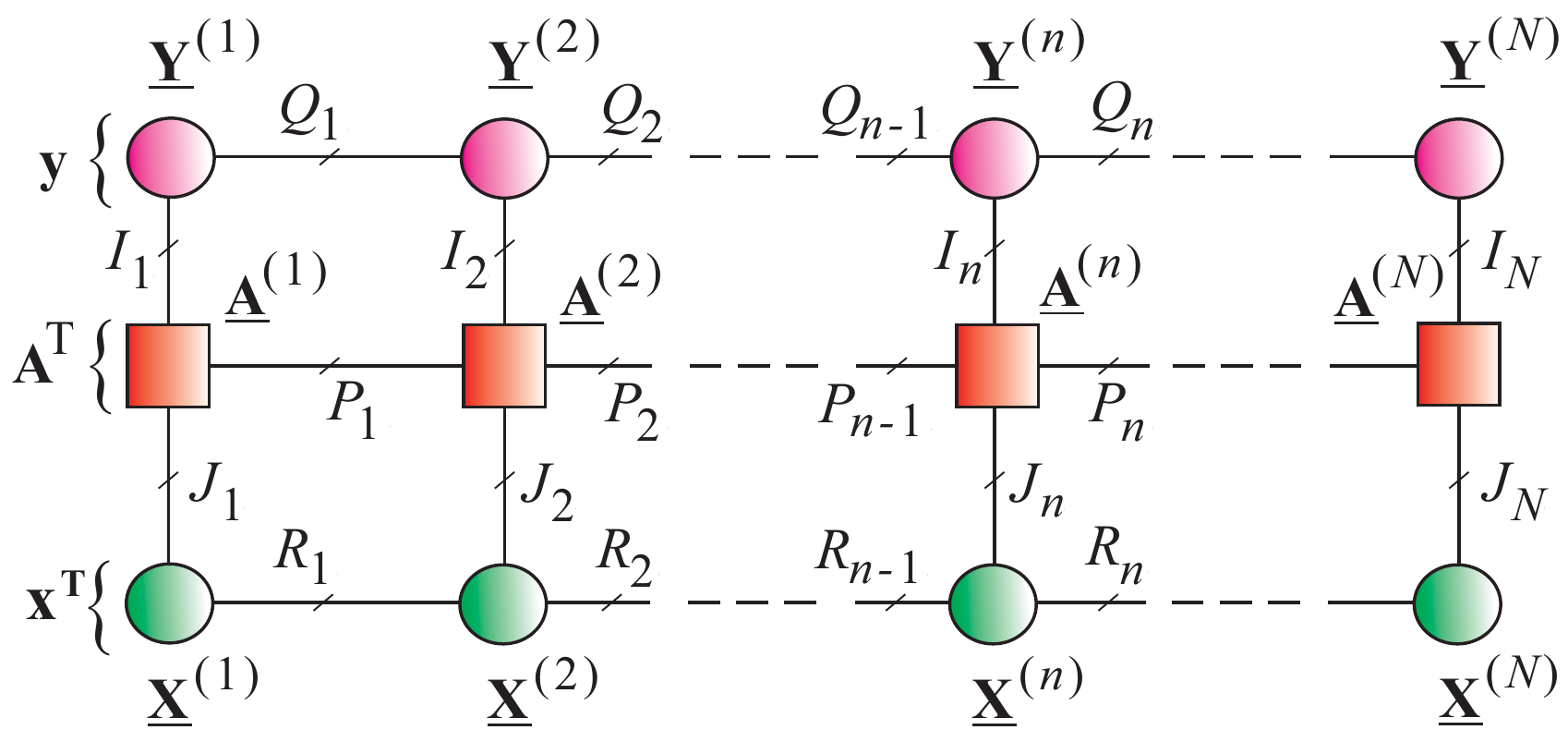}
\caption{Simplified tensor network scheme for solving systems of linear equations with a  huge non-symmetric matrix $\bA$.}
\label{Fig:ATAX=Y}
\end{figure}

Consider  a huge system of linear algebraic equations
in TT formats:
\be
\bA \bx \cong \by
\ee
or equivalently (if a matrix $\bA$ is not symmetric positive-definite)
\be
\bA^T \bA \bx \cong \bA^T \by
\ee
where $\bA \in \Real^{I \times J}$, (with $I \geq J$),  $\by \in \Real^I$
 and a matrix $\bA^T \bA \in \Real^{J \times J}$  is  a symmetric positive-definite matrix which
 does not need to be explicitly computed (see Fig. \ref{Fig:ATAX=Y}).
 The objective is to find the  vector $\bx \in \Real^J$ in a TT
 format.

To solve this problem in the Least Squares (LS) sense, we  minimize the following cost function
\be
J(\bx) &=&  ||\bA \bx -\by||^2_2= (\bA\bx-\by)^T(\bA\bx-\by) \nonumber \\
&=& \bx^T \bA^T \bA \bx - 2 \bx^T \bA^T \by +\by^T \by,
\ee
which can be simplified to
\be
J(\bx) = \bx^T \bA^T \bA \bx - 2 \bx^T \bA^T \by.
\ee
Using the TT representation of a matrix $\bA$ and vectors $\bx$ and $\by$
\cite{Dolgov2013alternating,Dolgov2013alternating2,Oseldets-Dolgov-lin-syst12}, we have:
\be
\underline \bA &=& \llbracket \underline\bA^{(1)}, \underline\bA^{(2)}, \ldots, \underline\bA^{(N)} \rrbracket \notag \\
\underline \bX &=& \llbracket \underline\bX^{(1)}, \underline\bX^{(2)}, \ldots,\underline \bX^{(N)} \rrbracket \\
\bY &=&  \llbracket \bY^{(1)}, \bY^{(2)},\ldots, \bY^{(N)} \rrbracket \notag
\ee
and  upon applying the frame equation $\bx=\bX_{\neq n} \bx^{(n)}$ with  frame matrices
\be
\bX_{\neq n} = (\bX^{<n}_{(n)})^T \otimes \bI_{I_n} \otimes (\bX^{>n}_{(1)})^T \in \Real^{J_1 J_2  \cdots  J_{N} \times R_{n-1} J_n R_n}, \notag
\label{frame-matricesX2}
\ee
the cost function can be  written as
\be
J(\bx) = J(\bX_{\neq n} \bx^{(n)}) &=&  (\bx^{(n)})^T \bX_{\neq n}^T \bA^T \bA \bX_{\neq n} \bx^{(n)} \notag \\
&&- 2 (\bx^{(n)})^T \bX_{\neq n}^T \bA^T  \by.
\ee

This converts the  problem of solving a  large-scale system of linear equations into
 to solving  smaller system of algebraic equations iteratively
\be
 \bar \bA^{(n)} \bx^{(n)} \cong \by^{(n)}, \quad n=1,2,\ldots,N,
\ee
where $\bx^{(n)} \in \Real^{R_{n-1} J_n R_n}$ and
\be
\bar \bA^{(n)} &=&\bX_{\neq n}^T \bA^T \bA \bX_{\neq n} \in \Real^{R_{n-1}  J_n  R_n \times R_{n-1}  J_n  R_n }, \notag \\
\by^{(n)} &=& \bX_{\neq n}^T \bA^T \by   \in \Real^{R_{n-1} J_n R_{n}},
\ee
under condition that  cores are suitably left and right orthonormalized.

Of course, we cannot perform such matrix multiplications  explicitly, but in TT formats, i.e., via iterative contraction of cores in the tensor network  shown in Fig. \ref{Fig:ATAX=Y}.


The computations of a huge full vector $\bx$ or $\bA\bx$ or a matrix  $\bA^T \bA$ are not possible due to their extremely large sizes. Via tensorization,  by representing them in TT/QTT formats, and iterative contractions of cores, we can avoid the curse of dimensionality.

An assumption that data admits low-rank TT/QTT approximation is a key factor in this approach.  However,  for data with  weak structure the TT rank could be still large, which makes the calculation difficult or even impossible. The way how TT ranks are chosen and adapted during the algorithm is  very important and various  approaches to solve large structured linear systems have been proposed i
\cite{Grasedyck-rev,Dolgov2013alternating,Dolgov2013alternating2,Dolgov-Oseledets-LSyst11,DolgovGMERES,Oseledets2011dmrg}.

\begin{figure*}[ht!]
(a)
\begin{center}
\includegraphics[width=14.4cm,height=5.0cm]{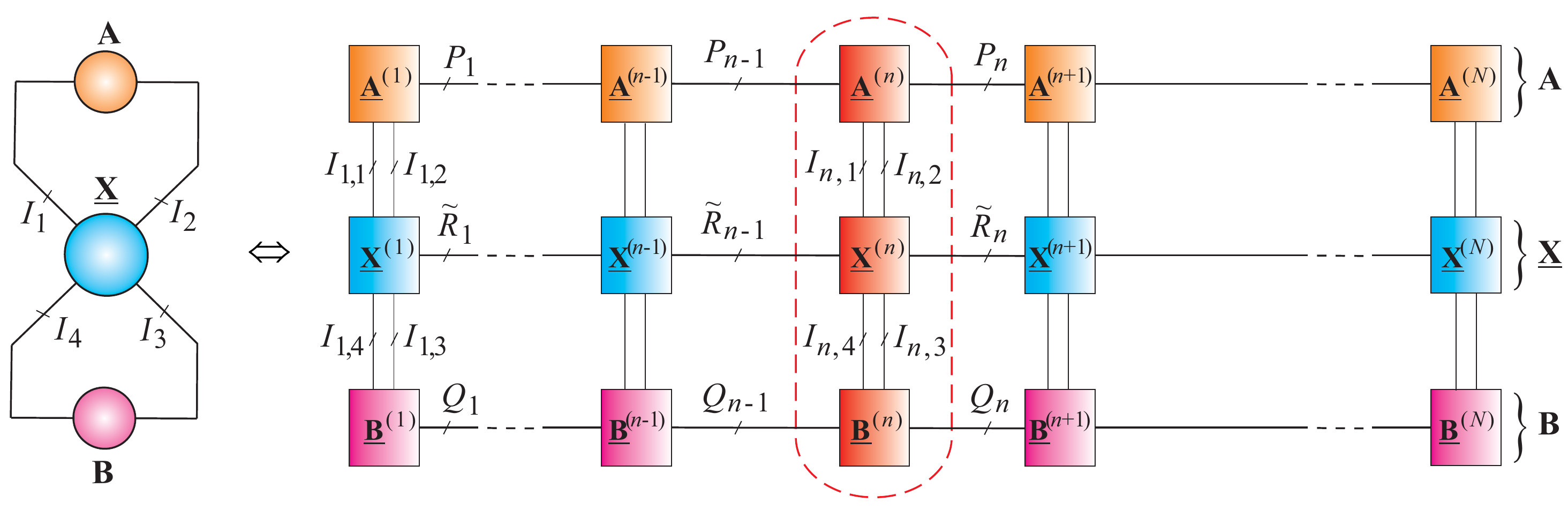}\\
\end{center}
(b)
\begin{center}
\includegraphics[width=14.4cm,height=9.6cm]{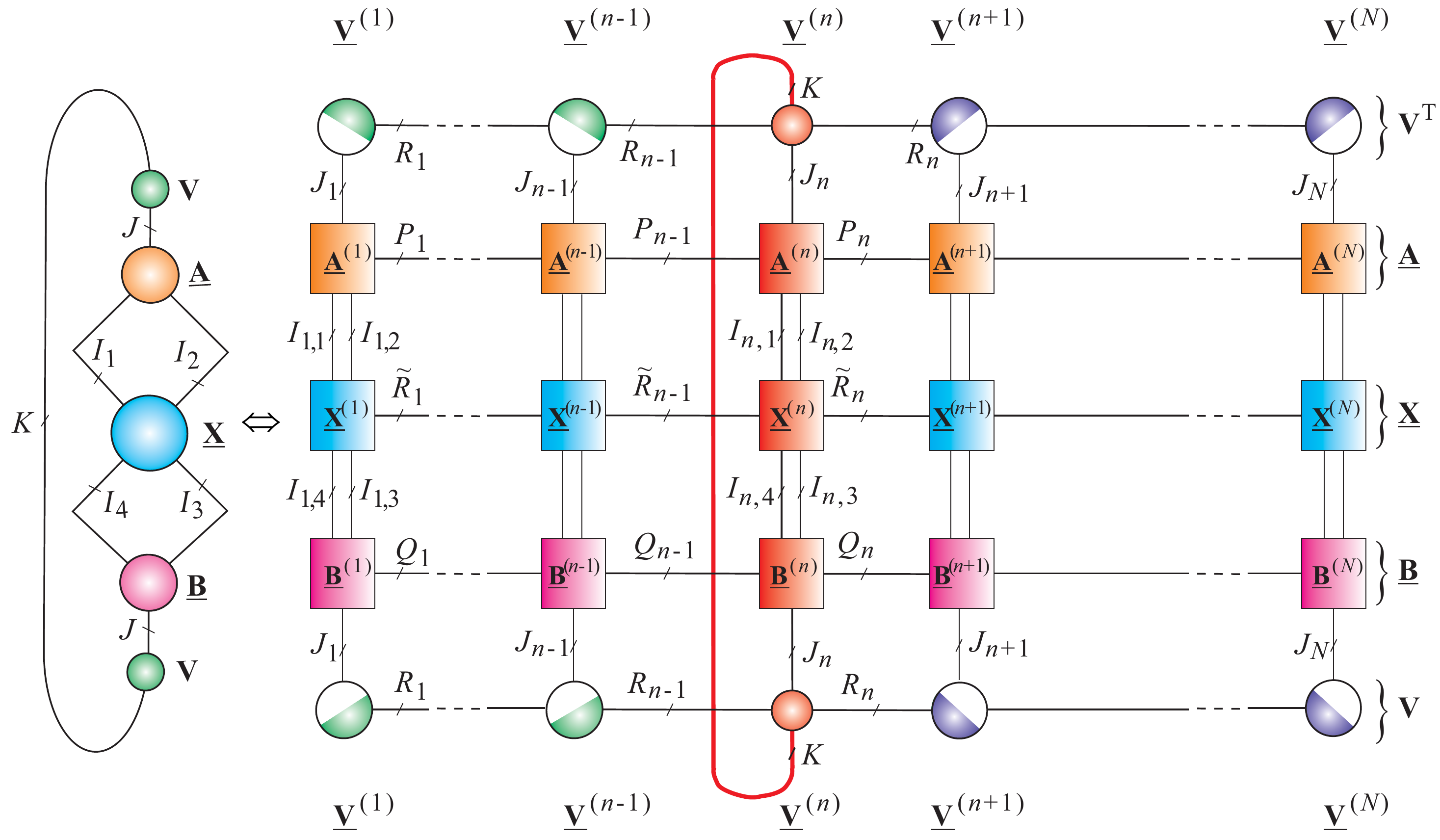}
\end{center}
\caption{Representation of tensor traces in   tensor train formats (see also Fig. \ref{Fig:trace}). These models arise in some optimization problems, in which we need to maximize the tensor traces subject to additional constraints imposed on  matrices.}
\label{Fig:TPO-traces}
\end{figure*}

{\bf Remark.} Some applications admit the use  of even more complex TT networks with higher-order cores as illustrated in Fig. \ref{Fig:TPO-traces} (a) and (b), for  which we can exploit
biorthonormality constraints \cite{Huang2011biorth}.

\subsection{\bf Software and Algorithms for Tensor Networks and Tensor Decompositions}

Tensor decompositions and tensor networks algorithms require sophisticated software libraries, which are only now being developed.

For standard TDs (CPD, Tucker models) the Tensor Toolbox for MATLAB,  originally developed by Kolda and Bader, provides several general-purpose commands and special facilities for handling sparse, dense, and structured standard TDs  \cite{tensortoolbox}, while the $N$-Way Toolbox for Matlab, by Andersson and Bro, has been  developed mostly for Chemometrics \cite{Nwaytoolbox}.
Moreover, we recently  developed the TDALAB (\url{http://bsp.brain.riken.jp/TDALAB}) and TENSORBOX  (\url{http://www.bsp.brain.riken.jp/~phan}), which provides user-friendly interface and advanced algorithms
 for basic tensor decompositions: Tucker and CPD  \cite{tdalab,tensorbox}.

The  Tensorlab toolbox developed by Sorber, Van Barel and De Lathauwer builds upon a  complex optimization framework and offers efficient numerical algorithms for computing the CPD, Block term Decomposition (BTD) or constrained Tucker decompositions.
The toolbox includes a library of many  constraints (e.g., nonnegativity, orthogonality) and offered the possibility to combine and jointly factorize dense, sparse and incomplete tensors \cite{Sorber-tensorlab}.

Similar to the CPD and/or Tucker decompositions, the TT and HT decompositions are often
based on generalized unfolding matrices $\bX_{[n]}$, and  a good approximation in a  decomposition for a given TT/HT-rank can be obtained from the SVDs of the unfolding matrices.
In practice, we avoid the explicit construction of these matrices and the SVDs when truncating a tensor via the TT decomposition to lower TT-rank. Such truncation algorithms for TT are described in \cite{OseledetsTT11}.
HT algorithms that avoid the explicit computation of
these SVDs when truncating a tensor that is already in HT decomposition are discussed in \cite{Grasedyck-rev,kressner2012htucker,espigtensorcalculus}.

In \cite{oseledets2010tt}  Oseledets
 proposed for TT decomposition a new approximative  formula in
which a $N$th-order data tensor is interpolated using  special form of Cross-Approximation, a modification of the CUR algorithm.  The total number of entries and
the complexity of the interpolation algorithm depend linearly on the order of data tensor $N$, so the developed algorithm  does not suffer from the curse of dimensionality. The TT-Cross-Approximation is
analog to the SVD/HOSVD  like algorithms for TT/MPS, but uses  adaptive
cross-approximation instead of the computationally more expensive SVD.

The TT Toolbox  developed by Oseledets (\url{http://spring.inm.ras.ru/osel/?page_id=24})  focusses on TT and QTT structures, which deal with the curse of dimensionality \cite{oseledets2012tt}.
The Hierarchical Tucker toolbox  by Kressner and Tobler  \cite{kressner2012htucker,kressner2014htucker} (\url{http://www.sam.math.ethz.ch/NLAgroup/htucker_toolbox.html}) and Tensor  library by Handschuh, Waehnert and Espig, focus mostly on HT and TT tensor networks, while TensorCalculus by Espig at al. is a C++ library is for more general tensor networks \cite{espigtensorcalculus}.

In  quantum physics and chemistry, a number of related software packages have been developed in the context
of DMRG techniques for simulating quantum networks; see for example intelligent Tensor (iTensor) by Stoudenmire and White \cite{iTensor14}. The iTensor Library is an open source C++ library for rapidly developing and applying tensor network algorithms. The iTensor is competitive with other available codes when performing  basic DMRG calculations, but due to its flexibility it is especially well suited for developing next-generation tensor network algorithms such as PEPS.

Another promising software is  the Universal Tensor Network Library  (Uni10)  developed in C++  by
Yun-Da Hsieh  and Ying-Jer Kao  (from the National Taiwan University) which
provides algorithms for  performing contraction of a complicated tensor network  with easy to use interface (\url{http://uni10.org/about.html}). The library is
 geared toward more complex tensor networks such as PEPS and MERA.


The problems related with  optimization and improvements of several
existing algorithms for TDs and TNs is an active area of research (see for example
\cite{Zhou-NTD14,Phan2012-Hess,Sorber-tensorlab}).

\section{\bf Conclusions}

Tensor networks, which can be considered as generalization and extension of tensor decompositions, are promising  tools for analysis of big data, especially, for wide family of large-scale optimization problems due to their extremely good compression abilities and distributed processing of data (cloud computing). Moreover, TNs have the  ability to address both the strong and the weak coupling between variables, and  to deal with incomplete and noisy data.
In fact, TDs have already found  application in generalized multivariate regression, multi-way blind source separation,  sparse representation and coding, feature extraction, classification, clustering and data assimilation \cite{Caiafa-Cichocki-CUR,Caiafa2012-NC,Qibin-HOPLS,NIPSQibin,PhanHALS2011,Phan2010TF,Zhou-NTD14}.

From a more general perspective,  the main concept  for big data analytic is to apply
 a suitable tensorization of the data  and  to perform an approximate decomposition in TT/QTT
  formats.  By constructing a suitable
tensor network we can perform all matrix/vectors operations in tensor network formats.
The use of the virtual tensorization or quantization (QTT) allows us to treat  more
efficiently  very large-scale data \cite{QTT-Tucker,QTT-Laplace,Kazeev2013LRT}.

In this paper, we have illuminated that  tensor networks, especially tensor trains,  are
very promising tools  for big data optimization problems, and  have illustrated
the natural and distributed  representations offered by tensor networks for a selected class of optimization problems.
This framework can be extended to a broader class of optimization problems, especially for extremely large-scale and untractable numerical problems.

In this approach a  large-scale optimization problem is transformed into a   set of small-scale linked optimization problems, each over a relative small group of unknown variables,
which are grouped via TT decompositions  and are represented by low-dimensional cores.
In other words, by representing  data in TT format we are able to turn a specific class of optimization problem into  local tractable subproblems, which have the same structure or type as the original huge optimization problem. This allows us to apply  any  efficient numerical algorithm to  local optimization problems.

The presented approach will work if only two assumptions are satisfied:
\begin{enumerate}

\item  The structured data can
be represented in TT formats that admit sufficiently good low-rank approximations.

\item Approximate solutions are acceptable \cite{Kressner-Uschmajew2014}.

\end{enumerate}

Challenging problems related to low-rank tensor approximations
remain that need to be addressed include:
\begin{itemize}

\item Current implementations of tensor train decomposition and tensor contractions still require a number of tuning parameters, e.g., approximations accuracy, TT ranks estimation. Improved
   and semi-automatic  TT approximation accuracy criteria, TT rank adaption and control and {\it a priori} errors bounds need to be developed. Particularly, the  unpredictable accumulation of rounding error and TT-rank explosion problem should be better understood and solved \cite{Savostyanov2014exact}.

\item Convergence analysis tools for TT algorithms should be developed and we need to better understand convergence properties of such algorithms.


\item  As the complexity of big data increases, this requires more efficient iterative algorithms for their computing, extending beyond the ALS/MALS, DMRG, SVD/QR and CUR/Cross-Approximation class of algorithms.

\item Theoretic and methodological approaches are needed to
determine what kind of constraints should be imposed on factor matrices/cores in order to extract desired hidden (latent) variables with meaningful physical interpretation.


  \item Generalizations of Tensor Train models to more sophisticated  tensor networks should be developed to fully integrate complex systems and optimization problems (e.g., a system simulating the biological  molecule structure) \cite{Savostyanov2014exact}.

\item Investigating the uniqueness of  various  TN models  and optimality properties, or lack thereof,
are needed and this may lead to faster and/or more reliable algorithms.

\item Special techniques  are needed to save and process huge ultra large-scale tensors which occupy peta-bytes memory.



\end{itemize}

In summary, TNs  is a fascinating and perspective area of research with  many potential applications in optimization problems for massive big data sets.



%


\end{document}